\input amstex
\documentstyle{amsppt}
\nologo
\TagsOnRight
\magnification 1200
\voffset -1cm
\input epsf

\leftline{\eightrm Tr\. Mat\. Inst\. Steklova {\eightbf 266} (2009),
149--183}
\vskip -12pt
\rightline{\eightrm Proc\. Steklov Inst\. Math\. {\eightbf 266} (2009),
142--176}
\hrule
\vskip 30pt

\def\fig#1{\smallskip\centerline{\epsffile{#1.EPS}}}
\def\on#1{\expandafter\def\csname#1\endcsname{{\operatorname{#1}}}}
\on{id} \on{Hom} \on{eq} \on{Fr} \on{Int} \on{Cl} \on{secat} \on{Tor} \on{Ext}
\def\fn#1{\expandafter\def\csname#1\endcsname{\operatorname{#1}\def\shortcut
{\if(\next{}\else\,\fi}\futurelet\next\shortcut}} \fn{coker} \fn{mod} \fn{rel}
\fn{lk} \fn{st} \fn{im}
\def\R{\Bbb R} \def\?{@!@!@!@!@!} \def\RP{\R\? P} \def\z#1#2{#1\?#2\?}
\def\Z{{\Bbb Z}\def\shortcut{\if/\next{}\z\fi}\futurelet\next\shortcut}
\let\emb\hookrightarrow \let\imm\looparrowright
\let\tl\tilde \let\but\setminus \let\x\times
\let\eps\varepsilon \let\tta\vartheta \let\phi\varphi \let\theta\vartheta
\let\Cup\smallsmile \let\Cap\smallfrown
\def\invlim{\lim\limits_{\longleftarrow}\,}
 \def\Q{\Bbb Q}
\loadbold \def\bigast{{\,\boldkey *\,}}
\def\Th{\operatorname{T}}
\let\imp\Rightarrow \def\H{\Cal H}
\def\vl{\hskip3pt {\ssize\mid}\hskip3pt }

\topmatter %%%%%%%%%%%%%%%%%%%%%%%%%%%%%%%%%%%%%%%%%%%%%%%%%%%%%%%%%%%%%%%%%%%%
\thanks Supported by Russian Foundation for Basic Research Grant No.
08-01-00663-a; Mathematics Branch of the Russian Academy of Sciences Program
``Theoretical Problems of Contemporary Mathematics''; and President Grant
MK-5411.2007.1.
\endthanks
\address Steklov Mathematical Institute, Division of Geometry and Topology;
Gubkina ~8, Moscow 119991, Russia \endaddress
\email melikhov\@mi-ras.ru \endemail

\title The van Kampen obstruction and its relatives \endtitle

\author Sergey A. Melikhov \endauthor

\abstract
We review a cochain-free treatment of the classical van Kampen obstruction
$\theta$ to embeddability of an $n$-polyhedron into $\R^{2n}$ and consider
several analogues and generalizations of $\theta$, including an extraordinary
lift of $\theta$ which in the manifold case has been studied by J.-P. Dax.
The following results are obtained.

$\bullet$ The $\bmod 2$ reduction of $\theta$ is incomplete, which answers
a question of Sarkaria.

$\bullet$ An odd-dimensional analogue of $\theta$ is a complete obstruction
to linkless embeddability (=``intrinsic unlinking'') of the given
$n$-polyhedron in $\R^{2n+1}$.

$\bullet$ A ``blown up'' $1$-parameter version of $\theta$ is a universal
type $1$ invariant of singular knots, i.e.\ knots in $\R^3$ with a finite number
of rigid transverse double points.
We use it to decide in simple homological terms when a given integer-valued
type $1$ invariant of singular knots admits an integral arrow diagram
(=\,Polyak--Viro) formula.

$\bullet$ Settling a problem of Yashchenko in the metastable range, we
obtain that every PL manifold $N$, non-embeddable in a given $\R^m$,
$m\ge\frac{3(n+1)}2$, contains a subset $X$ such that no map $N\to\R^m$ sends
$X$ and $N\but X$ to disjoint sets.

$\bullet$ We elaborate on McCrory's analysis of the Zeeman spectral sequence to
geometrically characterize ``$k$-co-connected and locally $k$-co-connected''
polyhedra, which we embed in $\R^{2n-k}$ for $k<\frac{n-3}2$ extending
the Penrose--Whitehead--Zeeman theorem.

\endabstract
\endtopmatter%%%%%%%%%%%%%%%%%%%%%%%%%%%%%%%%%%%%%%%%%%%%%%%%%%%%%%%%%%%%%%%%
\vskip-15pt
\centerline{\bf Contents}
{\eightrm
1. Introduction

2. A glimpse of underground algebraic topology

3. The van Kampen obstruction

4. Linkless and panelled embeddings

5. Arrow diagram formulas for type 1 invariants of singular knots

6. Extraordinary van Kampen obstruction

7. Embeddability versus disjoinability

8. Polyhedra whose subsets satisfy partial Alexander duality

}

\medskip
\def\epsfsize#1#2{0.15\hsize}
\fig{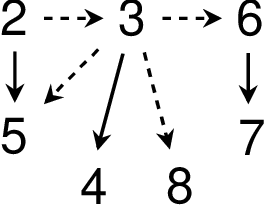}

Dependence of sections.
Dashed arrows stand for `consult as the need arises'.

\document
\head 1. Introduction \endhead

This paper is intended to be the first in a series devoted to extracting
obstructions to embeddability and (especially) knot and link invariants from
cohomological invariants of maps between configuration spaces and their various
compactifications.
In the present paper we shall consciously limit ourselves to a few very basic
situations, and merely try to put them ``in order''.

After some algebraic preparation in \S2 (which is mostly needed for
the purposes of \S5 and can be skipped on the first reading) we start in \S3
by reviewing some known facts about the van Kampen obstruction $\theta(X)$
in an invariant language.
The only treatment of $\theta(X)$ of this kind found by the author in
the literature is in \v Svarc's classical paper \cite{Sch} on the sectional
category of a fibration.
\footnote{Beware that it occurs on page 69 of the 84-page paper among other
applications.
The author is grateful to M. Grant for calling his attention to this reference,
which was overlooked in the initial version of the present paper (posted on
the arXiv in 2006).}
A clear treatment of the $\bmod 2$ reduction of $\theta(X)$ in invariant terms
has been given by Bestvina, Kapovich and Kleiner \cite{BKK} (see also
a refinement in \cite{SV}).

As a byproduct, we reidentify the cohomology group containing $\theta(X)$:
it was correctly identified in the classical papers of Shapiro and Wu, but
much of the modern literature, starting with \cite{FKT} and including, for
instance, \cite{Va2}, replicates a sign error, which instead of the order $2$
element $\theta(K_5)$, say, leads to a non-invariant element of infinite order.

Thus in \S3 we recall:

\roster
\item"$\bullet$" several definitions of the van Kampen obstruction
$\theta(X)$, which obstructs embeddability of the $n$-polyhedron $X$ in
$\R^{2n}$;

\item"$\bullet$" its $1$-parameter version $\theta(f,g)$, which is
an obstruction to isotopy of the two embeddings $f,g\:X\emb\R^{2n+1}$;

\item"$\bullet$" why $\theta(K^n)$ is nonzero, where $K^n$ is either
the $n$-skeleton of the $(2n+2)$-simplex or the join of $n+1$ copies
of the $3$-point set;

\item"$\bullet$" a proof that $\theta(X)$ is complete when $n>2$.
\endroster

The remaining content of the paper is organized as follows.

\definition{Van Kampen obstruction}
K. Sarkaria noticed that if $X$ is a graph and the $\bmod 2$ reduction of
$\theta(X)$ is zero, then $\theta(X)$ is zero, and he asked whether this
is true of every $n$-polyhedron $X$ \cite{Sa2} (see also \cite{Sk2} concerning
an error in \cite{Sa2}).
In Example 3.6 we show that this is not so for each $n>1$.

In Theorem 4.4 we find that for each $n$ (including $n=2$), $\theta(X)$
is a complete obstruction to the existence of an embedding $X\emb\R^{2n+1}$
extending to a map of the cone $f\:CX\to\R^{2n+1}$ such that $f^{-1}(X)=X$.

We also give a recipe for computation of $H^{2n}(\bar X)$, where
$\bar X=(X\x X\but\Delta)/(\Z/2)$, the group that contains $\theta(X)$
(Theorem 3.8).
In particular, this provides an easily verifiable sufficient condition for
the vanishing of $\theta(X)$.
\enddefinition

\definition{Linkless embeddings}
An odd-dimensional analogue $\eta(X)$ of $\theta(X)$ turns out to be
a complete obstruction to linkless embeddability of the $n$-polyhedron $X$
into $\R^{2n+1}$ (Theorem 4.2).%
\footnote{Added in v5: This is not quite true of $\tta(X)$ as defined in 
the published version of this paper, but true of its modified version $\tta'(X)$.}
We call an embedding $g$ of an $n$-polyhedron $X$ in $\R^{2n+1}$
{\it linkless} if every two disjoint closed subpolyhedra of $g(X)$ can be
separated by an embedded $(m-1)$-sphere; variants of this definition are
discussed in \S4.
(Nonexistence of linkless embeddings is also known as ``intrinsic linking'' in
the literature.)

In fact, the principal initial motivation of this paper, coming from \cite{MS},
was to find a method of ``local'' computation of $\theta(X)$, which would take
into account such information as $\theta(L)$ and $\eta(L)$, where $L$ runs over
the links of all simplices of some triangulation of $X$.
While this is not quite achieved at this point, some simple connections
between embeddability and local linking phenomena are given by Lemma 3.7,
which is folklore, and by Proposition 4.6.
The latter is used to give yet another proof that the $n$-skeleton of the
$(2n+3)$-simplex admits no linkless embeddings in $\R^{2n+1}$ (Example 4.7)
--- based on the fact that the $(n+1)$-skeleton of the $(2n+4)$-simplex does not
embed in $\R^{2n+2}$ even ``modulo 2''.
\enddefinition

\definition{Type 1 invariants of singular knots in $\R^3$}
In the theory of finite type invariants of knots, the study of type $m$
invariants of knots (up to smooth isotopy) is often broken down in steps
by integrating the successive Vassiliev derivatives of the invariant
(see, in particular, the ``actuality tables'' in \cite{Va1}).
At the initial step, one looks at the $(m-1)$th Vassiliev derivatives of type
$m$ invariants of knots.
These can be viewed as type 1 invariants of singular knots, that is, of smooth
immersions $S^1\looparrowright\R^3$ with prescribed double points, considered
up to smooth regular homotopy with prescribed double points.
(More detailed definitions are given in \S5.)
Notice that the smoothness of the regular homotopy guarantees that the
germs of the four strands going out of each double point of the knot lie
on an infinitesimal $2$-disk in $\R^3$, and in particular their cyclic
order on this disk does not change under the regular homotopy.
Thus the relation of smooth regular homotopy with prescribed double points
is strictly stronger than that of piecewise-smooth regular homotopy with
prescribed double points.
It is easy to exhibit specific singular knots that are distinct but equivalent
up to piecewise-smooth regular homotopy.

We construct (or rather reconstruct from \cite{Va2}) a configuration space of
the singular knot, using blowups to keep track of the infinitesimal $2$-disks
at the vertices.
The resulting ``smooth'' analogue $\zeta$ of the van Kampen graph planarity
obstruction $\theta$ is found to be a complete obstruction to planarity of
a singular knot, with its cyclic order of edges around each vertex
(Corollary 5.13).
The $1$-parameter version of $\zeta$ is found to yield a universal type $1$
invariant of singular knots (Corollary 5.3).

Next we study, prompted by Vassiliev's work \cite{Va2}, the question of
existence of arrow diagram formulas (also known as Polyak--Viro formulas
or ``combinatorial formulas'') with integer coefficients.
In fact, \S5 originated from the author's extended review for Mathematical
Reviews \cite{M3}, where the results of Vassiliev's 50-page paper \cite{Va2},
with the omission of explicit calculations and of generalization to higher
dimensions, were reproved in much less space and by elementary means, without
using Vassiliev's sophisticated machinery.
In the present paper we elaborate on this development and fill in some
details.
In particular, we correct an inaccuracy in the description of the configuration
space $\bar\Theta$ given in \cite{M3} and eliminate the usage of Vassiliev's
calculations, making our approach to his results self-contained.

More specifically, among other results, we decide in simple homological terms
whether a given type $1$ invariant of singular knots has an integral arrow
diagram formula (Theorem 5.6) and express Vassiliev's partial obstruction
$\propto(v)$ to existence of an integral arrow diagram formula in terms of
$\zeta$ (Theorem 5.10).
\enddefinition

\definition{Embeddings in the metastable range via generalized cohomology}
A ho\-motopy-theoretic classification of embeddings (up to isotopy) of
a compact polyhedron $X$ in $\R^m$ in the metastable range is provided by
the well-known Haefliger--Weber Criterion 3.1.
In \S6, we combine it with some basic equivariant algebraic topology to obtain
a more concrete classification in terms of generalized cohomology.
Namely, the set of isotopy classes of embeddings $X\emb\R^m$, if nonempty,
is in a natural bijection with an equivariant%
\footnote{Equivariant generalized cohomology reduces to non-equivariant one
by a standard procedure, which we will recall.}
stable cohomotopy group of $\tl X_+:=(X\x X\but\Delta)\sqcup\,$basepoint
(Theorem 6.6).
We also define a complete obstruction $\Theta^m(X)$ in generalized cohomology
to embeddability of $X$ in $\R^m$ in the metastable range (Theorem 6.3).
It is an element of an equivariant stable cohomotopy group of $\tl X_+$, and
when $m=2n$, it is identified with $\theta(X)$ by the Hurewicz isomorphism.

On the day of the final deadline for submission of the final version of this
paper to the special issue of the journal, the author discovered a paper by
J. P. Dax \cite{Dax}, where essentially the same results as those of \S6 have
been obtained in the case of embeddings of smooth manifolds.
The restriction to manifolds is essential for Dax's arguments, but modulo
\cite{We}, \cite{Har} and \cite{BRS}, they could have been easily modified
(in fact, simplified) so as to apply to polyhedra.
Regrettably Dax's work was not as influential as Hatcher--Quinn's
\cite{HQ; Theorem 2.3} (compare \cite{KW; Corollary G} and \cite{Kl}), where
a closely related problem is approached with a slightly different philosophy,
leading to a rather frighteningly looking answer.
For certain manifolds, a geometric description of Dax's group structure on
the set of isotopy classes of embeddings in $\R^m$ was found by A. Skopenkov
\cite{S2}.

A controlled version $\Theta(f)$ of $\Theta(X)$ is a complete obstruction to
$C^0$-approxim\-ability by embeddings of the map $f\:X\to\R^m$ in
the metastable range (Criterion 7.1).
$\Theta(f)$ generalizes both Skopenkov's cohomology obstruction \cite{RS} and
Akhmetiev's obstruction in skew-framed bordism \cite{Ah}.
\enddefinition

\definition{Towards understanding the notion of embeddability}
An intriguing connection between embeddability and what apparently is
some kind of local linking is given by Yashchenko's PhD thesis along with
our affirmative solution of the metastable case of his problem.
\enddefinition

\proclaim{Yashchenko's Theorem 1.1 {\rm (see \cite{Ya})}}
If a compactum $N$ does not admit a map to $\R^m$ with countable singular
set, there exists a subset $X\i N$ such that no map $f\:N\to\R^m$ sends $X$
and $N\but X$ to disjoint sets.
\endproclaim

By the {\it singular set} of a (continuous) map $f\:N\to\R^m$ we mean
$S_f=\{x\in N\mid f^{-1}(f(x))\ne\{x\}\}$.

\proclaim{Yashchenko's Problem 1.2 \cite{Ya}}
If a smooth manifold $N$ admits a map to $\R^m$ with countable singular
set, does it embed in $\R^m$?
\endproclaim

As Yashchenko notes in \cite{Ya}, the answer would be negative for any
compact $n$-polyhedron, non-embeddable in $\R^{2n}$.

Maps with countable singular sets abound.
Whenever $n<m$, one can easily construct a sequence of embeddings
$g_n\:S^n\to\R^m$ that uniformly converge to a map $f$ with countable dense $S_f$
(compare \cite{RR+}).
It is equally easy to construct a sequence of immersions $g_n\:S^1\imm\R^2$
(with winding numbers $w(g_n)=2^n$, say) that uniformly converge to a map
$f\:S^1\to\R^2$ whose restriction to any open subset is non-approximable
by embeddings, and whose $S_f$ is the set of all dyadic rational points
of $S^1=\R/\Z$.
More generally, if two $(n-1)$-spheres link nontrivially in $S^{m-1}$,
any embedding of an $n$-manifold $N$ into $\R^m$ can be modified by
connect-summing with increasingly small copies of the suspension over
the link into a map $f\:N\to\R^m$ whose restriction to any open subset is
non-approximable by embeddings, and whose $S_f$ is a countable dense set.

Using the controlled extraordinary van Kampen obstruction, we show that a map of
a PL $n$-manifold in $\R^m$, $m\ge\frac{3(n+1)}2$, with a countable (possibly
dense) singular set can be arbitrarily closely approximated by maps with finite
singular sets (Theorem 7.2).
In Theorem 7.4 we eliminate the finite singular set by a standard technique.
These results combine with Yashchenko's theorem to yield

\proclaim{Corollary 1.3} If a PL $n$-manifold $N$ does not embed in some $\R^m$,
$m\ge\frac{3(n+1)}2$, then $N$ contains a subset $X$ such that no map
$N\to\R^m$ sends $X$ and $N\but X$ to disjoint sets.
\endproclaim

The trouble with this result is that it does not tell us the answer to

\proclaim{Problem 1.4} Can this $X$ be chosen to be a subpolyhedron?
\endproclaim

\definition{Polyhedra satisfying partial Alexander duality}
In the last section we generalize the Penrose--Whitehead--Zeeman theorem
that $k$-connected PL $n$-manifolds embed in $\R^{2n-k}$ to $n$-polyhedra that
might be termed `homologically $k$-co-connected and locally homologically
$k$-co-connected'.
We find that such polyhedra coincide with those whose subsets satisfy
the first $k$ Alexander duality homomorphisms, and also (up to a shift by
one dimension) with those homologically $k$-connected polyhedra whose
$i$-cycles enjoy some kind of ``homological transversality'' with respect
to subpolyhedra of codimension $<i$, for all $i<k$ (Theorem 8.3).
This provides several related answers to a question of V. M. Buchstaber, who
queried the author regarding generalizations of the Penrose--Whitehead--Zeeman
theorem to polyhedra.
\enddefinition

\subhead Acknowledgements \endsubhead
I am grateful to V. M. Buchstaber and E. Kudryavtseva for stimulating questions;
to P. Akhmetiev, E. V. Shchepin, A. Skopenkov, M. Skopenkov, S. Tarasov and 
V. A. Vassiliev
for useful discussions; to Yu\. B. Rudyak and R. Sadykov for referring me to
\cite{M+} and to M. Grant for referring me to the literature on sectional
category.
P. Akhmetiev helped me to simplify the proof of Theorem 7.2 and to understand
how the extraordinary van Kampen obstruction fits with the viewpoint of
geometric immersion theory (in particular, Example 6.5 is due to him).
V. Manturov's continued interest in hearing the content of \S5 at his seminar
(despite a series of various organizational obstacles that delayed this for
almost two years) was crucial for convergence of preparation of the paper.
I am deeply grateful to M. Bestvina, J. Dydak and D. Sinha for inviting me
to present parts of this work at their respective Universities in early
Spring 2007.%
\footnote{Added in v5: I am also grateful to R. Nikkuni and V. Turchin 
for useful feedback.}

\head 2. A glimpse of underground algebraic topology \endhead

In this section we recall the notions of homological and cohomological
transfers of a double covering and of the Euler class of the associated
line bundle, and use them to give a geometric description of homomorphisms in
the Smith sequences of a free action of $\Z/2$.
A convenient language to discuss the transfers and the Euler class (which
language in particular has a name for the natural generality where they are
defined) is the geometric understanding of cohomology from the unique book of
Buoncristiano, Rourke and Sanderson \cite{BRS}.
The reader is hereby warned that mainstream algebraic topology (as exemplified
by J. F. Adams' review of \cite{BRS}) considers the ideas contained in this book
to be dangerous (for one's ``balance between the particular and the general''), and
is openly uncomfortable with the fact that its publication could not have been
prevented.

Local coefficient systems are discussed in the textbooks by Hilton--Wylie,
Spanier and Hatcher; and in more detail in those by G. Whitehead and
Davis--Kirk, not to mention sheaf theory texts such as Bredon's \cite{Bre}.

\definition{Pseudo-manifolds, local coefficients}
By an {\it $m$-pseudo-manifold} we shall understand a finite simplicial complex
where all simplices have dimension $\le m$ and every $(m-1)$-simplex is a facet
of precisely two $m$-simplices.
An {\it $m$-pseudo-manifold with boundary} is a pair of finite simplicial
complexes $(M,\partial M)$ where $\partial M$ is an $(m-1)$-pseudo-manifold and
the double $M\x\{0\}\cup_{\partial M\x\{0\}=\partial M\x\{1\}} M\x\{1\}$ is an
$m$-pseudo-manifold.

An {\it orientation} of an $n$-simplex, $n>0$ is an orbit of the alternating
group $A_{n+1}$ acting on the set of numberings of its vertices by $1,\dots,n+1$.
By convention, the $0$-simplex has two orientations $+$ and $-$.
An oriented $n$-simplex $\sigma$ that is a facet of an oriented $(n+1)$-simplex
$\tau$ is its {\it oriented facet} if either $n>0$ and the first $n+1$ vertices
of some representative numbering of the vertices of $\tau$ constitute
a representative numbering of the vertices of $\sigma$, or $n=0$ and
the orientation of $\sigma$ equals the sign of the unique permutation converting
a representative numbering of the vertices of $\tau$ into one where $\sigma$ is
numbered $1$.
An {\it orientation} of an $m$-pseudo-manifold with boundary is a choice of
orientations of its $m$- and $(m-1)$-simplices such that each oriented
$(m-1)$-simplex is an oriented facet of precisely one oriented $m$-simplex.

A {\it local coefficient system} on a polyhedron $P$ is a collection of
abelian groups $G_x$ indexed by $x\in P$ along with isomorphisms
$h_l\:G_{l(0)}\to G_{l(1)}$ for every path $l\:I\to P$ such that each $h_l$
depends only on the homotopy class of $l$ relative to its endpoints.
(Here $I=[0,1]$.)
If $p\:Q\to P$ is a double covering, it determines an integral (i.e.\ with
each $G_x\simeq\Z$) local coefficient system $\Z_p$ on $P$, called the
{\it orientation sheaf} of $p$, that associates the nontrivial automorphism
of $\Z$ precisely to those loops that do not lift to $Q$.
Clearly, every integral local coefficient system is of this type.

A simplicial map $f$ of an $m$-pseudo-manifold $M$ into a simplicial complex
$K$ is called a {\it singular} $m$-pseudo-manifold in $K$.
If $\Cal O$ is an integral local coefficient system on the underlying polyhedron
$P$ of $K$, an {\it $\Cal O$-orientation} of $f$ is a choice of orientations of
$m$- and $(m-1)$-simplices of $M$ such that a loop $l$ in the dual $1$-skeleton
of $M$ intersects an odd number of ``disoriented'' $(m-1)$-simplices if and
only if $f(l)$ acts nontrivially on $\Z$.
(An $(m-1)$-simplex is considered ``disoriented'' if it fails to be an oriented
facet of precisely one oriented $m$-simplex.)

We define $H_i(P;\Z/2)$ (resp\. $H_i(P;\Cal O)$), where $P$ is a polyhedron
(and $\Cal O$ an integral local coefficient system on $P$), to be the group of
($\Cal O$-oriented) singular pseudo-bordism classes of ($\Cal O$-oriented)
singular $i$-pseudo-manifolds in some fixed triangulation of $P$.
It is well-known, and easy to see that this is equivalent to the usual
definition (compare \cite{BRS}, \cite{Ah}), so in particular does not depend
on the choice of a triangulation of $P$.
An $m$-pseudo-manifold $M$ is $\Cal O$-orientable iff $H^m(M;\Cal O)$ contains
no torsion; in which case its $\Cal O$-orientation is equivalent to a choice of
a set of generators of $H_m(M;\Cal O)$ representable by cycles with
disjoint supports.%
\footnote{Added in v5: Here the ``cycles'' are meant to be simplicial, and by the support of 
an $m$-cycle we mean a subset of the set of $m$-simplices of $M$, rather than a subset of $M$ itself.}
\enddefinition

\definition{Pseudo-comanifolds}
We call a PL map $f\:W\to P$ between polyhedra a {\it $k$-pseudo-comanifold}
over a triangulation $K$ of $P$ if the preimage $f^{-1}(\sigma)$ of every
$i$-simplex $\sigma$ of $K$ is triangulated by an orientable
$(i-k)$-pseudo-manifold $M$ with boundary $\partial M$ such that
$f^{-1}(\partial\sigma)$ coincides with the underlying polyhedron of
$\partial M$.
For instance, if $T$ is the three-page book $(\{0\}*\{1,2,3\})\x\R$, then every
embedded $1$-pseudo-comanifold $G$ in $T$ is a graph with vertices of degrees
$1$ and $3$, where the degree $3$ vertices coincide with the intersection
points of $G$ and the binding $\{0\}\x\R$, and the degree $1$ vertices
coincide with the intersection points of $G$ and the page edges $\{1,2,3\}\x\R$.
Thus embedded $1$-pseudo-comanifolds in $T$ can be alternatively described
as transversal point-inverses of maps $T\to S^1$ (using that they are always
co-oriented, see below).
Pseudo-comanifolds originate in \cite{BRS}, where they are called
``mock bundles with codimension two singularities''.
See also \cite{Fe}, \cite{FS}, \cite{CG}, \cite{Gra}, \cite{De}.
Using the Poincar\'e duality, pseudo-comanifolds over a triangulation $K$ of
a polyhedron $P$ can also be viewed as pullbacks over $P$ of PL maps of
pseudo-manifolds with boundary $(M,\partial M)\to(N,\partial N)$ that are
transverse to $K$, where $N$ is a PL manifold neighborhood of a copy of $P$
in some Euclidean space.

If $f\:W\to P$ is a $k$-pseudo-comanifold over a triangulation $K$ of $P$,
a {\it co-orientation} of $f$ over a simplex $\sigma$ of $K$ is a choice of
bijection $b_{\sigma}$ between the two orientations of $\sigma$ and the two
orientations of%
\footnote{Added in v5: each connected component of}
$f^{-1}(\sigma)$.
If $(\sigma,o_\sigma)$ is an oriented facet of $(\tau,o_\tau)$, a co-orientation
$b_\sigma$ of $f$ over $\sigma$ is said to {\it cohere} with a co-orientation
$b_\tau$ of $f$ over $\tau$ if $(f^{-1}(\sigma),b_\sigma(o_\sigma))$ is an
oriented submanifold of the oriented boundary of $(f^{-1}(\tau),b_\tau(o_\tau))$.
A {\it co-orientation} of $f$ over $K$ is a coherent choice of its co-orientations
over all simplices of $K$.
If $\Cal O$ is an integral local coefficient system on $P$,
an {\it $\Cal O$-co-orientation} of $f$ over $K$ is choice of its co-orientations
$b_\sigma$ over all simplices $\sigma$ of $K$ such that for every sequence of
simplices $\sigma=\sigma_1,\dots,\sigma_k=\sigma$ of $K$ where for each $i$,
either $\sigma_i$ is a facet of $\sigma_{i+1}$ or vice versa, the number of
incoherences between consecutive $b_{\sigma_i}$ is odd if and only if
the corresponding loop in the $1$-skeleton of the barycentric subdivision of
$K$ acts nontrivially on $\Z$.

We define $H^k(P;\,\Z/2)$ (resp\. $H^k(P;\Cal O)$), where $P$ is a polyhedron
(and $\Cal O$ a local coefficient system on $P$), to be the group of
($\Cal O$-co-oriented) pseudo-cobordism classes of ($\Cal O$-co-oriented)
$k$-pseudo-comanifolds over some fixed triangulation of $P$.
It is easy to see that this is equivalent to the usual definition
(cf\. \cite{BRS}) so in particular does not depend on the choice of
a triangulation of $P$.
If $W\i P$ is an embedded $k$-pseudo-comanifold, it is $\Cal O$-co-orientable
iff $H_k(P,P\but W;\,\Cal O)$ contains no torsion; in which case its
$\Cal O$-co-orientation is equivalent to a choice of a set of generators of
$H^k(P,P\but W;\,\Cal O)$ representable by cocycles with disjoint supports.
Note that if $f\:W\to P$ is a pseudo-comanifold, where $P$ is finite-dimensional,
there exists an embedded pseudo-comanifold $\bar f\:W\to P\x\R^N$ for some $N$
that projects onto $f$.
\enddefinition

\definition{Transfers, products, Euler class}
Let $\phi\:P\to Q$ be a PL map between polyhedra and let $\Cal O$ be a local
coefficient system on $Q$.
Suppose first that $\phi$ is triangulated by a simplicial map $K\to L$.
If $f$ is a singular $i$-pseudo-manifold in $K$, then $\phi f$ is a singular
$i$-pseudo-manifold in $L$, so $[f]\mapsto [\phi f]$ defines the induced
homomorphism $\phi_*\:H_i(P;\phi^*\Cal O)\to H_i(Q;\Cal O)$.
On the other hand, if $f\:W\to Q$ is an $\Cal O$-co-oriented $i$-pseudo-comanifold
over $L$, the pullback $\phi^*(f)\:V\to P$ is a $\phi^*\Cal O$-co-oriented
$i$-co-pseudo-manifold over $K$, cf\. \cite{BRS; bottom of p.\ 23}.
This defines the induced homomorphism
$\phi^*\:H^i(Q;\Cal O)\to H^i(P;\,\phi^*\Cal O)$.
$$\CD
V@>\phi^*(f)>>P\\
@VVV@V\phi VV\\
W@>f>>Q
\endCD$$
Now suppose that $\phi$ is an $\Cal F$-co-oriented $k$-pseudo-comanifold
over a triangulation $L^!$ of $Q$, where $\Cal F$ is a local coefficient system
on $Q$, and the $\phi$-preimages of simplices of $L^!$ are subcomplexes of
a triangulation $K^!$ of $P$.
If $f$ is an $\phi^*\Cal O$-co-oriented $i$-pseudo-comanifold over $K^!$, then
$\phi f$ is an $\Cal O$-co-oriented $(i+k)$-pseudo-comanifold over $L^!$, so
$[f]\mapsto [\phi f]$ defines the {\it transfer}
$\phi_!\:H^i(P;\,\phi^*\Cal O)\to H^{i+k}(Q;\,\Cal O\otimes\Cal F)$.
On the other hand, if $f\:M\to L^!$ is an $\Cal O\otimes\Cal F$-oriented
$i$-pseudo-manifold in $L^!$, then the pullback $\phi^*(f)\:N\to K^!$ will be
a singular $\phi^*\Cal O$-oriented $(i-k)$-pseudo-manifold in $K^!$,
cf\. \cite{BRS; II.1.2}.
This defines the {\it transfer}
$\phi^!\:H_i(Q;\,\Cal O\otimes\Cal F)\to H_{i-k}(P;\,\phi^*\Cal O)$.
Using the techniques of \cite{BRS}, it is easy to show that the latter
transfer is well-defined; for the former transfer and the induced maps
this has been done in \cite{BRS}.
Note that passing between $L$ and $L^!$ generally requires making
perturbations, by the virtue of transversality (see \cite{BRS}).

For the cup- and cap-products, we have $\alpha\Cup[\phi]=\phi_!\phi^*(\alpha)$
\cite{BRS; II.2.2}; and
$\beta\Cap[\phi]=\phi_*\phi^!(\beta)$ \cite{BRS; p.\ 29}.
The {\it Euler class} $e(i)$ of an arbitrary $\Cal O$-co-oriented
$k$-pseudo-comanifold $i\:B\to E$ is defined to be
$i^*i_!([\id_B])\in H^k(B;\Cal O)$, where the co-oriented $0$-comanifold
$\id_B$ represents a generator of $H^0(B)$.
If $i$ is the zero cross-section of a $k$-vector bundle $\xi\:E\to B$, then
$e(i)\in H^k(B;\Cal O_\xi)$, where $O_\xi$ is the orientation sheaf of $\xi$,
and we call $e(i)$ the {\it Euler class} of $\xi$ and denote it by $e(\xi)$.
When $\xi$ is orientable, this is its usual Euler class; in general, its $\bmod2$
reduction is the Stiefel--Whitney class $w_k$ (see \cite{BRS; p.\ 26}).
By definition, $e(\xi)$ can be represented by a pseudo-comanifold that is
the zero set $B\cap s(B)$ of a generic cross-section $s(B)$ of $\xi$.
\enddefinition

\definition{Euler class of a line bundle}
Let $p\:E\to B$ be the double covering and $\lambda_p$ the associated line bundle.
The Euler class $e(\lambda_p)\in H^1(B;\,\Z_p)$ is always an element of order
two, since it is induced from $e(\lambda_q)\in H^1(\RP^\infty;\,\Z_q)\simeq\Z/2$
by a classifying map of $\lambda_p$, which also induces $p$ from the double
covering $q\:S^\infty\to\RP^\infty$.
Nevertheless, $e(\lambda_p)$ should not be confused with its $\bmod 2$ reduction,
the first Stiefel--Whitney class $w_1(\lambda_p)\in H^1(B;\Z/2)$.
The latter is less informative in general.
For instance, $w_1(\zeta)^2=0\ne e(\zeta)^2$ for the double covering
$\zeta\:\RP^3\to L(4,1)$.
\enddefinition

\remark{Remarks} (i) If $\lambda$ is a line bundle over a compact polyhedron $P$
and $w_1(\lambda)=0$, then $e(\lambda)=0$.
For otherwise $e(\lambda)$ would be divisible by $2$ and so would be in
the Bockstein image of $H^0(P;\,\Z_\lambda/4)=0$.
This explains how reduction of coefficients $\bmod 2$ ``does not lead to a loss
of information'' in \cite{M2; remark after Lemma 6}.

(ii) If $\lambda$ is a line bundle over a $d$-manifold $M$ and
$w_1^d(\lambda)=0$, then $e^d(\lambda)=0$.
Indeed, $H^d(M;\Z_\lambda^{\otimes d})$ is a direct sum of copies of $\Z$ and
$\Z/2$ (one summand for each compact component of $M$) and so does not contain
elements of order $4$.
This fills a minor gap in \cite{M2; proof of Theorem 5}.
\endremark

\proclaim{Lemma 2.1} If $i\:B\emb E$ is an embedded $\Cal O$-co-oriented
$k$-pseudo-comanifold, then $i^*i_!(\alpha)=\alpha\Cup e(i)$ and
$i^!i_*(\beta)=\beta\Cap e(i)$.
\endproclaim

The author is grateful to E. Kudryavtseva for calling his attention to
the problem of computation of the ``reverse transfer compositions'' $i^*i_!$
and $i^!i_*$.
When $i$ is an immersion, Lemma 2.1 can be combined with Herbert's formula
(see \cite{EG}).

\demo{Proof} By considering a regular neighborhood $R$ of $i(B)$ in $E$ we may
assume that $i^*$ is an isomorphism and so $\alpha=i^*(\gamma)$
for some $\gamma$.
Then we have $$i^*i_!(\alpha)=i^*i_!(i^*\gamma\Cup[\id_B])=
i^*(\gamma\Cup i_![\id_B])=i^*\gamma\Cup i^*i_![\id_B]=\alpha\Cup e(i)$$
using the formula $p_!(p^*\xi\Cup\zeta)=\xi\Cup p_!\zeta$ from \cite{BRS; II.2.6}.
The two homological analogues of this formula read:
$$p_*(p^!\xi\Cap\zeta)=\xi\Cap p_!\zeta\tag{$*$}$$ and
$p_*(\zeta\Cap p^*\xi)=p_*\zeta\Cap\xi$, the latter being the familiar
functoriality of $\Cap$-product.
On the other hand, the functoriality of $\Cup$-product is proved
(from the definition $\alpha\Cup[\phi]=\phi_!\phi^*(\alpha)$) by a cubical
diagram, which in the case of $\Cap$-product yields
$$p^!(\xi\Cap\zeta)=p^!\xi\Cap p^*\zeta.\tag{$**$}$$
Using formulas ($*$) and ($**$), we obtain: if $\beta=i^!(\delta)$ for some
$\delta$, then $$i^!i_*(\beta)=i^!i_*(i^!\delta\Cap[\id_B])=
i^!(\delta\Cap i_![\id_B])=i^!\delta\Cap i^*i_![\id_B]=\beta\Cap e(i).$$
It remains to observe that $i^!$ may be assumed to be an isomorphism
(specifically, the Thom isomorphism) by considering the pair $(R,\Fr R)$.
\qed
\enddemo

\definition{Equivariant cohomology and Smith sequences}
Let $P$ be a polyhedron endowed with a free PL action of
$\Z/2=\left<t\mid t^2\right>$, and let $p\:P\to P/t$ be the double covering.
Let us fix a module $M$ over the group ring $\Lambda:=\Z[\Z/2]$.
Let $C_*$ be the simplicial chain complex of some equivariant triangulation
of $P$, where simplices from the same orbit are oriented coherently.
Then $C_*$ is a free $\Lambda$-module.

We recall that equivariant homology and cohomology groups are defined by
$H_i^{\Z/2}(P;M)=H_i(C_*\otimes_\Lambda M)$ and
$H^i_{\Z/2}(P;M)=H^i(\Hom_\Lambda(C_*;M))$ (compare \cite{Bro}).%
\footnote{Added in v5: See e.g.\ Hatcher's textbook for a more thorough discussion.}
Note that $C_*\otimes_\Lambda\Lambda\simeq C_*$, whereas
$\Hom_\Lambda(C_*;\Lambda)$ is isomorphic to the simplicial cochain complex
$C^*:=\Hom_\Z(C_*;\Z)$.
Hence $H_i^{\Z/2}(P;\Lambda)\simeq H_i(P;\Z)$ and
$H^i_{\Z/2}(P;\Lambda)\simeq H^i(P;\Z)$.

Let us consider the trivial $\Lambda$-module $\Z$, which is isomorphic to
$\Lambda/(t-1)$, and the augmentation ideal $I$ of $\Lambda$, which is
isomorphic to $\Lambda/(t+1)$.
By definition, $\Hom_\Lambda(C_*;\Z)$ (resp\. $\Hom_\Lambda(C_*;I)$) is
the subcomplex of $C^*$ consisting of all cochains $f$ satisfying $f(c)=f(tc)$
(resp\. $f(c)=-f(tc)$).
It follows that $H^i_{\Z/2}(P;\Z)\simeq H^i(P/t;\Z)$ and
$H^i_{\Z/2}(P;I)\simeq H^i(P/t;\Z_p)$.
Similarly, $C_*\otimes_\Lambda\Z$ (resp\. $C_*\otimes_\Lambda I$) is
the quotient complex of $C_*$ by $c=tc$ (resp\. $c=-tc$).
Hence $H_i^{\Z/2}(P;\Z)\simeq H_i(P/t;\Z)$ and
$H_i^{\Z/2}(P;I)\simeq H_i(P/t;\Z_p)$.
\footnote{More generally, every module $M$ over $\Lambda$ gives rise to a local
coefficient system $\Cal F_M$ on $P/t$ whose stalks are isomorphic to
the underlying abelian group of $M$, and the action of $g\in\pi_1(P/t)$ on
the stalks is given by the action of $G_*(g)\in\pi_1(\RP^\infty)\simeq\Z/2$,
where $G\:P/t\to\RP^\infty$ is a classifying map of $p$.
We have $\Cal F_{M\oplus N}\simeq\Cal F_M\oplus\Cal F_N$,
$\Cal F_{M\otimes N}\simeq\Cal F_M\otimes\Cal F_N$\footnotemark\
and $H_i^{\Z/2}(P;M)\simeq H_i(P/t;\Cal F_M)$,
$H^i_{\Z/2}(P;M)\simeq H^i(P/t;\Cal F_M)$.
In particular, we obtain the Vietoris--Begle isomorphism for double coverings:
$H_i(P/t;\Cal F_\Lambda)\simeq H_i(P;\Z)$,
$H^i(P/t;\Cal F_\Lambda)\simeq H^i(P;\Z)$.}
\footnotetext{Added in v5: Here $M\otimes N$ denotes tensor product over $\Z$ with 
the diagonal action of $\Z/2$.}

The coefficient sequences $0\to I@>>>\Lambda@>>>\Z\to 0$ and
$0\to\Z@>>>\Lambda@>>>I\to 0$ give rise to long exact sequences,
which can be written as follows using the above identifications of groups
and straightforward identifications of homomorphisms:
$$\gather
\dots@>>>H^n(P/t;\Z_p)@>p^*>>H^n(P)@>p_!>>H^n(P/t)@>i^*i_!>>\dots\\
\dots@>>>H^n(P/t)@>p^*>>H^n(P)@>p_!>>H^n(P/t;\Z_p)@>i^*i_!>>\dots\\
\dots@>>>H_n(P/t;\Z_p)@>p^!>>H_n(P)@>p_*>>H_n(P/t)@>i^!i_*>>\dots\\
\dots@>>>H_n(P/t)@>p^!>>H_n(P)@>p_*>>H_n(P/t;\Z_p)@>i^!i_*>>\dots
\endgather$$
(here $i$ denotes the inclusion of $P/t$ onto the zero cross-section of
the line bundle $\lambda_p$ associated with the double covering $p\:P\to P/t$;
constant integer coefficients are omitted).
These are known as the {\it Smith sequences}.
They can also be identified with the exact sequences of the relative
mapping cylinder of the double covering $p$.
By Lemma 2.1 the connecting homomorphisms in cohomology and homology can be
described respectively as $\Cup$- and $\Cap$-multiplication by the Euler class
$e(\lambda_p)$.
\enddefinition

\proclaim{Lemma 2.2} Let $M=H_n(P)$.
Then $H^1(\Z/2;\,M)=0$ (resp.\ $H_1(\Z/2;\,M)=0$) iff each element of
$H_n(P/t;\,\Z_p)$ (resp.\ $H_n(P/t)$) is either divisible or annihilated by
$e(\lambda_p)$, and each element of $H_n(P/t)$ (resp.\ $H_n(P/t;\,\Z_p)$) is
either not divisible or not annihilated by $e(\lambda_p)$.
In addition, ``$H_n$'' may be replaced by ``$H^n$'' throughout.
\endproclaim

This follows trivially from the definitions, which for convenience of
the reader we do recall (compare \cite{Bro}).

\demo{Proof} We recall that $H_*(\Z/2;\,M)=\Tor_*^\Lambda(\Z,M)$ and
$H^*(\Z/2;\,M)=\Ext^*_\Lambda(\Z,M)$ are the derived functors of
$H_0(\Z/2;\,M)=\Z\otimes_\Lambda M$ and $H^0(\Z/2;\,M)=\Hom_\Lambda(\Z,M)$
viewed as functors of $\Z$.
Thus we start with a projective resolution of $\Z$ over $\Lambda=\Z[\Z/2]$,
namely let us take the chain complex of the simplest $\Z/2$-equivariant
cell complex structure on $S^\infty$:
$$\Z@<\eps<<\Lambda@<1-t<<\Lambda@<1+t<<\Lambda@<1-t<<\Lambda@<1+t<<\dots$$
Now we drop the initial $\Z$, apply the functor, and take homology; thus
$H_*(\Z/2;\,M)$ is the homology of the chain complex
$$M@<1-t<<M@<1+t<<\dots$$
and $H^*(\Z/2;\,M)$ is the homology of the cochain complex
$$M@>1-t>>M@>1+t>>\dots$$
So $H_1(\Z/2;\,M)=\ker(1-t)/\im(1+t)$ and
$H^1(\Z/2;\,M)=\ker(1+t)/\im(1-t)$.

Now let $M=H_n(P)$, and consider $H^1(\Z/2;\,M)$; the other cases are similar.
We have $1+t=p^!p_*(x)$ with coefficients in $\Z$, and $1-t=p^!p_*$ with
coefficients in $\Z_p$.
Then $x\in\ker(1+t)\but\im(1-t)$ iff either $y:=p_*(x)\in H_n(P/t)$ is nonzero
but lies in $\ker p^!$, or $x=p^!(z)$ for some $z\in H_n(P/t;\,\Z_p)$ that
does not lie in $\im p_*$.
Hence $\ker(1+t)\ne\im(1-t)$ iff either some $y\in H_n(P/t)$ is divisible
and annihilated by $e(\lambda_p)$ or some $z\in H_n(P/t;\,\Z_p)$ is not
divisible and not annihilated by $e(\lambda_p)$. \qed
\enddemo

\remark{Remark 2.3 (added in v5)}
The generalized Smith sequence, found in \cite{Bre; \S2.19 and \S5.20},
can be quickly obtained as follows.
Note that it is a special case of the generalized Smith--Gysin sequence,
found in \cite{Bre; \S4.12.3} and used e.g.\ in \cite{M1'}.

Let $M$ be a $\Lambda$-module and $M_0$ its underlying abelian group, and
let $M'$ denote $M\otimes_\Z I$ with the diagonal action of $\Z/2$.
By a well-known lemma \cite{Bro; III.5.7}, $M\otimes_\Z\Lambda$ with the diagonal
action of $\Z/2$ is isomorphic to the induced module $M_0\otimes_\Z\Lambda$ via
$m\otimes g\mapsto g^{-1}m\otimes g$.
Then the same argument used above to construct the usual Smith sequences shows that
the short exact sequence $0\to M\to M\otimes_\Z\Lambda\to M'\to 0$
of $\Lambda$-modules gives rise to long exact sequences
$$\gather
\dots\to H^n(P;\,M_0)\to H^n(P/t;\,\Cal F_{M'})\to H^{n+1}(P/t;\,\Cal F_M)\to\dots\\
\dots\to H_n(P;\,M_0)\to H_n(P/t;\,\Cal F_{M'})\to H_{n-1}(P/t;\,\Cal F_M)\to\dots
\endgather$$
These may be convenient, but do not really contain additional information:
for $M=\Z$ and $M=I$ we get the usual Smith sequences constructed earlier
in this section, and for $M=\Lambda$ the generalized Smith sequences split into
short exact sequences, $0\to H\to H\otimes_\Z\Lambda\to H\otimes_\Z I\to 0$,
where $H=H^n(P)$ or $H_n(P)$.
On the other hand, every $\Lambda$-module is a direct sum of $G_1$, $G_2\otimes I$
and $G_3\otimes\Lambda$ for some abelian groups $G_1$, $G_2$ and $G_3$
(see \cite{CR; Theorem 74.3}).
\endremark

\head 3. The van Kampen obstruction \endhead

If $X$ is a polyhedron, the factor exchanging involution $(x,y)\mapsto(y,x)$
on $X\x X$ restricts to a free action of $\Z/2$ on the {\it deleted product}
$X\x X\but\Delta$.
We write $\tl X$ for the $\Z/2$-space $(X\x X\but\Delta,\,\Z/2)$, and $\bar X$
for the quotient $\tl X/(\Z/2)$.
The $m$-sphere $S^m\i\R^{m+1}$ endowed with the antipodal involution
$x\mapsto -x$ will be denoted $S^m_-$.

If $g\:X\emb\R^m$ is an embedding, its {\it Gauss map}
$\tl g\:\tl X\to S^{m-1}_-$ is defined by
$(x,y)\mapsto\frac{g(x)-g(y)}{||g(x)-g(y)||}$.
In other words, $\tl g$ is the restriction $\tl X\to\widetilde{\R^m}$ of
$g\x g$ followed by the obvious equivariant homotopy equivalences
$\widetilde{\R^m}\simeq\R^m_-\but\{0\}\simeq S^{m-1}_-$.
The Gauss map descends to $\bar g\:\bar X\to\RP^{m-1}$, which assigns to
an unordered pair $\{x,y\}$ of distinct points of $X$ the $1$-subspace
$\left<g(x)-g(y)\right>$ of $\R^m$.

\definition{Sectional category}
It is easy to see that the existence an equivariant map from a $\Z/2$-space $K$
into $S^{m-1}_-$ is equivalent to the existence of a cross-section of the sphere
bundle $p_m\:K\x_{\Z/2} S^{m-1}_-\to K/(\Z/2)$.
This bundle is the Whitney join of $m$ copies of the double covering
$p_1\:K\to K/(\Z/2)$, and it follows that $p_m$ has a cross-section if and
only if there exists a cover of $P$ by $m$ open sets $U_i$ such that $p_1$
has a cross-section over each $U_i$ (see \cite{CF}, \cite{Sch}, \cite{Ja} or
\cite{CL+} for details).
The largest integer $m$ that does not satisfy any of these three equivalent
properties is called the {\it sectional category} $\secat p_1$
(cf\. \cite{CL+}). \footnote{In older literature, ``sectional category'' used
to mean $1+\secat p_1$, which is the least integer satisfying any of
the three equivalent properties (cf\. \cite{Ja}).
Originally, $1+\secat p_1$ was known as the Krasnosel'sky--\v Svarc genus
of $p_1$ \cite{Sch} and $\secat p_1$ itself as Yang's B-index or Conner--Floyd
co-index of the $\Z/2$-space $K$ \cite{CF}.
Furthermore, $\secat p_1$ equals the category of a classifying map
$K\to\RP^\infty$ of $p_1$ (see \cite{CL+}).}
\enddefinition

The largest integer $m$ such that a polyhedron $X$ does not embed in $\R^m$
is therefore bounded below by the sectional category of the double covering
$\tl X\to\bar X$.
In the metastable range the two numbers are equal:

\proclaim{Haefliger--Weber Criterion 3.1 \cite{We} {\rm (see also \cite{Ad})}}
Let $X$ be an $n$-polyhedron.

(a) Given a $\Z/2$-map $\Phi\:\tl X\to S^{m-1}_-$, if $m\ge\frac{3(n+1)}2$,
there exists an embedding $g\:X\emb\R^m$; moreover, $\tl g$ can be chosen
$\Z/2$-homotopic to $\Phi$.

(b) If $m>\frac{3(n+1)}2$, two embeddings $f,g\:X\emb\R^m$ are
ambient isotopic if and only if $\tl f$ and $\tl g$ are $\Z/2$-homotopic.
\endproclaim

\demo{Proof of the case $m=2n$ in (a)}
Fix some triangulation of $X$ and embed $X^{(n-1)}$ by general position.
Then $\tl g$, which is now defined on $\widetilde{X^{(n-1)}}$, is
equivariantly homotopic to the restriction of $\Phi$ since $S^{2n-1}$ is
$(2n-2)$-connected.
Let $\sigma_1,\sigma_2,\dots$ be the $n$-simplices of $X$.
Extend the embedding to $\sigma_1$ by general position.
We want $\tl g$, which is now also defined on $(\sigma_1\x X^{n-1}\cup
X^{n-1}\x\sigma_1)\setminus\Delta$, to be equivariantly homotopic to
the restriction of $\Phi$.
(Not just to prove the moreover part, but also to prepare for the
following steps.)
This is easy to achieve: if the two maps differ on $\sigma_1\x\tau$ for
some $(n-1)$-simplex $\tau$ by a degree $l$ map, pick a small $S^n$
winding around $\tau$ with linking number $l$, and amend $\sigma_1$ by
taking a connected sum with this $S^n$ along a small tube.

Suppose that $Y:=X^{(n-1)}\cup\sigma_1\cup\dots\cup\sigma_{i-1}$ is embedded
via $g$, and that $\tl g$ where defined is equivariantly homotopic to
the restriction of $\Phi$.
Let $Z$ be $Y$ minus the interior of the simplicial neighborhood of
$\sigma_i$.
Then the restriction of $\tl g$ to $Z\x\partial\sigma_i$ extends to
$Z\x\sigma_i$, and so is null-homotopic.
On the other hand, this map is homotopic to the composition
$Z\x\partial\sigma_i
@>\psi\x\id_{\partial\sigma_i}>>S^n\x\partial\sigma_i
@>\tl\chi|_{S^n\x\partial\sigma_i}>>S^{2n-1}$, where
$\psi$ is the composition of the inclusion $Z\i\R^{2n}\but\partial\sigma_i\i
S^{2n}\but\partial\sigma_i$ and a deformation retraction of
$S^{2n}\but\partial\sigma_i$ onto a small $S^n$ winding around $\partial\sigma_i$
with linking number $1$; and $\chi$ is the inclusion
$S^n\sqcup\partial\sigma_i\emb\R^{2n}$.
Since $\tl\chi|_{S^n\x\partial\sigma_i}$ has degree one,
$\psi\x\id_{\partial\sigma_i}$ must be trivial in $2n$-cohomology.
Then $\psi\:Z\to S^n$ is trivial in $n$-cohomology, so by the Hopf
classification theorem, $\psi$ is null-homotopic.
Hence the inclusion of $Z$ in $S^{2n}\but\partial\sigma_i$ is null homotopic.

If $H$ is a homeomorphism of $\R^{2n}$ fixing $\partial\sigma_i$ and sending
$Z^{(n-1)}$ into a small ball $B^{2n}$, by dimensional reasons there is
a null-homotopy $h_t$ of $H(Z)$ in $\R^{2n}\but\partial\sigma_i$ such that
$Z^{(n-1)}$ stays in $B^{2n}$.
If $D^n$ is a disk bounded by $\partial\sigma_i$ and disjoint from $B^{2n}$,
the null-homotopy $H^{-1}h_t$ of $Z$ is such that $Z^{(n-1)}$ stays disjoint
from $g(\sigma_i):=H^{-1}(D^n)$.
Hence $g(\sigma_i)$ has zero intersection number with each $n$-simplex of
$Z$.

Then $g(\sigma_i)$ can be made disjoint from each $n$-simplex $\tau$ of $Z$
by the Whitney trick.
Specifically, every pair of double points between $\sigma_i$ and $\tau$ with
opposite signs can be eliminated by the price of adding a $1$-handle to
$\sigma_i$, that is, cutting out $B^n\x S^0$ and gluing in $S^{n-1}\x B^1$.
Since $n>2$, every circle in the modified $\sigma_i$ bounds a $2$-disk in
$\R^{2n}\but Z$, hence all the $1$-handles can be cancelled without
introducing new self-intersections.

Finally, $\sigma_i$ can be made disjoint from every adjacent $n$-simplex
$\tau$ of $Y$ by the Penrose--Whitehead--Zeeman trick.
Specifically, if $g$ sends $p\in\sigma_i$ and $q\in\tau$ to the same point in
$\R^{2n}$, let $J\i\sigma_i\cup\tau$ be an arc with endpoints $p,q$ and
with $J\cap\sigma_i\cap\tau$ consisting of precisely one point $v$.
Then a small regular neighborhood of $\rlap{\vrule height 2.5pt depth-1.75pt width20pt}{g(J)}$%
\footnote{Added in v5: Since $n>2$, the embedded $1$-sphere $g(J)$ bounds 
an embedded $2$-disk $D$ that meets $g(Y\cup\sigma_i)$ only $\partial D$.
Then a small regular neighborhood of $D$ in $\R^m$ is a ball $B^{2n}$.}
in $\R^m$ is a ball $B^{2n}$,
and its preimage in $\sigma_i\cup Y$, which is a regular neighborhood of $J$
in $\sigma_i\cup Y$, is a cone $v*P$.
Redefining $g$ on $v*P$ by mapping it conically into
$B^{2n}\cong v*\partial B^{2n}$ eliminates the double point $g(p)=g(q)$.
This yields an extension of $g$ to $\sigma_i$, and it can be adjusted like
before so that $\tl g$ is still equivariantly homotopic to the restriction
of $\Phi$. \qed
\enddemo

\definition{Van Kampen obstruction \cite{vK}, \cite{Sh}, \cite{Wu}}
Let $G\:\bar X\to\RP^\infty$ be any map classifying the line bundle $\lambda$
associated to the double covering $\tl X\to\bar X$.
The {\it van Kampen obstruction} $\tta(X)\in H^{2n}(\bar X;\Z)$ is defined to be
$G^*(\xi)$, where $\xi$ is the generator of $H^{2n}(\RP^\infty;\Z)\simeq\Z/2$.
Since $G$ is unique up to homotopy, $\tta(X)$ is well-defined.
Given an embedding $g\:X\emb\R^{2n}$, we have $\tta(X)=\bar g^*(\xi)=0$.
So $\tta(X)$ obstructs embeddability of $X$ into $\R^{2n}$.
\enddefinition

A useful interpretation of $\tta(X)$, due to \v Svarc \cite{Sch}, can be given
in terms of

\definition{Cohomological sectional category}
Recall from \S2 that if $t$ is a free PL involution on a polyhedron $P$,
the Euler class $e(\lambda_p)$ of the line bundle $\lambda_p$ associated
with the double covering $p\:P\to P/t$ is an element of $H^1(P/t;\Z_p)$,
where $\Z_p$ is the orientation sheaf of $p$.
Note that $e(\lambda_p)$ may be defined as the image of the generator
of $H^1(\RP^\infty;\Z_q)\simeq\Z/2$ under a classifying map of $\lambda_p$,
where $q\:S^\infty\to\RP^\infty$ is the double covering.
Indeed, this generator is nothing but $e(\gamma)$, where $\gamma=\lambda_q$ is
the tautological line bundle over $\RP^\infty$.

We denote the maximal $k$ such that $e(\lambda_p)^k\ne 0$ by $\secat_\Z(P,t)$.
(Note that $e(\lambda_p)^k\in H^k(P/t;\,\Z_p^{\otimes k})$, where
$\Z_p^{\otimes k}$ is isomorphic to $\Z_p$ when $k$ is odd and to $\Z$ when $k$
is even, cf.\ \S2.)
This number is also known as the Conner--Floyd cohomological co-index of
$(P,t)$ over $\Z$ \cite{CF} or James' ``Euler index'' of $p$
\cite{Ja}.%
\footnote{Added in v5: As I learned from S. Parsa, the same number is also known 
as the Smith index.}
\enddefinition

Now we are ready for \v Svarc's intrinsic description of the van Kampen
obstruction.
Since $\xi=e(\gamma)^{2n}$, we have $\tta(X)=e(\lambda)^{2n}=e(2n\lambda)$,
the Euler class of the orientable vector bundle
$(2n)\lambda\:\tl X\x_{\Z/2}\R^{2n}_-\to\bar X$.
Thus $\tta(X)\ne 0$ if and only if $\secat_\Z(\tl X)=2n$.
By similar arguments, if $X$ embeds in $\R^m$, then $\secat_\Z(\tl X)<m$.

\proclaim{Theorem 3.2 (Shapiro \cite{Sh}, Wu \cite{Wu})}
Let $X$ be an $n$-polyhedron.
If $\tta(X)=0$ and $n>2$, then $X$ embeds in $\R^{2n}$.
\endproclaim

\demo{Proof}
Suppose that $G^*(\xi)=0$ for some map $\bar X\to\RP^\infty$
classifying $\lambda$.
We may assume that $G$ sends the $(2n-1)$-skeleton of $\bar X$ into
$\RP^{2n-1}$.
Since $\pi_1(\RP^{2n-1})$ acts trivially on $\pi_{2n-1}(\RP^{2n-1})$,
by definition $G^*(\xi)$ is the primary obstruction to the existence of
a map $F\:\bar X\to\RP^{2n-1}$, coinciding with $G$ on $\bar X^{(2n-2)}$.
By the non-homotopically-simple obstruction theory (as described e.g\. in
the Hilton--Wylie textbook), such an $F$ exists; if $n>1$, it still
classifies $\lambda$.
By the covering theory, $F$ lifts to an equivariant map $\tl X\to S^{2n-1}_-$.
Since $n>2$, by the Haefliger--Weber Criterion 3.1(a), $X$ embeds into
$\R^{2n}$. \qed
\enddemo

\demo{Alternative proof}
The existence of an equivariant map $\tl X\to S^{2n-1}$ is
obviously equivalent to the existence of a cross-section of the bundle
$\tl X\x_{\Z/2} S^{2n-1}_-\to\bar X$ (with fiber $S^{2n-1}$).
The primary obstruction to the latter is the Euler class of this bundle, which
is well-known to be complete in this situation.
But it is the same as the Euler class $e(2n\lambda)$, which as we have seen
above coincides with $\theta(X)$. \qed
\enddemo

When $n=2$, $\tta(X)$ is incomplete \cite{FKT}.
When $n=1$, $\theta(X)$ is complete, due to the classical Kuratowski--Pontryagin
Theorem: either $X$ embeds in $\R^2$ or it contains a copy of either $K_5$
(the $1$-skeleton of the $4$-simplex) or $K_{3,3}$ (the join of the $0$-skeleta
of $2$ copies of the $2$-simplex).
Indeed, $\theta(K_5)$ and $\theta(K_{3,3})$ are non-zero by the following

\example{Example 3.3}
If $J$ is an $n$-dimensional join of $n_i$-skeleta of
$(2n_i+2)$-simplices, then $\theta(J)\ne 0$.
Indeed, suppose that $\theta(J)=0$.
Then by the proof of Theorem 3.2, $\tl J$ admits an equivariant map to
$S^{2n-1}_-$.
Hence so does the {\it simplicial deleted product} $\tl J_s$, that is,
the subpolyhedron of $\tl J$ consisting of the products of all pairs of
disjoint simplices of $J$.

Let us consider the {\it simplicial deleted join} $\hat J_s$, that is,
the subpolyhedron of the join $J*J$ consisting of the joins of all pairs
of disjoint simplices of $J$.
The suspension $\Sigma\tl J_s$ is $\Z/2$-homeomorphic to
$\hat J_s/(J*\emptyset\cup\emptyset*J)$.
Hence $\hat J_s$ admits an equivariant map to $\Sigma\tl J_s$, and therefore to
$\Sigma S^{2n-1}_-=S^{2n}_-$.
But this cannot be by the Borsuk--Ulam Theorem, since $\hat J_s$ is well-known
to be equivariantly homeomorphic to $S^{2n+1}_-$ \cite{F1}, \cite{F2},
\cite{Gr} (see also \cite{Ro} and \cite{dL}+\cite{Mat}).%
\footnote{The construction of this homeomorphism deserves a thorough analysis,
which in turn has many other interesting consequences (compare remark at the end
of \S4); this will be the subject of a forthcoming paper by the author.
(Added in v5: arXiv:1103.5457.)}
\endexample

\remark{Remark} In what follows, we will freely use the well-known fact
that for any polyhedron $X$ with a fixed triangulation, $\tl X$ equivariantly
deformation retracts onto $\tl X_s$ (see \cite{Hu} or \cite{S1; 5.3.c}; beware
that the proof given in \cite{Sh} is incorrect).
\endremark

\definition{$1$-Parameter van Kampen obstruction}
Let $f,g\:X\emb\R^{2n+1}$ be embeddings.
We define $\tta(f,g)$ to be the primary obstruction to the existence of a
homotopy between $\bar f,\bar g\:\bar X\to\RP^{2n}$.
In more detail, let us pick a homotopy $h\:\bar X\x I\to\RP^\infty$ between
their compositions with the inclusion $\RP^{2n}\i\RP^\infty$.
Let $\Psi$ be a generator of $H^{2n+1}(\RP^\infty,\RP^{2n};\,\Z_q)\simeq\Z$.
Then $\tta(f,g)$ is defined to be the image of $h^*(\Psi)$ under the Thom
isomorphism
$H^{2n+1}(\bar X\x I,\bar X\x\partial I;\,\Z_p)\simeq H^{2n}(\bar X;\,\Z_p)$,
where $p\:\tl X\to\bar X$ is the double covering.

Clearly, $\tta(f,g)$ is independent of the choice of $h$, and vanishes if $f$
and $g$ are isotopic.
It follows from Criterion 3.1(b) that the converse holds when $n>1$; this was
first proved in \cite{Wu}.
When $n=1$, the converse holds if $X\x I$ is allowed to be amended by taking
connected sum with tori (attached to open $2$-disks in $X\x I$) \cite{Ta1}.
Furthermore, when $n=1$, it is shown in \cite{ST} (see also \cite{Ta1; 5.2}) that
$\tta(f,g)=0$ iff $\tta(f|_K,g|_K)=0$ for each subgraph $K$ of $X$ homeomorphic
to $K_5$, $K_{3,3}$ or $S^1\sqcup S^1$.

The group $H^{2n}(\bar X;\,\Z_p)$ also contains $\bar f^*(\Xi)$ and
$\bar g^*(\Xi)$, where $\Xi$ is a generator of $H^{2n}(\RP^{2n};\,\Z_q)\simeq\Z$.
Up to a sign, we may assume that $\delta^*(\Xi)=2\Psi$.
It then follows that $2\tta(f,g)=\bar g^*(\Xi)-\bar f^*(\Xi)$.
When $n=1$, there is no $2$-torsion in $H^{2n}(\bar X;\Z_p)$ \cite{Ta1}, so
in this case $\tta(f,g)$ reduces to $\bar g^*(\Xi)-\bar f^*(\Xi)$.
\enddefinition

\definition{Yang index \& unoriented van Kampen obstruction}
If $t$ is a free involution on a polyhedron $P$ and $p\:P\to P/t$ is
the double covering, the first Stiefel--Whitney class
$w_1(\lambda_p)\in H^1(\bar X;\Z/2)$ is generally less informative than
$e(\lambda_p)$ (see \S2).

The {\it Yang index} $\secat_{\Z/2}(P,t)$ is the maximal $k$ such that
$w_1(\lambda_p)^k\ne 0$ (see \cite{CF}).
Thus $\secat_{\Z/2}(\tl X)=2n$ iff $\tta(X)$ has trivial $\bmod2$ reduction.
If $M$ is a closed manifold, $\secat_{\Z/2}(\tl M)$ equals the minimal $d$
such that the normal Stiefel--Whitney classes $\bar w_i(M)=0$ for $i\ge d$
\cite{Wu'}, \cite{Sch; Ch.\ VII, \S2} (see also \cite{Wu; paper II}; for
a geometric proof of the inequality $\secat_{\Z/2}(\tl M)\ge d$ see \cite{Mc4}).

We note that if $\tl X$ has Yang index $<2n-1$ then $\tta(X)=0$.
Indeed, since the generator $\xi\in H^{2n}(\RP^{2n};\Z)$ is of order two,
it is the Bockstein image of the generator
$w_1(\gamma)^{2n-1}\in H^{2n-1}(\RP^{2n};\,\Z/2)$.
Hence $\tta(X)$ is the Bockstein image of $w_1(\lambda_p)^{2n-1}$.
(This also yields an alternative intrinsic definition of $\tta(X)$.)
\enddefinition

\example{Example 3.4} Krushkal \cite{Kr} relates the intersection pairing of
a $4$-thickening of a $2$-polyhedron $X$ with the image of $\tta(X)$ in
$H^4(\tl X;\Q)$ under the homomorphism induced by the $2$-covering
$\tl X\to\bar X$.
This image is identically trivial since $\tta(X)$ is of order two.
In fact, since $\tta(X)$ is divisible by $e(\lambda)$, the image of $\tta(X)$
in the integral cohomology $H^4(\tl X;\Z)$ is zero as well by
the Smith sequence.
\endexample

\definition{Geometric definition of $\tta(X)$}
Pick an embedding $g\:X\emb\R^{2n+1}$.
Consider the equivariant map $\tl g\:\tl X\to S^{2n}$, assigning to a pair
$(x,y)$ of distinct points of $X$ the unit vector in the direction from
$g(x)$ to $g(y)$.
Then $\tta(X)$ can be identified with
$\tl g^*_\eq(\xi_\eq)\in H^{2n}_{\Z/2}(\tl X;\Z)$, induced from
the generator $\xi_\eq\in H^{2n}_{\Z/2}(S^{2n};\Z)\simeq\Z/2$.

If $g$ is generic, it projects onto an immersion $f\:X\imm\R^{2n}$.
We may assume that $f$ maps the vertices of some triangulation of
$X$ in general position and extends linearly to each simplex.
Consider the cochain $\tta_f$ on $\tl X$ assigning to a product
$\sigma\x\tau$ of distinct oriented $n$-simplices of the triangulation
of $X$ the algebraic number of their intersections in $\R^{2n}$ (i.e\.
$0$ if they are mapped disjointly and $1$ or $-1$ otherwise, according
as the orientation of $\sigma\x\tau$ matches or not that of $\R^{2n}$).
Then $\tta_f(\tau\x\sigma)=(-1)^n\tta_f(\sigma\x\tau)$, but also%
\footnote{This point was apparently missed in \cite{FKT}, leading to
an incorrect identification of the group $H^{2n}(\bar X;\Z)$ that contains
$\theta(X)$ as $H^{2n}(\bar X;\Z_p)$ for odd $n$.}
$t(\tau\x\sigma)=(-1)^n\sigma\x\tau$ (as chains).
Thus $\tta_f$ lies in the subcomplex $\Hom_\Lambda(C_*(\tl X_s),\,\Z)$ of
$C^*(\tl X_s)$, and it is immediate that $\tta_f=\tl g^*_\eq(\Xi_\eq)$, where
$\Xi_\eq$ is a cocycle representative of $\xi_\eq$ with support in two
antipodal $2n$-simplices of $S^{2n}\but S^{2n-1}$.
In conclusion, $\tta(X)=[\tta_f]$.

\smallskip
Using the language of \S2, the idea of construction of $\tta_f$ can be
formulated in a more elegant way.
Let $\Delta_f/t$ be the set of unordered pairs $\{x,y\}$ of points of $X$ such
that $f(x)=f(y)$ and $x\ne y$.
We have $\Delta_f/t=\bar g^{-1}(pt)$, where $pt\in\RP^{2n}$ is the vertical
direction of $\R^{2n+1}$.
If $f$ is generic, $\Delta_f/t$ is an embedded $2n$-pseudo-comanifold in
$\bar X$, which is co-oriented once $pt\in\RP^{2n}$ is.
Thus $[\Delta_f/t]=\bar g^*([pt])=\tta(X)\in H^{2n}(\bar X)$.
\enddefinition

\example{Example 3.5} The $n$-skeleton $K$ of the $(2n+2)$-simplex
$\Delta^n*\Delta^{n+1}$ is, apart from the simplex $\Delta^n$, already
contained in the $2n$-sphere $\partial\Delta^n*\partial\Delta^{n+1}$.
This simplex can be mapped conewise onto $c*\partial\Delta^n$, where $c$ is
a point in the interior $U$ of some $n$-simplex in $\partial\Delta^{n+1}$.
This yields a map $K\to\R^{2n}$ with one generic double point between
non-adjacent simplices.

More generally, an $n$-dimensional join $J$ of $n_i$-skeleta of
$(2n_i+2)$-simplices $\Delta^{n_i+1}*\Delta^{n_i}$ is, apart from the join of
the simplices $\Delta^{n_i}$, already contained in the $2n$-sphere
$(\bigast\partial\Delta^{n_i+1})*\partial(\bigast\Delta^{n_i})$ (since
$\partial(A*B)$ contains $A*\partial B$).
The remaining simplex can be mapped conewise onto
$c*\partial(\bigast\Delta^{n_i})$, where $c$ is a point in the interior of
some $n$-simplex in $(\bigast\partial\Delta^{n_i+1})$.
This again yields a map $J\to\R^{2n}$ with one generic double point
between non-adjacent simplices.
\endexample

\example{Example 3.6} The $\bmod2$ reduction of $\tta(G)$ is a complete
obstruction to embeddability of the graph $G$ into $\R^2$, since it is nonzero
for $G=K_5$ and $K_{3,3}$ by Example 3.3.
We will show that this is not the case in every other dimension $n>1$.

In the notation of the preceding example, cut out a small $n$-ball from
$U\but c$, and reattach this cell $B^n$ back to its former boundary in $K$
by some degree two map $S^{n-1}\to S^{n-1}$.
Denote the result by $X$.
Then $\tta(X)=2\alpha$, where $\alpha$ is represented by a cocycle with
support in $B^n\x\Delta^n\cup\Delta^n\x B^n$.

It remains to show that $\tta(X)\ne 0$.
Suppose, on the contrary, that $X$ embeds in $\R^{2n}$.
Let $g$ be the restriction of this embedding to
$X\but\Int B^n=K\but\Int B^n$.
If $l$ is the linking number between the $(n-1)$-sphere $g(\partial B^n)$ and
the $n$-sphere $g(\partial\Delta^{n+1})$, then $2l=0$ by the construction.
Hence $l=0$, contradicting the following
\endexample

\proclaim{Lemma 3.7} Let $X$ be an $n$-polyhedron with a fixed triangulation
and such that $\tta(X)\ne 0$.
Suppose that for some $n$-simplex $\sigma$ of $X$ we are given an embedding
$g\:X\but\Int\sigma\emb\R^{2n}$, and let $Z$ denote $X$ minus the interior
of the simplicial neighborhood of $\sigma$.
Then $\lk(g(\partial\sigma),g(Z))\ne 0\in H^n(Z)$.
\endproclaim

This linking number is defined to be the image of a fixed generator under
$$H^{2n-1}(S^{2n-1})@>\tl g^*>>H^{2n-1}(\partial\sigma\x Z)\simeq H^n(Z).$$
If $X=K$ and $\sigma=\Delta^n$, then $Z$ is the $n$-sphere
$\partial\Delta^{n+1}$ and we clearly get the usual linking number.

\demo{Proof}
Suppose that the linking number is zero, and write $Y=X\but\Int\sigma$.
Then the restriction of $\tl g\:\tl Y\to S^{2n-1}$ to $\partial\sigma\x Z$
is null-homotopic by the Hopf classification theorem.
Therefore $\tl g$ equivariantly extends to
$\tl Y\cup\sigma\x Z\cup Z\x\sigma$, which contains the simplicial deleted
product $\tl X_s$.
Hence $\tta(X)=0$, a contradiction. \qed
\enddemo

\remark{Remark} The following explicit computation of $\tta(X)$ has been
proposed recently \cite{MTW}, based on the well-known fact that every set of
at most $2n+1$ distinct points on the ``moment curve'' $\gamma(\R)\i\R^{2n}$,
$\gamma(t)=(t,t^2,\dots,t^{2n})$, is linearly independent.
Fixing some injection $\phi\:X^{(0)}\to\R$ of the $0$-skeleton of some
triangulation of $X$ and extending $X^{(0)}@>\phi>>\R@>\gamma>>\R^{2n}$
linearly to the simplices of $X$, one thus obtains a specific generic map
$f\:X\to\R^{2n}$.
It follows (see \cite{MTW}) that $\tta(X)=[\tta_f]$ is the class of the cochain
$c\in C^{2n}_{\Z/2}(\tl X_s)$ defined by $c(\sigma\x\tau)=1$ or $(-1)^n$ if
the triple $(\R,\phi(\sigma^{(0)}),\phi(\tau^{(0)}))$ is homeomorphic to
$(\R,2[n],2[n]+1)$ or $(\R,2[n]+1,2[n])$ respectively, where
$[n]=\{0,1,\dots,n-1\}$, and by $c(\sigma\x\tau)=0$ in all other cases.
As remarked in \cite{MTW}, the same result can also be obtained directly
from \v Svarc's interpretation of $\theta(X)$.
\endremark
\medskip

Of course, describing an explicit cocycle representing $\tta(X)$ still leaves
the question of deciding whether its class is trivial in $H^{2n}(\bar X)$.

\proclaim{Theorem 3.8} Let $X$ be an $n$-polyhedron with a fixed
triangulation and a choice of orientation of simplices.
Let $\Gamma$ be the graph with vertices corresponding to the $n$-simplices
of $X$ and edges to non-adjacent pairs of $n$-simplices.
Then $$H^{2n}(\bar X)\simeq\coker\left[C_1(\Gamma;\Z)@>h>>
\bigoplus_{\sigma\in\Gamma^{(0)}}H^n(X\but N\sigma)\right],$$
where $N$ stands for open simplicial neighborhood and $h$ sends
a non-adjacent pair $(\sigma^n,\tau^n)$ to the difference between
$[\sigma]\in H^n(X\but N\tau)$ and $(-1)^n[\tau]\in H^n(X\but N\sigma)$.
\endproclaim

It is convenient to think of each group $H^n(X\but N\sigma)$ as being put
into the vertex of $\Gamma$ corresponding to $\sigma$, so that we may regard
$h$ as the ``boundary homomorphism''.
(This can be made precise: $\coker h$ is the $0$-homology of $\Gamma$
with the pre-cosheaf coefficients.)

\example{Example 3.9} If $X=K_5$, then $\Gamma$ is the Petersen graph $P$.
If we orient the edges of $K_5$ according to a representation of $K_5$
as the union of two cycles of length $5$, then $5$ of the cycles of length $3$
will be oriented coherently and the other $5$ not.
Let us fix the generator of each $H^1(\text{cycle of length }3)$ corresponding
to the orientations of the $3$ edges, if they are coherent, and to those of
the $2$ coherently oriented edges otherwise.
This puts a $\Z$ with a preferred generator into each vertex of $P$.
Each edge of $P$ corresponds to an isomorphism between the groups in
its endpoints; due to the $(-1)^n$ multiplier, it identifies the preferred
generators iff precisely one of the vertices corresponds to a coherently
oriented cycle of length $3$.
This puts a sign onto each edge of $P$, and it is easy to check that
there is a cycle in $P$ which carries an odd number of minus signs.
Thus some preferred generator gets identified with negative itself,
and we conclude that $H^2(\bar K_5)\simeq\Z/2$.
\endexample

\demo{Proof of Theorem 3.8} We have
$H^{2n}(\bar X)\simeq H^{2n}(\bar X_s,\bar Y_s)\simeq
H^{2n}_c(\bar X_s\but\bar Y_s)$, where $Y=X^{(n-1)}$, and
$\bar X_s$ denotes $\tl X_s/t$.
Now
$\bar X_s\but\bar Y_s=\bigcup_{\sigma^n<X}W_\sigma$, where
$W_\sigma=[(\Int\sigma)\x (X\but N\sigma)\cup
(\Int\sigma)\x (X\but N\sigma)]/t$.
Since $W_\rho\cap W_\sigma\cap W_\tau=\emptyset$ for any pairwise distinct
$\rho,\sigma,\tau$, we have the Mayer--Vietoris exact sequence
$$\bigoplus_{\{\sigma^n,\tau^n\}}H^{2n}_c(W_\sigma\cap W_\tau)\to
\bigoplus_{\sigma^n}H^{2n}_c(W_\sigma)\to
H^{2n}_c(\bar X_s\but\bar Y_s)\to 0.$$
We have $H^{2n}_c(W_\sigma)\simeq H^n(X\but N\sigma)$
by the K\"unneth formula, and
$H^{2n}_c(W_\sigma\cap W_\tau)\simeq H^{2n}_c(\Int\sigma\x\Int\tau)\simeq\Z$
when $\sigma$ and $\tau$ are non-adjacent and $W_\sigma\cap W_\tau=\emptyset$
otherwise.
Finally, given orientations on $\sigma$ and $\tau$ such that
$[\tau]\in H^n(X\but N\sigma)$ and $[\sigma]\in H^n(X\but N\tau)$ coincide
with some fixed generators, then $$[\sigma\x\tau]\in
H^{2n}_c((\Int\sigma)\x (X\but N\sigma))\simeq H^{2n}_c(W_\sigma)$$ and
$$[\tau\x\sigma]\in H^{2n}_c((\Int\tau)\x (X\but N\tau))\simeq
H^{2n}_c(W_\tau)$$ may be assumed to coincide with fixed generators.
At the same time their preimages $$[\sigma\x\tau]\in
H^{2n}_c(\Int\sigma\x\Int\tau)\simeq H^{2n}_c(W_\sigma\cap W_\tau)$$
and $$[\tau\x\sigma]\in H^{2n}_c(\Int\tau\x\Int\sigma)\simeq
H^{2n}_c(W_\sigma\cap W_\tau)$$ differ by $(-1)^n$ in the right hand group.
\qed
\enddemo

The same argument works to prove

\proclaim{Addendum 3.10} Let $\tl\Gamma$ be the $2$-cover of $\Gamma$ induced
from the $2$-covering $S^1\to\RP^1$ by a map $f\:\Gamma\to\RP^1$ such that
$f^{-1}(\RP^0)=\Gamma^{(0)}$.
Then $$H^{2n}(\tl X)\simeq\coker\left[C_1(\tl\Gamma;\Z)@>\tl h>>
\bigoplus_{\sigma\in\tl\Gamma^{(0)}}H^n(X\but N\sigma)\right],$$
where $\tl h$ lies over $h$.
Moreover, $H^{2n}(\bar X)@>p^*>>H^{2n}(\tl X)$ is given by
$C_1(\bar X)@>\pi_!>> C_1(\tl X)$.
\endproclaim

It follows that $\ker p^*$, which is a group containing $\tta(X)$, is
generated by elements $a\in H^n(X\but N\sigma)$ in the vertices
$\sigma\in\Gamma^{(0)}$ that are homologous (with pre-cosheaf coefficients)
to $-a$ via loops in $\Gamma$ that do not lift to $\tl\Gamma$.
Note that a (graph-theoretic) cycle in $\Gamma$ lifts to $\tl\Gamma$ iff
it has an even length.
\medskip

For further computations of cohomology of deleted product see
\cite{Bau}, \cite{Um}, \cite{Ni}, \cite{BF} and references there.

\head 4. Linkless and panelled embeddings\footnotemark \endhead
\footnotetext{This section was considerably altered in v5.}

An interesting theme emerged in combinatorial embedding theory in the early 1980s when J. H. Conway
and C. Gordon, and independently H. Sachs observed that every embedding of the complete graph $K_6$ 
in $\R^3$ links some pair of disjoint (graph-theoretic) cycles of $K_6$ with an odd linking number.
Amazingly, Robertson, Seymour and Thomas proved that if a graph $G$ embeds in $\R^3$ in such a way 
that no pair of disjoint cycles is linked with an odd linking number, then it also admits an embedding
in $\R^3$ that is {\it panelled} (or ``flat'') in the sense that the image of every cycle of $G$ bounds 
a disk in $\R^3$ whose interior is disjoint from the image of $G$ \cite{RST}.
Let us note that every panelled embedding is {\it knotless}, that is, contains no nontrivially knotted cycles.

Robertson, Seymour and Thomas also showed that if two panelled embeddings of a graph in $\R^3$
are inequivalent then they differ already on some subgraph that is isomorphic to
a subdivision of one of the Kuratowski graphs $K_5$, $K_{3,3}$ \cite{RST}.
In particular, if $G$ is planar, then all its panelled embeddings in $\R^3$ are equivalent.
Most famously, Robertson, Seymour and Thomas proved that a graph admits a panelled embedding in $\R^3$
if and only if it has no minor among the seven graphs known as the {\it Petersen family} \cite{RST}.

The ``panelled'' in the latter result can be replaced, as we just discussed, by ``admitting an embedding
that links no pair of disjojnt cycles with an odd linking number'' and also by any intermediate property,
including the following one.
We say that an embedding $g$ of an $n$-polyhedron $X$ in $\R^m$ is {\it linkless} if for every two disjoint 
closed subpolyhedra of $g(X)$, one is contained in an $m$-ball disjoint from the other one.
Thus a linkless embeding of a graph in $\R^3$ need not be knotless, but it may not contain a non-split link of 
any number of cycles (such as the Whitehead link or the Borromean rings).%
\footnote{But not conversely: an embedding of a graph in $\R^3$ whose restriction to any disjoint union of cycles is 
a split link does not have to be linkless, see \cite{FN; Figure 2.1 (2)}.}
Consequently, it may not contain a link of any number of cycles that is nontrivial up to PL isotopy. 

\proclaim{Lemma 4.1} A panelled embedding of a graph in $\R^3$ is linkless.%
\footnote{The published version of the present paper also claimed a converse: a linkless, knotless embedding
of a graph in $\R^3$ is panelled; the last step in the proof of this claim was erroneous. 
The author is grateful to R. Nikkuni for pointing out that this assertion is false, as shown, for instance, by
Kinoshita's $\Theta$-curve \cite{Ki; Figure 1}, which is a knotless (and, trivially, linkless) embedding of 
$K_{3,2}$ in $\R^3$ that is not equivalent to the standard embedding, and consequently is not panelled.}
\endproclaim

\demo{Proof} Let $g\:G\emb\R^3$ be a panelled embedding.
 If $H$ is a subgraph of a graph $G$, by $\bar H$ we shall denote the union of all 
edges of $G$ disjoint from $H$.
It suffices to prove that for each subgraph $H$ of $G$, there exists
a $3$-ball containing $g(H)$ and disjoint from $g(\bar H)$.
Let $T$ be a spanning forest of $H$ (i.e., a union of spanning trees
of the components of $H$) and let $E_1,\dots,E_r$ be the edges of $H$
that are not in $T$.
Let $C_i$ be the unique graph-theoretic cycle in $T\cup E_i$.
Note that each $A_{ij}:=C_i\cap C_j$ is connected, for it is the intersection
of $T\cap C_i$ and $T\cap C_j$ which are connected and whose union contains no
cycles.
Let $H^+$ be the $2$-complex obtained from $H$ by glueing up $C_1,\dots,C_r$
by $2$-cells $D_1,\dots,D_r$, then $H^+$ collapses onto $T$, which is
collapsible onto finitely many points.
Hence it suffices to prove that $g|_H$ extends to an embedding
$h\:H^+\emb\R^3\but g(\bar H)$, for a regular neighborhood of its image will
then be a union of $3$-balls, which can be connected by thin tubes to get 
the desired $3$-ball.
The existence of $h$ follows from a well-known lemma of B\"ohme and H. Saran (see \cite{B\"o}
or \cite{RST; (2.8)}); for convenience of the reader, we include a short proof.

Suppose that $g|_H$ has been extended to an embedding
$h_i\:H\cup D_1\cup\dots\cup D_{i-1}\emb\R^3\but\bar H$.
By the hypothesis, $g|_{C_i}$ extends to an embedding $d\:D_i\emb\R^3$
such that $d(D_i\but C_i)$ is disjoint from $g(G)$.
We may assume that $d$ is transverse to each $h_i(D_j)$, relative to
the boundary.
Since $A_{ij}$ is connected, every connected component of
$h_i(D_j\but A_{ij})\cap d(D_i)$ is
null-homologous as a locally finite cycle in $h_i(D_j\but A_{ij})$.
Then at least one such component $S$ is innermost in $h_i(D_j\but A_{ij})$,
that is, the unique graph-theoretic cycle $S'$ in $S\cup g(A_{ij})$ bounds
a disk $D_S$ in $h_i(D_j)$ whose interior is disjoint from $d(D_i)$.
Redefine $d$ by replacing the disk bounded by $S'$ in $d(D_i)$ with $D_S$,
and pushing it slightly off $h_i$ (relative to $g(A_{ij})$).
This decreases the number of components in $h_i(D_j\but A_{ij})\cap d(D_i)$
at least by one, so if we continue this process, it will terminate with a new
$d$ that combines with $h_i$ to yield the desired embedding $h_{i+1}$. \qed
\enddemo

Let us say that an embedding $g$ of an $n$-polyhedron $X$ into $\R^{2n+1}$
is {\it $Y$-panelled,} where $Y$ is an $(n+1)$-polyhedron, $Y\supset X$, if
$g$ extends to a map $f\:Y\to\R^{2n+1}$ such that $f^{-1}(X)=X$.
A $Y$-panelled embedding of $X$ into $\R^{2n+1}$ yields an equivariant map
$\tl X^Y\to S^{2n}$, where $\tl X^Y=(X\x Y\cup Y\x X)\but\Delta_X$.
Let $\lambda_r$ be the line bundle associated with the double covering
$r\:\tl X^Y\to\bar X^Y$, where $\bar X^Y=\tl X^Y/t$.
Thus $$\eta(X,Y):=e(\lambda_r)^{2n+1}\in H^{2n+1}(\bar X^Y;\Z_p)$$
vanishes if $X$ admits a $Y$-panelled embedding into $\R^{2n+1}$.

Let us call a subpolyhedron $S$ of an $n$-polyhedron $X$ {\it linkable} if
$S\x (X\but S)$ admits an essential map to $S^{2n}$ (or equivalently, 
$H^n(S)\otimes H^n(X\but S)$ is nonzero).
Let us call a simplicial complex $K$ {\it bounded} if every linkable
subcomplex of $K$ is homeomorphic to a quotient of a PL $n$-manifold $M$ with boundary
by some identification on $\partial M$.

For $n>1$ it is not hard to see that an embedding of $g$ an $n$-polyhedron $X$
in $\R^{2n+1}$ is linkless if it is $X_+$-panelled, where $X_+$ is obtained
from $X$ by glueing up by cones all linkable subpolyhedra of $X$, triangulated
by subcomplexes of a fixed triangulation of $X$.%
\footnote{Indeed, if $P$ and $Q$ are disjoint subcomplexes of $X$, then either
$P\x Q$ admits no essential map to $S^{2n}$ or $g(P)$ is null-homotopic in
the complement to $g(Q)$.
In both cases $\tl g|_{P\x Q}$ is null-homotopic, so by the
Haefliger--Weber Criterion 3.1, $g|_{P\sqcup Q}$ is equivalent to the embedding
$h\:P\sqcup Q\emb\R^{2n+1}$, obtained by combining $e_1g|_P$
and $e_2g|_Q$, where $e_1,e_2\:B^{2n+1}\to\R^{2n+1}$ are embeddings with disjoint
images.}
It is also not hard to see that the converse holds if $n>1$ and $X$ has a bounded
triangulation.%
\footnote{Indeed, let $S$ be a linkable subcomplex of a bounded triangulation of $X$.
Since $g$ is linkless, $g(S)$ lies in a ball $B$ disjoint from $g(\bar S)$, 
where $\bar S$ is the union of all simplices disjoint from $S$.
Then $g|_S$ extends to a map $f\:CS\to B$, and since $S$ is the quotient of a manifold by
some identification on the boundary, using the Penrose--Whitehead--Zeeman trick
it is easy to achieve that $f^{-1}(S)=S$ (thus, $g|_S$ is $CS$-panelled) and that
$f(CS)$ is disjoint from $g(X\but S)$.}
Thus for $X$ with a bounded triangulation and $n>1$, $\eta(X):=\eta(X,X_+)$ vanishes if $X$ admits 
a linkless embedding in $\R^{2n+1}$.

\proclaim{Theorem 4.2} Let $X$ be an $n$-polyhedron and $Y$ an
$(n+1)$-polyhedron containing $X$.
\smallskip

(a) Let $n>1$.
Then $X$ admits a $Y$-panelled embedding in $\R^{2n+1}$ iff $\eta(X,Y)=0$.
\smallskip

(b) $X$ admits a linkless embedding in $\R^{2n+1}$ if $\eta(X)=0$.
If $X$ has a bounded triangulation,%
\footnote{The published version of the present paper erroneously omitted this hypothesis.
To see that it cannot be omitted, note that $\eta(K_5\sqcup K_5)\ne 0$ by part (a) 
and Theorem 4.4 below.
The author is grateful to V. Turchin for calling his attention to the proof of
Theorem 4.2(b).
A more satisfactory correction, imposing no restrictions on $X$, appears in 
a subsequent paper by the author, arXiv:1103.5457v2,
where the proof of Theorem 4.6 shows that $X$ admits a linkless embedding in
$\R^{2n+1}$ if and only if the following cohomology class $\eta'(X)$ vanishes.
Fix some triangulation $K$ of $X$.
If $S$ is triangulated by a subcomplex of $K$, let $\bar S$ be the union
of all simplices of $K$ disjoint from $S$.
Let $\tl X_s^+$ be obtained from $\tl X_s$ by glueing up by cones all products
$S\x\bar S$ that admit an essential map into $S^{2n}$.
Then $\eta'(X)$ is the $(2n+1)$th power of the Euler class of the line bundle
associated with the double cover $\tl X_s^+\to\bar X_s^+$.
}
then the converse also holds.
\endproclaim

\demo{Proof. (a)} It follows from the proof of the Haefliger--Weber Criterion 3.1
that if $\tl X^Y$ admits an equivariant map to $S^{2n}$ and $n>1$, then $X$
admits a $Y$-panelled embedding into $\R^{2n+1}$. \qed
\enddemo

\demo{(b)} The case $n>1$ follows from (a) and the preceding remarks.
Suppose that $n=1$.
If $\eta(X)=0$, the proof of (a) works to construct
an embedding $X\emb\R^3$ such that every pair of cycles has zero linking number.
By \cite{RST}, this suffices to re-embed $X$ linklessly in $\R^3$. 
Conversely, if $X$ admits a linkless embedding in $\R^3$, then by \cite{RST} it
admits a panelled embedding $g$ in $\R^3$.
If additionally $X$ has a bounded traingulation, then every linkable subcomplex
of this triangulation is homeomorphic to $S^1$.
Hence $g$ is $X_+$-panelled, and thus $\eta(X)=0$. \qed
\enddemo

\example{Example 4.3} Consider the complete graph $K_6$ as the $1$-skeleton of
the cone $CK_5$.
Since $CK_5$ contains no pair of disjoint $2$-simplices, every embedding of
$K_6$ in $\R^3$ may be considered ``$CK_5$-linkless''.
However, $K_6$ admits no $CK_5$-panelled embedding in $\R^3$, since
$\secat_\Z(\widetilde{CK_5})=3$ (see Example 3.3).
\endexample

\proclaim{Theorem 4.4} For $n>0$, an $n$-polyhedron $X$ admits a $CX$-panelled
embedding into $\R^{2n+1}$ if and only if $\tta(X)=0$.
\endproclaim

Note that the case $n=2$ is not an exception here.

\demo{Proof}
Note that
$\tl X^Y_s:=\bigcup\{A\x B, B\x A\mid A\in X, B\in Y,\, A\cap B=\emptyset\}$
is equivariantly homotopy equivalent to $\tl X^Y$.
On the other hand $\tl X^{CX}_s$ coincides with the simplicial deleted
product $\widetilde{CX}_s$.
By Lemma 4.5 below, for $n>0$ there is an isomorphism
$H^{2n+1}_{\Z/2}(\widetilde{CX}_s;\,I)\simeq
H^{2n+1}_{\Z/2}(\Sigma(\tl X_s);\,I)$ induced by an equivariant map, $I$
being the augmentation ideal of $\Z[\Z/2]$, so
$\secat_\Z(\widetilde{CX}_s)=\secat_\Z(\Sigma(\tl X_s))$.
It is not hard to see that $\secat_\Z(\Sigma(\tl X_s))=\secat_\Z(\tl X_s)+1$,
cf.\ \cite{CF; (5.1)}.
It follows that $\secat_\Z(\tl X)=\secat_\Z(\tl X^{CX})-1$.
Hence $\tta(X)=0$ iff $\tl X^{CX}$ admits an equivariant map to $S^{2n}$.
This implies the ``only if'' part of the theorem.

The ``if'' part now follows from Theorem 4.2(a) when $n>1$.
For $n=1$ we note that if $X$ contains a subgraph $K$ homeomorphic to
either $K_5$ or $K_{3,3}$, then $\widetilde{CX}_s$ contains
$\widetilde{CK}_s$, which is equivariantly homeomorphic to $S^3$
(see Example 3.3).
This contradicts the existence of an equivariant map
$\widetilde{CX}_s\to S^2$.
So $X$ embeds into $\R^2$ and $CX$ embeds into $\R^3$. \qed
\enddemo

\proclaim{Lemma 4.5} There exists an equivariant map
$p\:\widetilde{CX}_s\to\Sigma(\tl X_s)$ whose relative mapping cylinder
has the same (generalized) equivariant cohomology as $(\Sigma X, pt)\x\Z/2$.
\endproclaim

\demo{Proof} We define $p$ by sending $X\x\{c\}$ and $\{c\}\x X$ onto
the suspension points and via the homeomorphism on their complement.
The relative mapping cylinder of $p$ equivariantly collapses onto
$(CX\cup_{\{c\}\x X}\widetilde{CX}_s\cup_{X\x\{c\}}CX,\,\widetilde{CX}_s)$,
which is excision-equivalent to $(CX\sqcup CX,\,X\sqcup X)$. \qed
\enddemo

The following is a generalization of Lemma 3.7.

\proclaim{Proposition 4.6} Let $K$ be a simplicial complex triangulating
an $n$-polyhedron $X$.
Let $S$ and $L$ be the subpolyhedra of $X$ triangulated by the star $\st(v,K)$
and the link $\lk(v,K)$ of some vertex $v$ of $K$.
Set $Y=(X\but S)\cup L$.

(a) If $g\:Y\emb\R^{2n}$ is an embedding and
$(\bar g|_{\bar L^Y})^*\:H^{2n-1}(\RP^{2n-1})\to H^{2n-1}(\bar L^Y)$
is the zero map, then $\tta(X)=0$.

(b) If $Y$ embeds into $\R^{2n}$ and $L$ admits a $Y$-panelled embedding into
$\R^{2n-1}$, then $\tta(X)=0\pmod 2$.
\endproclaim

\demo{Proof. (a)}
Since $H^{2n-1}(\RP^{2n-1})$ maps onto $H^{2n}(\RP^\infty,\RP^{2n-1})$,
the latter has trivial image in $H^{2n}(\bar S^X,\bar L^Y)$.
This image is precisely the obstruction to extending $\bar g|_{\bar L^Y}$ to
$\bar S^X$.
Since $\bar Y\cup\bar S^X=\bar X$, we conclude that $\tta(X)=0$. \qed
\enddemo

\demo{(b)} Given a $g\:Y\emb\R^{2n}$, our
$(\bar g|_{\bar L^Y})^*\:H^{2n-1}(\RP^{2n-1};\,\Z/2)\to H^{2n-1}(\bar L^Y;\,\Z/2)$
factors as
$H^{2n-1}(\bar L^Y;\,\Z/2)@>G_*>>H^{2n-1}(\RP^\infty;\,\Z/2)@>\simeq>>
H^{2n-1}(\RP^{2n-1};\,\Z/2)$, where $G\:\bar L^Y\to\RP^\infty$ is
any classifying map.
Hence modulo $2$, $(\bar g|_{\bar L^Y})^*$ does not depend on the choice of
$g$, and the proof of (a) applies. \qed
\enddemo

\example{Example 4.7} If $X$ is the $n$-skeleton of $\Delta^{2n+2}$,
and $v$ any its vertex, then $L$ is the $(n-1)$-skeleton of $\Delta^{2n+1}$,
and $Y=L_+$, the $n$-skeleton of $\Delta^{2n+1}$.
Since $\tta(X)$ is nonzero even modulo $2$ (see \S3), we conclude that $L$
admits no linkless embedding into $\R^{2n-1}$.
This was originally proved in \cite{LS; Corollary 1.1} and \cite{Ta2} by
different arguments.
\endexample

\remark{Remark}
Interestingly, an argument of converse type was found in \cite{Sk1}, where
non-embeddability of a certain polyhedron is derived from non-existence of
linkless embeddings of its links of vertices.
\endremark

\remark{Remark}
Let $K$ be an $n$-dimensional simplicial complex.
If $S$ is a subcomplex of $K$, let $S'$ be the union of all simplices of $K$
disjoint from $S$.
Consider all subcomplexes $S_1,\dots,S_r$ of $K$ such that $S_i\x S_i'$
admits an essential map to $S^{2n}$.
Let $\hat K$ be obtained from the simplicial deleted product $\widetilde{CK}_s$
by attaching $CS_i\x CS_i'$ to each $(CS_i)\x S_i'\cup S_i\x (CS_i')$.
If $L$ is the $n$-skeleton of the $(2n+3)$-simplex, then $\hat L$ is
$\Z/2$-homeomorphic to $S^{2n+2}_-$.
If $G$ is one of the seven graphs Petersen graphs, then $\hat G$ is
$\Z/2$-homeomorphic to $S^4_-$ in six cases, and $\Z/2$-homotopy equivalent but
not homeomorphic to $S^4_-$ in the remaining case.
This provides yet another proof, in the spirit of Example 3.3, that $L$ does not
linklessly embed in $\R^{2n+1}$ and $G$ does not linklessly embed in $\R^3$.
The details will appear in a forthcoming paper of the author, where these
homeomorphisms are unified in one general result.
\endremark

\head 5. Arrow diagram formulas for type $1$ invariants of singular knots \endhead

Let us fix a closed $1$-manifold $M$.
All results of this section are interesting already in the case $M=S^1$.
A {\it chord diagram} $\Theta$ on $M$ is a collection of $m$ unordered pairs
of distinct points of $M$.
Each pair is called a {\it chord}.
When $M=S^1$, $\Theta$ is called {\it reducible} if it can be partitioned into
two nonempty chord diagrams $\Theta_1$, $\Theta_2$ such that every chord of
$\Theta_1$ is unlinked with every chord of $\Theta_2$.
Associated to $\Theta$ is the equivalence relation $R_\Theta$ on $M$, where
$xR_\Theta y$ iff $x=y$ or $\{x,y\}\in\Theta$.
A {\it $\Theta$-knot} is a smooth immersion $\phi\:M\looparrowright\R^3$ that
factors through an embedding of the graph $G:=M/R_\Theta$ and has no
self-tangencies at the double points of $\phi$.
We consider $\Theta$-knots up to {\it $\Theta$-isotopy}, that is smooth regular
homotopy through $\Theta$-knots.

If $v$ is an invariant of $\Theta$-knots with values in an abelian group,
the {\it Vassiliev derivative} $v'$ is an invariant of $\Theta^+$-knots for each
chord diagram $\Theta^+$ obtained by adding one chord $\{x,y\}$ to $\Theta$.
It is defined by $v'(k)=v(k_+)-v(k_-)$, where $k_+$ and $k_-$ are
$\Theta$-knots, $C^1$-close to the $\Theta^+$-knot $k$ and such that
$(k_+'(x),\,k_+'(y),\,k_+(x)-k_+(y))$ is a positive frame and
$(k_-'(x),\,k_-'(y),\,k_-(x)-k_-(y))$ is a negative frame.
(``Positive'' and ``negative'' clearly do not depend on the choice of
the ordering of $x$ and $y$ and the choice of the orientation of $M$.)
The invariant $v$ is a {\it type $k$ invariant} if the $(k+1)$th derivative
$v^{(k+1)}$ vanishes identically.

We will be concerned with type $1$ invariants of $\Theta$-knots, which
include (the restrictions of) the $m$th Vassiliev derivatives of all type
$m+1$ invariants of genuine knots.
The derivative of a type $1$ invariant $v$ of $\Theta$-knots may be regarded
as a function $v'(\{x,y\})$ of the extended chord diagram
$\Theta^+=\Theta\cup\{\{x,y\}\}$.

\definition{Configuration spaces $\tl\Theta$ and $\bar\Theta$}
For convenience, we shall assume hereafter that each component of $M$
contains at least one chord (the general case can be treated similarly).
Then $G=M/R_\Theta$ is a graph, possibly with loops and multiple edges.
Let $\tl\Theta_0=G\x G\but\bigcup C\x C$, where $C$ runs over all open
edges of $G$, and let $\bar\Theta_0=\tl\Theta_0/(\Z/2)$.
On the other hand, let us fix some triangulations of $M$ and $G$ such that
the quotient map $q\:M\to G$ is simplicial.
Let us consider $(q\x q)(\tl M_s)$, the image of the simplicial deleted product
of $M$ in $G\x G$.
Let $\tl\Theta_0^s=(q\x q)(\tl M_s)\cap\tl\Theta_0$ and
$\bar\Theta_0^s=\tl\Theta_0^s/(\Z/2)$.

Each of $\tl\Theta_0$, $(q\x q)(\tl M_s)$ and $\tl\Theta_0^s$ meets the diagonal
of $G\x G$ precisely in points of the type $(w,w)$, where $w=q(x)=q(y)$ is
a vertex of $G$ covered by the points of some chord $\{x,y\}$ of $\Theta$.
For every such $(w,w)$, let us consider its links $L_{(w,w)}$ and
$L^s_{(w,w)}$ respectively in $\tl\Theta_0$ and in $\tl\Theta_0^s$.
Each of them is endowed with a free involution, induced from
the factor-exchanging involution of $G\x G$.
Clearly, $L^s_{(w,w)}$ is $\Z/2$-homeomorphic to the link of ${(w,w)}$ in
$(q\x q)(\tl M_s)$.

Let us temporarily use the following notation: $[\vl]=\{1,3\}$, $[-]=\{2,4\}$
and $[+]=[\vl]\cup[-]$.
The link of $w$ in $G$ is homeomorphic to $[+]$, consequently the link of
$(w,w)$ in $G\x G$ is homeomorphic to $[+]*[+]$.
Then $L_{(w,w)}$ is $\Z/2$-homeomorphic to the simplicial deleted join
$\widehat{[+]}_s$, and $L^s_{(w,w)}$ to its subset $[\vl]*[-]\cup[-]*[\vl]$.
Now the pair $(\widehat{[+]}_s,\,[\vl]*[-]\cup[-]*[\vl])$ is in turn
$\Z/2$-homeomorphic to (all edges of $I^3$, all horizontal edges of $I^3$),
with involution induced from the antipodal involution on $\partial I^3$.
The latter homeomorphism sends the first copy of $[\vl]$ to the vertices
of some diagonal of the top face of the cube, the first copy of $[-]$ to
the vertices of the perpendicular diagonal of the bottom face of the cube,
the second copies of $[\vl]$ and $[-]$ symmetrically, and all edges
linearly.
Thus the star $(w,w)*L_{(w,w)}$ of $(w,w)$ in the second derived subdivision
of $\tl\Theta_0$ is homeomorphic to the cone over the $1$-skeleton of the
$3$-cube, by a homeomorphism throwing $(w,w)*L^s_{(w,w)}$ to the cone over
the boundaries of the two horizontal faces of the cube.

Let $\tl\Theta$ be obtained from $\tl\Theta_0$ by removing, for every vertex
$w$ as above, the interior of the star $(w,w)*L_{(w,w)}$, identified with
the cone over the $1$-skeleton of the cube, and replacing it with the $4$
vertical faces of the cube.
The factor exchanging involution on $\tl\Theta_0$ carries over to a free
involution on $\tl\Theta$ by considering the antipodal involution on
the boundary of the cube.
We set $\bar\Theta=\tl\Theta/(\Z/2)$ and denote the double covering by
$p\:\tl\Theta\to\bar\Theta$.
\enddefinition

\remark{Remark} Let us discuss the geometry of $\tl\Theta^s\i\tl\Theta$ and
$\bar\Theta^s\i\bar\Theta$, which are similarly obtained from
$\tl\Theta_0^s\i\tl\Theta_0$.
At each $(w,w)$, we have in $\tl\Theta^s$ a singularity of the type
$C(S^1\sqcup S^1)$, which is resolved in the standard way, by glueing
in the handle $S^1\x I$.
The involution on this handle is orientation-reversing, so the quotient
is homeomorphic to the M\"obius band.
Thus $\bar\Theta^s$ is obtained from $\bar\Theta_0^s$ by blowing it up at
all points of the type $\{w,w\}$.

Returning to $\bar\Theta$, notice that the open star of $\{w,w\}$ in
$\bar\Theta_0$ is homeomorphic to
$$\R\x\R\x\{0\}\quad\cup\quad\R\x\{0\}\x[0,\infty)\quad\cup\quad\{0\}\x\R\x
(-\infty,0]$$
where $\R\x\R\x\{0\}$ is identified with the open star of $\{w,w\}$ in
$\bar\Theta_0^s$.
Then $\bar\Theta$ is obtained from $\bar\Theta_0$ by blowing up
the $\R\x\R\x\{0\}$ at the origin and then regluing $\R\x\{0\}\x[0,\infty)$ and
$\{0\}\x\R\x (-\infty,0]$ back along the images of $\R\x\{0\}\x\{0\}$ and
$\{0\}\x\R\x\{0\}$ under the blowup, and repeating for each $\{w,w\}$.

It follows that $\bar\Theta$ is a simplicial analogue of the space
$B_2(\Theta)$ from \cite{Va2}, which collapses onto $\bar\Theta$.
An attempt to define $\bar\Theta$ occurs in \cite{M3}, where $\bar\Theta^s$ was
defined instead.
\endremark

\proclaim{Lemma 5.1} Let $v$ be an integer-valued type $1$ invariant of
$\Theta$-knots.
Let $z_v^\#=\sum v'(\{c,d\})C\x D$, where the sum runs over all products $C\x D$
of components of $M\but\bigcup\Theta$ and each $c\in C$ and $d\in D$.
Then $z_v^\#$ is a skew-invariant $2$-cycle in $G\x G$, which lies in
$\tl\Theta_0$.
Moreover, it ``blows up'' to a unique skew-invariant $2$-cycle $z_v$ in
$\tl\Theta$, which may be identified with an element of $H_2(\bar\Theta;\,\Z_p)$.
\endproclaim

By a skew-invariant element of a $\Z[\Z/2]$-module $M$ we mean an $x\in M$
such that $t(x)=-x$, where $t$ is the generator of $\Z/2$.

\demo{Proof} Since $t(C\x D)=-D\x C$ as $2$-chains, $z_v^\#$ is a skew-invariant
chain in $\tl\Theta_0$.

An edge of the product cell structure on $\tl\Theta_0$ is of the form
$p\x C$ or $C\x p$, where $p$ is a point of $M/R_\Theta$ corresponding to
a chord of $\Theta$, and $C$ is a $1$-simplex of $T_\Theta$.
Such an edge is incident to four $2$-cells in $\tl\Theta_0$ of the form
$D\x C$ (respectively, $C\x D$) where $D$ is a $1$-simplex of $T_\Theta$.
Hence $\partial z_v^\#=0$ is a reformulation of the well-known (and obvious)
``four-term relation''.

The assertion that $z_v^\#$ lies in $\tl\Theta_0$ is a direct consequence
of the well-known (and obvious) ``one-term relation'', also known as the
``framing independence relation''.

To prove the final assertion, let us consider
$h\:H_2(\bar\Theta;\,\Z_p)\to H_2(\bar\Theta,\mu;\,\Z_p)$, where $\mu$ is
the union of the M\"obius bands obtained upon the blowup of the chords.
Since $\mu$ is homotopically $1$-dimensional, $h$ is injective.
Since $H_1(\mu;\,\Z_p)=0$ as well%
\footnote{This is the only point where the proof fails for $\Z/2$-valued type
$1$ invariants.}%
, it is an isomorphism.
Now since $(\tl\Theta,\nu)$, where $\nu$ is the preimage of $\mu$, is
$2$-dimensional, its $2$-dimensional cohomology classes are identified with
its $2$-cycles, and therefore the homomorphism $p^!$ from
$H_2(\bar\Theta,\mu;\,\Z_p)\simeq H_2^{\Z/2}(\tl\Theta,\nu;\,I)$
to $H_2(\tl\Theta,\nu)$ is an embedding onto the subgroup of skew-invariant
elements (see \S2).
Let $\tl\Theta^+$ be the ``bordism'' between $\tl\Theta_0$ and $\tl\Theta$,
obtained by adding the cone over each component of $\nu$, and let $\nu^+$ be
the union of these cones.
Then $H_2(\tl\Theta,\nu)\simeq H_2(\tl\Theta^+,\nu^+)$ by excision, and at
the same time $(\tl\Theta^+,\nu^+)$ admits an equivariant deformation retraction
onto $(\tl\Theta_0,\nu_0)$, where $\nu_0=\tl\Theta_0\cap\Delta_G$.
Since $\nu_0$ is $0$-dimensional,
$H_2(\tl\Theta_0,\nu_0)\simeq H_2(\tl\Theta_0)$. \qed
\enddemo

If $k\:M/\Theta\emb\R^3$ is a $\Theta$-knot, we may consider its Gauss maps
$\tl k\:\tl\Theta\to S^2$ and $\bar k\:\bar\Theta\to\RP^2$.
(The details of the construction are left as an exercise for the reader.)
We define the parametric van Kampen obstruction
$\zeta(k_1,k_2)\in H^2(\bar\Theta;\,\Z_p)$ to existence of a $\Theta$-isotopy
between the $\Theta$-knots $k_1,k_2\:M/\Theta\emb\R^3$ to be the first
obstruction to homotopy of their Gauss maps
$\bar k_1,\bar k_2\:\bar\Theta\to\RP^2$.

\proclaim{Theorem 5.2} If $v$ is an integer-valued type $1$ invariant of
$\Theta$-knots and $k_1$ and $k_2$ are $\Theta$-knots,
$v(k_1)-v(k_2)=\zeta(k_1,k_2)[z_v]$.
\endproclaim

\demo{Proof} Pick a generic homotopy $h_t$ between $k_1$ and $k_2$.
Let $k^i$ be the $\Theta\cup\{\{c_i,d_i\}\}$-knots occurring in $h_t$ and
$\eps_i$ be their signs corresponding to the increase of $t$.
Then $v(k_1)-v(k_2)=\sum_i\eps_iv'(k_i)=\sum_i\eps_iv'(\{c_i,d_i\})$.

On the other hand, the Thom isomorphism
$t\:H^2(\bar\Theta;\,\Z_p)\simeq
H^3(\bar\Theta\times I,\bar\Theta\times\partial I;\,\Z_p)$
identifies $\zeta(k_1,k_2)$ with $\bar H^*(\chi)$, where
$H\:M/\Theta\to\R^4$ is a $\Theta$-isotopy projecting to $h$ and
$\chi$ generates $H^3(\RP^\infty,\RP^2;\,\Z_p)\simeq\Z$.
Clearly, $\bar H^*(\chi)=\sum_i\eps_i\omega(C_i,D_i)$, where $\omega(C_i,D_i)$
is the class of the cocycle assuming $1$ on the skew-equivariant
chain $C_i\x D_i+D_i\x C_i$ and $0$ on all other basic skew-equivariant
chains. \qed
\enddemo

\proclaim{Corollary 5.3} Let $k_0$ be a fixed $\Theta$-knot.
Then $u(k):=\zeta(k,k_0)$ is a universal type $1$ invariant of $\Theta$-knots.
\endproclaim

In more detail, $u$ takes values in $H^2(\bar\Theta;\,\Z_p)$ (which is free
abelian, since $\bar\Theta$ is $2$-dimensional) and every integer-valued type $1$
invariant $v$ can be recovered from $u$.
In fact, the same arguments work for invariants $v$ with values in any abelian
group containing no elements of order $2$.
(Elements of order $2$ do not fit in the proof of Lemma 5.1.)

\demo{Proof} By the proof of Theorem 5.2,
$u(k_+)-u(k_-)=\zeta(k_+,k_-)$ is the class of the cocycle assuming $1$ on
the skew-invariant chain $C\x D+D\x C$, where $C$ and $D$ are the components
of $M\but\bigcup\Theta$ being intersected, and $0$ on all other basic
skew-invariant chains.
So $u$ is a type $1$ invariant.
By Theorem 5.2, any invariant $v$ of $\Theta$-knots is a function of $u$,
namely $v(k)=v(k_0)+u(k)[z_v]$.
Thus $u$ is universal. \qed
\enddemo

Given a chord $\alpha$ of $\Theta$ and some choice of orientation on it
(i.e.\ an ordering of its points), along with a fixed orientation of $M$,
let $\alpha^+$ and $\alpha^-$ denote the images in $G$ of the oriented
half-circles $a^+,a^-\i M$ with $\partial a^+=\partial a^-=\alpha$.

\proclaim{Corollary 5.4} (a) The group $\Gamma_1(\Theta)$ of integer-valued
type $1$ invariants of $\Theta$-knots modulo type $0$ invariants is isomorphic
to $H^2(\bar\Theta;\,\Z_p)$.

(b) \cite{Va2; Theorem 1(b)} If $M=S^1$, $\Gamma_1(\Theta)$ is free abelian
of rank ${m\choose 2}+n$, where $m$ is the number chords and $n$ is the number
of irreducible factors of $\Theta$.
\endproclaim

\demo{Proof. (a)} A map $\Gamma_1(\Theta)\to H^2(\bar\Theta;\,\Z_p)$
is given by Lemma 5.1; by construction, it is a homomorphism.
By Theorem 5.2, this map has an inverse. \qed
\enddemo

\demo{(b)} By the proof of Lemma 5.1, $H_2(\bar\Theta;\,\Z_p)$ is isomorphic
to the subgroup of skew-invariant elements of $H_2(\tl\Theta_0)$.
By the K\"unneth formula, $H_2(G\x G)$ is freely generated by products of
the type $[\alpha_i^+\x\alpha_j^+]$, where $\alpha_1,\dots,\alpha_m$ are
all chords of $\Theta$ with arbitrary but fixed orientations, and
$\alpha_0^+$ is the image of $[M]$ in $G$.
The subgroup of skew-invariant elements of $H_2(G\x G)$ is freely generated by
elements of the type $[\alpha_i^+\x\alpha_i^+]$ and
$[\alpha_i^+\x\alpha_j^+]+[\alpha_j^+\x\alpha_i^+]$ with $j>i$, and so is of
rank $m+1+{m+1\choose 2}$.
Clearly, $H_2(G\x G,\tl\Theta_0)$ is free abelian of rank $2m$, and all
its elements are skew-invariant.
It is easy to see that an edge $e$ of $G$ is such that
$[e\x e]\in H_2(G\x G,\tl\Theta_0)$ lies in the image of $H_2(G\x G)$ iff
it is the unique common edge of some pair of cycles $a,b\in H_1(G)$.

The irreducible factors of $\Theta$ correspond to arcs $J_1,\dots,J_n$ in $M$
with endpoints in $\bigcup\Theta$ such that each pair of arcs is either disjoint
or one arc is contained in the interior of the other; the $i$-th irreducible
factor of $\Theta$ consists of all chords with both endpoints in $J_i$ but not
in any $J_j\i J_i$.
Clearly, an edge $e$ of $G$ is the unique common edge of some pair of cycles
iff it is connecting chords from the same irreducible component.
Let $G'$ be the graph with $n$ vertices, obtained from $G$ by contracting all
such edges.
If $G'$ has $k$ edges, $H_1(G')$ is of rank $k-n+1$.
Since no edge of $G'$ is the unique common edge of some pair of cycles,
$H_1(G')$ has a basis where no two elements have a common edge.
Let $N$ be the union of products of distinct edges of $G'$.
Then $H_2(G'\x G',N)$ has rank $k$ and the image of $H_2(G'\x G')$ in
$H_2(G'\x G',N)$ has rank $k-n+1$.
Hence $\coker(H_2(G'\x G')\to H_2(G'\x G',N))$ has rank $n-1$, and it
follows that $\coker(H_2(G\x G)\to H_2(G\x G,\tl\Theta_0))$ also has rank
$n-1$.
Hence $K:=\im(H_2(G\x G)\to H_2(G\x G,\tl\Theta_0))$ has rank $2m-n+1$.

The image of the skew-invariant subgroup of $H_2(G\x G)$ in $K$ contains
$2K$, hence has the same rank $2m-n+1$.
Since $H_2(\tl\Theta_0)$ maps injectively to $H_2(G\x G)$, we conclude that
the skew-invariant subgroup of $H_2(\tl\Theta_0)$ has rank
${m+1\choose 2}-m+n={m\choose2}+n$. \qed
\enddemo

\definition{Arrow diagram formulas}
A {\it one-arrow diagram} over $\Theta$ is an ordered pair $(C,D)$ of connected
components of $M\setminus\bigcup\Theta$.
It corresponds to the function $f_{(C,D)}$ on generic $\Theta$-knots, whose value
on $k\:M\looparrowright\R^3$ is the algebraic number of ways to draw
an upward oriented vertical segment in $\R^3$ whose initial endpoint
is in $k(C)$ and terminal in $k(D)$.
Each such segment is counted with the sign of the frame, consisting of this
oriented segment and the tangent vectors to $k$ at its endpoints.
In addition, the unique arrowless diagram over $\Theta$ corresponds to
the function assuming 1 on every $\Theta$-knot.

An {\it arrow diagram formula} for a type $1$ invariant $v$ of $\Theta$-knots is
a representation of $v$ as a sum $\sum_{(C,D)}f_{(C,D)}+\sum 1$ of the
functions corresponding to one-arrow and arrowless diagrams.
\enddefinition

\proclaim{Corollary 5.5 \cite{Va2; Theorem 5(1)}}
Every integer-valued type $1$ invariant $v$ of $\Theta$-knots admits an arrow
diagram formula with half-integer coefficients.
\endproclaim

\demo{Proof} Let $w=\sum_{(C,D)}v'(\{c,d\})(C,D)$, where $c\in C$, $d\in D$.
Clearly, $w'=2v'$, as long as we know that $w$ is an invariant of $\Theta$-knots.
To prove the latter, pick a suitable generator $\phi$ of $H^2(S^2)$.
Representing $\phi$ by a cocycle with support in the north pole, we see that
$w(k)=\tilde k^*(\phi)(p^![z_v])$, which is an invariant of $k$. \qed
\enddemo

\proclaim{Theorem 5.6} An integer-valued type $1$ invariant $v$ of $\Theta$-knots
admits an integral arrow diagram formula if and only if $[z_v]\Cap e(\lambda_p)=0$,
where $e(\lambda_p)\in H^1(\bar\Theta;\,\Z_p)$ is the Euler class of the line
bundle associated with the $2$-covering $p\:\tl\Theta\to\bar\Theta$.
\endproclaim

\demo{Proof} From the Smith sequence of $p$, whose connecting homomorphism is
$\cdot\Cap e(\lambda_p)$, we have $[z_v]=p_*(c)$ for some
$c\in H_2(\tilde\Theta)$.
Then $k\mapsto\tilde k^*(\phi)(c)$ is an arrow diagram formula for $v$
(up to a type $0$ invariant), since $\zeta(k_1,k_2)(p_*c)=
p^*\zeta(k_1,k_2)(c)=[\tilde k_2^*(\phi)-\tilde k_1^*(\phi)](c)$.

Conversely, by definition, an integral arrow diagram for $v$ is a $2$-chain
$c\in C_2(\tilde\Theta)$ such that $v(k)=\tilde k^*(\phi)(c)$ up to a type $0$
invariant.
Similarly to the proof of Lemma 5.1, $\partial c=0$.
Then by the preceding paragraph, $\zeta(k_1,k_2)(p_*[c])=\zeta(k_1,k_2)[z_v]$
for any $\Theta$-knots $k_1$, $k_2$.
Since $\zeta(k_1,k_2)$ can be an arbitrary skew-equivariant $2$-cochain
(see the proof of Corollary 5.3), $p_*([c])=[z_v]$ and therefore
$[z_v]\Cap e(\lambda_p)=0$. \qed
\enddemo

We call $\Theta$ {\it planar} if there exists a $\Theta$-knot in the plane,
$k\:M/\Theta\emb\R^2\i\R^3$.

\proclaim{Lemma 5.7} Let $\mu\i\bar\Theta$ be the union of the M\"obius bands
obtained upon the blowup of the chords of $\Theta$.

(a) If $M=S^1$, then $\Cap e(\lambda_p)\:H_2(\bar\Theta,\mu;\Z_p)\to
H_1(\bar\Theta,\mu)$ is trivial.

(b) If $M=S^1$ and $\Theta$ is irreducible,
$\ker [H_1(\bar\Theta)\to H_1(\bar\Theta,\mu)]/\text{\rm (odd torsion)}$
is cyclic.

(c) If $M=S^1$ and $\Theta$ is irreducible and planar, the cyclic group
from (b) is isomorphic to $\Z$, and consequently
$\Cap e(\lambda_p)\:H_2(\bar\Theta;\Z_p)\to H_1(\bar\Theta)$ is trivial.
\endproclaim

Part (a) replaces \cite{Va2; Lemma 4}, which is easier due to the fact
that genuine knots are replaced with long knots in \cite{Va2}.
The idea of proof of (b) is taken from \cite{Va2; proof of Lemma 6}.

\demo{Proof. (a)} By the proof of Lemma 5.1, $H_2(\tl\Theta,\nu)\simeq
H_2(\tl\Theta_0)$ over $\Lambda=\Z[\Z/2]$, where $\nu$ is the preimage of $\mu$
in $\tl\Theta$.
By the proof of Corollary 5.4(b), we have a short exact sequence of
$\Lambda$-modules
$$0\to H_2(\tl\Theta_0)\to H_2(G\x G)\to K\to 0,$$
where $K$ is isomorphic to the direct sum of $2m-n+1$ copies of the
augmentation ideal $I=\ker(\Lambda\to\Z)$.
The $\Z/2$-invariant subgroup $H^0(\Z/2;\,K)$ of $K$ is trivial, and from
the explicit description of generators given in the proof of 5.4(b),
the diagonal subgroup $D$ of $H_2(G\x G)=H_1(G)\otimes H_1(G)$ embeds onto
a subgroup $e(D)$ of $K$,
moreover $K/e(D)$ is isomorphic to the direct sum of $m-n$ copies of $I$.
Hence $H^0(\Z/2;\,K/e(D))=0$ and so $H^1(\Z/2;\,D)\to H^1(\Z/2;\,K)$ is injective.
On the other hand, $H^1(\Z/2;\,D)\to H^1(\Z/2;\,H_2(G\x G))$ is surjective,
since obviously $H_2(G\x G)/D$ is free over $\Lambda$.
Hence $H^1(\Z/2;\,H_2(G\x G))\to H^1(\Z/2;\,K)$ is injective.
Thus $H^1(\Z/2;\,H_2(\tl\Theta_0))$ is trivial.
Consequently, $H^1(\Z/2;\,H_2(\tl\Theta,\nu))$ is trivial.
The assertion now follows from Lemma 2.2. \qed
\enddemo

\demo{(b)}
The {\it intersection graph} of $\Theta$ has a vertex for every
chord of $\Theta$, and two vertices are connected by an edge iff the
corresponding chords are linked.
It is well-known and easy to see that $\Theta$ is irreducible iff its
intersection graph is connected.

Let $\alpha=\{x,y\}$ and $\beta=\{x',y'\}$ be a pair of linked chords of
$\Theta$ with some orientations.
Let us denote the union of the edges of $\alpha^+$ and $\beta^+$ by
$e_{\alpha^+\beta^+}$.
Then $\alpha^+\x\beta^+-\beta^+\x\alpha^+$ is a $\Z/2$-invariant cycle of
$G\x G$, which in fact lies in $\tl\Theta_0$.%
\footnote{Every such cycle reduced $\bmod 2$ gives rise to a $\Z/2$-valued
type $1$ invariant that does not lift to an integral-valued one.}
This cycle gives rise to an integral $2$-chain $R_{\alpha^+\beta^+}$ in
$\bar\Theta$ with
$\partial R_{\alpha^+\beta^+}=\pm c_\alpha\pm c_\beta+2\sum\pm c_{\gamma_i}$,
where $\gamma_i=\{x_i,y_i\}$ are chords of $\Theta$ with $(x_i,y_i)$
``lying'' in the support of $R_{\alpha^+\beta^+}$, and $c_\alpha$ denotes
the central curve of the M\"obius band corresponding to $\alpha$ (with some
choice of orientation).
Then each $\gamma_i$ is linked with either $\alpha$ or $\beta$ (or both),
and if $\gamma_i$ is linked with $\alpha$ (resp.\ $\beta$), then
$e_{\alpha^+\beta^+}$ contains $e_{\gamma_i^+\alpha^+}$ (resp.\
$e_{\gamma_i^+\beta^+}$) for a unique orientation of $\gamma_i$.
Inducting on the number of edges in $e_{\alpha^+\beta^+}$, we obtain
an equation of the form $k[c_\alpha]=l[c_\beta]$, where $k$ and $l$ are
odd integers.
Working modulo odd torsion, we may assume that $\gcd(k,l)=1$, and the assertion
follows. \qed
\enddemo

\remark{Remark} It is easy to see that a choice of orientations of $M$
and of a chord $\alpha$ determines an orientation of $c_\alpha$.
Furthermore, if all $\gamma_i$ are oriented as above, then it is not hard to see
that all signs in the formula for $\partial R_{\alpha^+\beta^+}$ are positive.
According to the assertion of \cite{Va2; Lemma 6}, it can be deduced from this
that $k,l=\pm1$, in which case $\ker [H_1(\bar\Theta)\to H_1(\bar\Theta,\mu)]$
contains no odd torsion.
\endremark

\demo{(c)} Clearly, $[c_\alpha]\Cap e(\lambda_p)\ne 0$, where $c_\alpha$ is
the central curve of the M\"obius band corresponding to $\alpha$.
Given a planar $\Theta$-knot $k\:G\emb\R^2$, it follows that the Gauss map
$\bar k$ sends $c_\alpha$ with nonzero degree to $\RP^1$.
Hence no multiple of $[c_\alpha]$ can be trivial in $H_1(\bar\Theta)$.

By (a), the image of $H_2(\bar\Theta;\Z_p)@>\Cap e(\lambda_p)>>H_1(\bar\Theta)$
is contained in the kernel of $H_1(\bar\Theta)\to H_1(\bar\Theta,\mu)$.
The latter kernel contains no elements of order $2$ by the above.
However, $2e(\lambda_p)=0$ as $\lambda_p$ is induced from the tautological
line bundle over $\RP^\infty$. \qed
\enddemo

\proclaim{Corollary 5.8 \cite{Va2; Theorem 5(2)}} If $\Theta$ is planar, every
integer-valued type $1$ invariant of $\Theta$-knots admits an integral
arrow diagram formula.
\endproclaim

\demo{Proof} If $\Theta$ has irreducible factors $\Theta^i$, we can pick
embeddings $e_i\:M/\Theta^i\to M/\Theta$ respecting the smooth structure at
each crossing of $M/\Theta^i$, and thus obtain embeddings
$e_{i*}\:H_2(\tl\Theta^i_0;\,\Z_p)\to H_2(\tl\Theta_0;\,\Z_p)$.
By the proof of Lemma 5.1, this yields embeddings
$E_i\:H_2(\bar\Theta^i;\,\Z_p)\to H_2(\bar\Theta;\,\Z_p)$.
By the proof of Corollary 5.4(b), $\Gamma_1(\Theta)$ is generated by
$E_i(\Gamma_1(\Theta^i))$ and by the classes of
$v(k)=\lk(k(\alpha_j^+),k(\beta_j^+))$, where $(\alpha_j,\beta_j)$ runs over
$\Theta^k\x\Theta^l$ for all $(k,l)$ with $k>l$, and each $\alpha_j$ and
$\beta_j$ are oriented so that $\alpha_j^+\cap\beta_j^+=\emptyset$.
Obviously (see e.g.\ the proof of Theorem 5.6), every such $v$ has
an integral arrow diagram formula.
On the other hand, by Theorem 5.6 and Lemma 5.7(c), each
$w\in\Gamma_1(\Theta^i)=H_2(\Theta^i;\,\Z_p)$ has an arrow diagram formula.
Then again by Theorem 5.6, $E_i(w)$ also has an arrow diagram formula. \qed
\enddemo

\proclaim{Lemma 5.9} (a) \cite{Va2; Theorem 7(1)} (a) If $v$ is a type $1$
invariant of $\Theta$-knots and $R$ is a reflection of $\R^3$, then
$v(k)+v(Rk)$ is a type $0$ invariant of $\Theta$-knots.

(b) \cite{Va2; Theorem 7(2)} Its parity $\propto(v)$ obstructs existence of
an integral arrow diagram formula for the integer-valued invariant $v$.
\endproclaim

In particular, $\propto(v^{(m)})$ is an obstruction to existence of integral
Polyak--Viro formulas for the type $m+1$ knot invariant $v$.

\demo{Proof. (a)} $v'(k_s)=-v'(Rk_s)$ for every singular knot $k_s$. \qed
\enddemo

\demo{(b)} If the reflection $R$ is taken to be in a vertical plane,
$A(Rk)=-A(k)$ for every one-arrow diagram function $A$ (only the signs of
the frames are changed). \qed
\enddemo

\proclaim{Theorem 5.10} If $M=S^1$ and $v$ is an integer-valued type $1$
invariant of $\Theta$-knots, $\propto(v)=[z_v]\Cap e(\lambda_p)^2$.
\endproclaim

\demo{Proof}
If $h$ is an isotopy with values in $\R^4$ between $\Theta$-knots $k_1$ and
$k_2$, we have
$2\zeta(k_1,k_2)=2t^{-1}\bar h^*(\chi)=t^{-1}\bar h^*(\delta^*\psi)=
\bar k_2^*(\psi)-\bar k_1^*(\psi)$, where
$\psi=p_!\phi$ generates $H^2(\R P^2;\Z_p)$, $\chi=\frac12\delta^*\psi$
and $t$ is the Thom isomorphism ($\chi$ and $t$ are as in the proof of
Theorem 5.2).
Now $\overline{Rk}^*(\psi)=-\bar k^*(\psi)$, thus
$$v(k)-v(Rk)=\zeta(k,Rk)[z_v]=\frac12[\overline{Rk}^*(\psi)-\bar k^*(\psi)]=
-\bar k^*(\psi)[z_v].$$
Since the $\bmod 2$ reduction of $\bar k^*(\psi)$ is $w_1(\lambda_p)^2$, we
obtain $\propto(v)=w_1(\lambda_p)^2[z_v]$.
Since $H_0(\bar\Theta;\Z_p)\simeq\Z/2$ due to $M=S^1$, we can identify
$\propto(v)=[z_v]\Cap w_1(\lambda_p)^2$ with $[z_v]\Cap e(\lambda_p)^2$. \qed
\enddemo

Let $a,b\:S^1\to G$ be a pair of oriented cycles in $M/R_\Theta$ with no common
edges.
A type one invariant $v_{ab}$ of $\Theta$-knots $k\:M/R_\Theta\emb\R^3$ is given
by the linking number between $C^1$-close pushoffs $ka_*$ and $kb_*$ of
$ka$ and $kb$, which are chosen as follows.
In those common vertices of $a(S^1)$ and $b(S^1)$ where $a$ and $b$ have corners,
we smoothen the corners of $ka$ and $kb$, and in each vertex $a(x)=b(y)$ of
transversal intersection between $a$ and $b$, we slide $ka$ and $kb$ away from
each other so that $(ka_*'(x),kb_*'(y),ka_*(x)-kb_*(y))$ form a positive frame.

Let $a,b\:S^1\to G$ be a pair of oriented graph-theoretic cycles in
$M/R_\Theta$ with disjoint edges.
We call them a {\it Manturov pair} if they have precisely one transversal
self-intersection.
(They may have any number of common vertices where each has a corner.)
Clearly, $v_{ab}(Rk)=\pm 1-v_{ab}(k)$, so $\propto(v_{ab})$ is nontrivial.

\example{Example 5.11} Let $\Theta_0$ be the unique irreducible diagram with two
chords on $M=S^1$.
Both basic type 1 invariants of $\Theta_0$-knots, which are yielded by the two
possible choices of Manturov cycles, have nonzero $\propto$.
However the basic type $3$ invariant $v_3$ of genuine knots has $v_3''$
equal to the sum of these two, so $\propto(v_3'')=0$.
In fact, $v_3$ is known to have an integral Polyak--Viro formula.
\endexample

\remark{Remark} It is computed in \cite{Va2; Example 5} that one of the three
basic invariants of knots of type $4$ satisfies $\propto(v_4''')\ne 0$.
In fact, $v_4'''=v_{ab}$ for a certain Manturov pair of cycles $a$, $b$.
Thus $v_4'''$ admits no integral arrow diagram formula, and in particular
$v_4$ admits no integral Polyak--Viro formula.
\endremark
\medskip

Obviously, if $M/\Theta$ admits a pair of Manturov cycles, $\Theta$ cannot
be planar.
Verifying Vassiliev's earlier conjecture, Manturov proved the converse, by
a nontrivial combinatorial argument.

\proclaim{Manturov's Criterion 5.12 \cite{Man}} Let $M=S^1$.
Then $\Theta$ is non-planar if and only if $M/R_\Theta$ contains
a Manturov pair of cycles.
\endproclaim

The van Kampen obstruction $\zeta(\Theta)\in H^2(\bar\Theta)$ to planarity of
$\Theta$ is induced from the generator of $H^2(\RP^\infty)$ by a classifying map
of the $2$-covering $p\:\tl\Theta\to\bar\Theta$.
Equivalently, $\zeta(\Theta)=e(\lambda_p)^2$, where
$e(\lambda_p)\in H^1(\bar\Theta;\Z_p)$ is the Euler class of the line bundle
associated with the $2$-covering $p\:\tl\Theta\to\bar\Theta$.

\proclaim{Corollary 5.13} Let $M=S^1$. Then $\Theta$ is planar iff
$\zeta(\Theta)=0$.
\endproclaim

\demo{Proof} If $k$ is a planar $\Theta$-knot,
$\bar k\:\bar\Theta\to\RP^1\i\RP^\infty$ classifies $p\:\tl\Theta\to\bar\Theta$,
hence $\zeta(\Theta)=0$.
If $\Theta$ is non-planar, let $v=v_{ab}$, where $a,b$ is a Manturov pair of
cycles.
By Theorem 5.10 we have $[z_v]\Cap\zeta(\Theta)=\propto(v)\ne 0$, hence
$\zeta(\Theta)\ne 0$. \qed
\enddemo

\proclaim{Corollary 5.14 \cite{Va2; Theorem 6}} Let $M=S^1$.
If $\Theta$ is non-planar, there exists an integer-valued type $1$ invariant
$v$ of $\Theta$-knots that admits no integral arrow diagram formula.
\endproclaim

\demo{Proof} Let $v=v_{ab}$, where $a,b$ is a Manturov pair.
Then $\propto(v)\ne 0$. \qed
\enddemo

\proclaim{Lemma 5.15} If $M=S^1$ and $\Theta$ is irreducible and non-planar,
the cyclic group from Lemma 5.7(b) is isomorphic to $\Z/2$.
\endproclaim

\demo{Proof} By Manturov's Criterion 5.12, there is a Manturov pair $a,b$, and
therefore a map $\bar\Theta^1\to\bar\Theta$, where $\Theta^1$ consists of
a single chord $\beta$ on $S^1\sqcup S^1$, with endpoints at distinct
components.
The component of $\bar\Theta^{1s}\i\bar\Theta^1$ containing the central curve
$c_\beta$ of the M\"obius band corresponding to $\beta$ is homeomorphic
to the blowup of $S^1\x S^1$ at a point.
Hence $2c_\beta$ bounds a punctured torus and $2[c_\beta]=0$.
Then $2[c_\alpha]=0\in H_1(\bar\Theta)$, where $\alpha$ is the transversal
intersection of $a$ and $b$. \qed
\enddemo

The following result was announced by Vassiliev \cite{Va2; Theorem 7(3)}, from
whom the author learned that his proof of this fact was inadvertently omitted
from \cite{Va2}.
The following proof was found by the author \cite{M3}.

\proclaim{Corollary 5.16} Suppose that $M=S^1$ and $\Theta$ is irreducible.
If $v$ is an integer-valued type $1$ invariant of $\Theta$-knots and
$\propto(v)=0\in\Bbb Z/2$, then $v$ admits an integral arrow diagram formula.
\endproclaim

\demo{Proof}
By Theorem 5.10, $[z_v]\Cap e(\lambda_p)^2=0$.
If $\Theta$ is planar, we are done by Corollary 5.8.
Otherwise Lemma 5.15 gives an explicit generator $[c_\alpha]$ of
$H_1(\bar\Theta)\simeq\Z/2$.
Clearly, $[c_\alpha]\Cap e(\lambda_p)\ne 0$.
Hence already $[z_v]\Cap e(\lambda_p)=0$.
By Theorem 5.6, $v$ admits an integral arrow diagram formula. \qed
\enddemo

\head 6. Extraordinary van Kampen obstruction \endhead

\definition{Equivariant stable cohomotopy}
Let $K$ be a $k$-polyhedron, homotopy equivalent to a compact polyhedron.
It is well known that if $k\le 2m-2$, the cohomotopy set $\pi^m(K):=[K,S^m]$
is an abelian group and $\Sigma\:\pi^m(K)\to\pi^{m+1}(\Sigma K)$ is
an isomorphism.
Hence the stable cohomotopy group
$\omega^m(K):=[\Sigma^\infty K,S^{m+\infty}]=[S^\infty\wedge K,S^{m+\infty}]$
is well-defined, where $\infty$ denotes a sufficiently large natural number,
and $m$ may be negative.
These groups form a generalized cohomology theory.

Throughout this section, $*$ will stand for the basepoint, and if $K$ is
an unpointed space, $K_+$ will denote the pointed space $K\sqcup *$.
We recall that the smash product of pointed spaces $P\wedge Q=P\x Q/P\vee Q$.
Note that when $P$ and $Q$ are compact and not of the form $K_+$, then
$P\wedge Q$ is the one-point compactification by $*$ of $(P\but *)\x (Q\but *)$.

Now suppose that $P$ is pointed and $G$ is a finite group acting on $P$
and fixing the basepoint.
We also assume that $P$ is $G$-homotopy equivalent to a compact polyhedron.
If $V$ is a finite-dimensional $\R G$-module, let $S^V$ be the one-point
compactification of the Euclidean space $V$ with the obvious action of $G$.
If $G$ acts trivially on $V$, one identifies $V$ with the integer $\dim V$.
The equivariant stable cohomotopy group
$$\omega^{V-W}_G(P):=[S^{W+V_\infty}\wedge P, S^{V+V_\infty}]_G$$ is
well-defined, where $V_\infty$ denotes a sufficiently large (with respect
to the partial ordering with respect to inclusion) finite-dimensional
$\R G$-submodule of the countable direct sum $\R G\oplus\R G\oplus\dots$
\cite{M+; p\. 81} (see also \cite{Ada}, \cite{Ko}, \cite{Hau}).
As emphasized in \cite{M+; IX.5 and XIII.1} and explained in
\cite{M+; pp\. 101--102}, one should not think of $V$, $W$ and $V_\infty$ as
abstract $\R G$-modules but rather as specific submodules of
$\R G\oplus\R G\oplus\dots$, in the same vein as $\pi_1(X)$ is technically
$\pi_1(X,*)$ if we care not only for the groups but also for their
homomorphisms.
However for our purposes it will be safe enough to ignore this point, as
is indeed done in later chapters of \cite{M+}, so that $\omega^*_G(P)$ is
effectively graded by the real representation ring of $G$.
\enddefinition

As explained in \cite{M+; p\. 100}, one cannot take $V_\infty$ to be any
$\R G$-module of sufficiently large dimension, even in the case $G=\Z/2$.
However:

\proclaim{Theorem 6.1} {\rm \cite{Hau} (cf.\ \cite{M+; IX.1.4}, \cite{Ada; 3.3})}
$[P,Q]_G\to [S^V\wedge P, S^V\wedge Q]_G$ is surjective if $\dim P^H\le\nu_H$
for all $H\i G$ and injective if $\dim P^H\le\nu_H-1$ for all $H\i G$, where
$\nu_H$ is the maximal integer such that

(1) either $V^H=0$ or $\pi_i(Q^H)=0$ for all $i$ with $2i+1\le\nu_H$; and

(2) either $V^K=V^H$ or $\pi_i(Q^K)=0$ for all $i\le\nu_H$ and all
$K\i H$.
\endproclaim

Let $kT+l$ denote $\R^k\x\R^l$ with the involution $(x,y)\mapsto (-x,y)$;
every $\R[\Z/2]$-module is of this form.

\proclaim{Corollary 6.2} $[P,S^{mT}]_{\Z/2}=\omega^{mT}_{\Z/2}(P)$ if
$\dim P\le 2m-2$ and $\dim P^{\Z/2}\le m-2$.
\endproclaim

\definition{Extraordinary van Kampen obstruction}
As in \S3, by $S^k_-$ we denote the $k$-sphere with the antipodal involution,
i.e\. the unit sphere in $(k+1)T$.
For $k\ge m$, we have that $S^k_-\but S^{m-1}_-$ is
$\Z/2$-homeomorphic to $S^{k-m}_-\x mT$.
Shrinking to points $S^{m-1}_-$ and each fiber $S^{k-m}_-\x\{pt\}$,
we get an equivariant map $\rho^k_m\:S^k_-\to S^{mT}$.

For any $k$-polyhedron $K$ with a free PL action of $\Z/2$ we have
an equivariant map $\phi_K\:K\to S^\infty_-$, which is unique up to equivariant
homotopy.
Here $\infty$ may be thought of as a sufficiently large natural number
(specifically, $k+1$ will do).
Given a compact $n$-polyhedron $X$ and a positive integer $m\ge n+2$, we let
$$\Theta^m(X):=[\rho^\infty_m\phi_{\tl X}]\in [\tl X_+,S^{mT}]_{\Z/2}=
\omega^{mT}_{\Z/2}(\tl X_+).$$
\enddefinition

\remark{Remark}
The geometric description of ordinary cohomology given in \S2 extends to
an {\it arbitrary} equivariant (with respect to a finite group) generalized
cohomology theory \cite{BRS}.
In particular, $\omega^{mT}_{\Z/2}(\tl X_+)$ is identified with the cobordism
group of $\Cal\Z_p^{\otimes m}$-co-oriented $m\lambda_p$-framed $m$-comanifolds
in $\bar X$, where ``comanifold'' is our synonym for ``mock bundle''
from \cite{BRS}, and an $m\lambda_p$-framing of an $m$-comanifold is
an isomorphism of its stable normal block bundle with the sum of $m$ copies of
the line bundle $\lambda_p$ (associated to the double covering
$p\:\tl X\to\bar X$).%
\footnote{Let us indicate a relation with a more traditional approach.
Let us consider the cobordism group of $\Cal\Z_p^{\otimes m}$-co-oriented
$m\lambda$-framed $m$-comanifolds in $\bar X$, where $\lambda$ is not
globally defined, but is a part of the data of the $m\lambda$-comanifold.
It is well-known that this group is isomorphic to
$[S^\infty*\bar X; S^\infty*(\RP^\infty/\RP^{q-1})]$ (see \cite{AE}).
Here $\RP^\infty/\RP^{q-1}$ is the Thom space of the bundle $q\gamma$ over
$\RP^{\infty-q}$, so the point-inverse $Q$ of the basepoint of this space is
$q\lambda$-framed in $\bar X$, where $\lambda$ is the pullback of $\gamma$
under the map $Q\to\RP^{\infty-q}$.}
Similarly to the geometric interpretation of $\tta(X)$, it is easy to see
that if $f\:X\to\R^m$ is a generic PL map that lifts to an embedding
$g\:X\emb\R^{2n+1}$, then $\Theta^m(X)$ is represented by the $m$-comanifold
$\Delta_f/t=\{\{x,y\}\in\bar X\mid f(x)=f(y)\}$ in $\bar X$, which is
$\Cal\Z_p^{\otimes m}$-co-oriented and $m\lambda$-framed as the preimage
of $\RP^{2n-m}$ under $\bar g\:\bar X\to\RP^{2n}$.
We also have $\Theta^m(X)=E(\lambda_p)^m$, where the equivariant Euler class
$E(\lambda_p)=i^*i_!([\id_{\bar X}])$ as in \S2.
\endremark

\proclaim{Theorem 6.3} Let $X$ be a compact $n$-polyhedron and
$m\ge\frac{3(n+1)}2$.
Then $X$ embeds in $\R^m$ if and only if $\Theta^m(X)=0$.
\endproclaim

This follows from the Haefliger--Weber Criterion 3.1 and

\proclaim{Lemma 6.4} Let $\Z/2$ act freely on a $k$-polyhedron $K$ and
let $k\le 2m-3$.
There exists an equivariant map $K\to S^{m-1}_-$ if and only if
$\rho^\infty_m\phi_K\:K\to S^{mT}$ is equivariantly null-homotopic.
\endproclaim

\demo{Proof} Let $S_m$ denote the ``cosphere'' $S^\infty_-/S^{m-1}_-$
with basepoint at the shrunk $S^{m-1}$.
Then $\rho^\infty_m$ is the composition $S^\infty_-@>q>>S_m@>p>>S^{mT}$,
where $q$ shrinks $S^{m-1}_-$ and $p$ shrinks all fibers
$S^{\infty-m}_-\x\{pt\}$ of its complement $S^{\infty-m}_-\x mT$.
All point-inverses of $p$ are weakly contractible, and it follows that $p$ is
a fibration.
Since $p$ is also equivariant, the given equivariant null-homotopy
$K\x I\to S^{mT}$ lifts to an equivariant homotopy
$F\:K\x I\to S_m$, which has to be a null-homotopy.
We are interested in equivariantly lifting it to $S^\infty_-$.
It certainly does lift over $K\x 0$ and over $A:=F^{-1}(S^{\infty-m}_-)$,
where $S^{\infty-m}_-$ is Hopf linked with the shrunk $S^{m-1}_-$.
Here $A$ may be assumed to be a $(k+1-m)$-polyhedron (after a small
perturbation of $F$).
Moreover, if $U$ is a $\Z/2$-invariant tubular neighborhood of
$S^{\infty-m}_-$ in the complement of the shrunk $S^{m-1}_-$,
then $F$ also lifts over $A^+:=K\x 0\cup F^{-1}(U)$.
Next we would like to extend this partial lift $A^+\to S^\infty_-$ over
the shadow of $A$ in $K\x I$, which is a $(k-m+2)$-polyhedron $B$.
The obstructions to this relative equivariant lifting problem lie in
the groups $H^i_{\Z/2}(B,A;\pi_{i-1}(S^{m-1}_-))$, which are all zero for
$i\le k-m+2$ since $k-m+2\le m-1$ by the hypothesis.
The obtained partial lift $B\cup A^+\to S^\infty_-$ of $F$ sends
the frontier of $B\cup A^+$ in $K\x I$ into the complement of
$S^{\infty-m}_-$.
Since $K\x I$ collapses onto $B\cup A^+$, this lift can be extended to
an equivariant map $K\x I\to S^\infty_-$ sending the exterior of
$B\cup A^+$ into the complement of $S^{\infty-m}_-$.
In particular, we get an equivariant map of $K\x\{1\}$ into
$S^\infty_-\but S^{\infty-m}_-$, which is equivariantly homotopy
equivalent to $S^{m-1}_-$. \qed
\enddemo

\remark{Remark} By the covering theory, there exists an equivariant map
$K\to S^{m-1}$ iff there exists a map $K/t\to\RP^{m-1}$ whose composition
with $\RP^{m-1}\i\RP^\infty$ is up to homotopy the classifying map $\psi_p$
of the line bundle associated to the $2$-covering $p\:K\to K/t$.
As shown by the following example, this is {\it not} equivalent to
null-homotopy of the composition
$K/t@>\psi_p>>\RP^\infty@>>>P_m$, where $P_m=\RP^\infty/\RP^{m-1}$.
The proof of Lemma 6.4 breaks down in this setting since $\pi_1(\RP^{m-1})$
is nontrivial as opposed to $\pi_1(S^{m-1}_-)$.
\endremark

\example{Example 6.5} (P. M. Akhmetiev)
Let $K/t$ be obtained from $\RP^{2n}$ by attaching a $(2n+1)$-cell along
the composition $S^{2n}@>>>S^{2n}@>>>\RP^{2n}$ of a degree $2$ map and
the $2$-covering.
Let $p\:K\to K/t$ be the universal covering.
Then $\psi_p^*\:H^{2n+1}(\RP^{\infty};\Z_p)\to H^{2n+1}(K/t;\Z_p)$ is
easily seen to be the inclusion $\Z/2\i\Z/4$, whence $\psi_p$ is not
homotopic to a map into $\RP^{2n}$.
On the other hand, $K/t@>\psi_p>>\RP^\infty@>>>P_{2n+1}$ may be
assumed to shrink $\RP^{2n}$ to the basepoint.
Since $\pi_{2n+1}(P_{2n+1})\simeq H_{2n+1}(P_{2n+1})\simeq\Z/2$, this
composition is null-homotopic.

Nevertheless, $K@>\phi_K>>S^\infty_-@>>>S_{2n+1}$
is not equivariantly null-homotopic since it sends
$H^{2n+1}_{\Z/2}(S_{2n+1};I)\simeq\Z$ onto the even
subgroup of $H^{2n+1}_{\Z/2}(K;I)\simeq\Z/4$.
This is not surprising as $S_{2n+1}\to P_{2n+1}$ is not a fibration.
\endexample

\definition{$\Theta^m(X)$ as cohomotopy Euler class}
Let $\Z/2$ act freely on a $k$-polyhedron $K$ with quotient $K/t$.
Consider the vector bundle $\xi\:(K\x mT)/t\to K/t$.
Let $\xi^{-1}$ be a vector bundle over $K/t$ such that $\xi\oplus\xi^{-1}$
is a trivial bundle of dimension $M$ say.
Let $\bar\xi^{-1}$ be the pullback of $\xi^{-1}$ over $K$.
Recall that the Thom space $\Th\xi$ is the one-point compactification of the
total space of $\xi$.
We have $\Th\bar\xi^{-1}\wedge S^{mT}=\Th(K\x M)=K_+\wedge S^M$.
Therefore
$$\omega^{mT}_{\Z/2}(K_+)=[K_+\wedge S^M,S^{mT+M}]_{\Z/2}\simeq
[\Th\bar\xi^{-1},S^M]_{\Z/2}=\omega^M_{\Z/2}(\Th\bar\xi^{-1})\simeq
\omega^M(\Th\xi^{-1}).$$
(See also \cite{Be} for a generalization of this isomorphism to the case of
actions with fixed points.)

In particular, this isomorphism identifies the extraordinary van Kampen
obstruction with an element of a non-equivariant stable cohomotopy group.
This element is well-known in the literature as the cohomotopy Euler class
of $\xi$ (see \cite{St}).
It can alternatively be defined as the class of the composition
$\Th\xi^{-1}\overset j\to\i\Th(\xi\oplus\xi^{-1})=(K/t)_+\wedge S^M\to S^M$.
If $\xi$ admits a nowhere vanishing cross-section, the inclusion
$K/t\i\Th\xi$ is null-homotopic, hence so is the inclusion $j$ in the above
composition.
The converse is known to hold when $k\le 2m-3$ \cite{Cr; Prop\. 2.4} (see also
\cite{Mey} for a nontrivial $\bmod p$ generalization).
This is in accordance with our Lemma 6.4, for it is easy to see that $\xi$
admits a nowhere vanishing cross-section if and only if $K$ admits an
equivariant map to $S^{m-1}$.
\enddefinition

\definition{Hurewicz homomorphism}
Let $\Z/2$ act freely on a $k$-polyhedron $K$ with quotient $K/t$.
Let us define $h\:\omega^{mT}_{\Z/2}(K_+)\to H^m_{\Z/2}(K;I^{\otimes m})$
by $h([\phi])=\phi^*(\xi^m)$, where
$\xi^m\in H^m_{\Z/2}(S^{mT};\pi_m(S^{mT}))\simeq\Z/2$ is the generator.
Clearly $h(\Theta^m(X))$ is the $m$th power of the Euler class of the line
bundle associated with the $2$-covering $\tl X\to\bar X$, in particular
$h(\Theta^{2n}(X))=\tta(X)$, under the isomorphism
$H^m_{\Z/2}(\tl X;I^{\otimes m})\simeq H^m(\bar X;\Z_p^{\otimes m})$.

One also has the natural transformation
$H\:\omega^{mT}_{\Z/2}(K_+)\to H^{mT}_{\Z/2}(K_+)$ into the equivariant
(representation graded) ordinary cohomology theory \cite{M+}, \cite{Ko}.
We claim that $H=h$.
Indeed, the above proof that $\omega^{mT}_{\Z/2}(K_+)=\omega^M(\Th\xi^{-1})$
works with any equivariant cohomology theory, in particular
$H^{mT}_{\Z/2}(K_+)\simeq H^M(\Th\xi^{-1})$.
On the other hand, using a twisted Thom isomorphism,
$H^M(\Th\xi^{-1})\simeq H^m(K/t;\Z_p^{\otimes m})$.
\enddefinition

\definition{Classification of embeddings} An isotopy classification of
embeddings in the metastable range is given by a generalized cohomology:
\enddefinition

\proclaim{Theorem 6.6} If a compact $n$-polyhedron $X$ embeds in $\R^m$ and
$m>\frac{3(n+1)}2$, the set of isotopy classes of embeddings of $X$ into
$\R^m$ is in bijection with $\omega^{mT-1}_{\Z/2}(\tl X_+)$.
\endproclaim

Note that
$\omega^{mT-1}_{\Z/2}(\tl X_+)\simeq\omega^{M-1}(\Th\xi^{-1})$ similarly to
the above.

\demo{Proof} By Corollary 6.2, $[\tl X_+\wedge S^1,S^{mT}]\simeq
\omega^{mT}_{\Z/2}(\tl X_+\wedge S^1)\simeq\omega^{mT-1}_{\Z/2}(\tl X_+)$.
There is an equivariant map $\phi\:\tl X\x I\to\tl X_+\wedge S^1$ sending
$\tl X\x\partial I$ to the basepoint, where $\Z/2$ acts trivially on $I$.
So any $\alpha\:\tl X_+\wedge S^1\to S^{mT}$, precomposed with $\phi$, is
an equivariant homotopy between two constant maps $\tl X\to *\i S^{mT}$.
Given an embedding $g_0\:X\emb\R^m$, one of these constant maps lifts
to $\tl g_0\:\tl X\to S^{m-1}\i S^\infty$.
By Lemma 6.7(b) below, the null-homotopy $\alpha\phi$ of this constant
map corresponds to an equivariant homotopy $\tl X\x I\to S^\infty$ between
$\tl g_0$ and some map $G\:\tl X\to S^{m-1}$.
By the Haefliger--Weber Criterion 3.1(a), $G$ is equivariantly homotopic
within $S^{m-1}$ to $\tl g_\alpha$ for some embedding $g_\alpha\:X\emb\R^m$.
If $\alpha'$ is homotopic to $\alpha$, the Haefliger--Weber Criterion 3.1(b)
along with Lemma 6.7(b) imply that $g_{\alpha'}$ is isotopic to $g_\alpha$.

Conversely, given any embedding $g\:X\emb\R^m$, $\tl g_0$ and $\tl g$
are equivariantly homotopic with values in $S^\infty$.
Projecting this homotopy to $S^{mT}$ we get a map $\tl X\x I\to S^{mT}$
that factors into the composition of $\phi$ and a map
$\alpha_g\:\tl X_+\wedge S^1\to S^{mT}$.
If $g'$ is isotopic to $g$, then $\tl g'$ and $\tl g$ are equivariantly
homotopic with values in $S^{m-1}$, therefore $\alpha_{g'}$ coincides with
$\alpha_g$. \qed
\enddemo

If $f,g\:K\to S^{m-1}_-$ are equivariant maps of an $k$-polyhedron with
a free $\Z/2$ action, they are related by an equivariant homotopy
$\phi_{f,g}\:K\x I\to S^\infty_-$, which is unique up to equivariant
homotopy $\rel K\x\partial I$.

\proclaim{Lemma 6.7} Let $\Z/2$ act freely on a $k$-polyhedron $K$ and
trivially on $I=[0,1]$.

(a) Suppose $k\le 2m-4$.
Equivariant maps $f,g\:K\to S^{m-1}_-$ are equivariantly homotopic iff
$\rho^\infty_m\phi_{f,g}\:K\x I\to S^{mT}$ is equivariantly null-homotopic
$\rel K\x\partial I$.

(b) Suppose $k\le 2m-3$.
For any equivariant map $f\:K\to S^{m-1}_-$ and any equivariant
self-homotopy $H\:K\x I\to S^{mT}$ of the constant map $K\to *\i S^{mT}$
there exists an equivariant map $g\:K\to S^{m-1}_-$ such that
$\rho^\infty_m\phi_{f,g}$ and $H$ are equivariantly homotopic
$\rel K\x\partial I$.
\endproclaim

The proof is similar to that of Lemma 6.4.

\proclaim{Theorem 6.8} Let $X$ be a compact $n$-polyhedron, and let
$S=N\but\Delta_X$, where $N$ is the second derived neighborhood of $\Delta_X$
in a $\Z/2$-invariant triangulation of $X\x X$.

If $m\ge\frac{3(n+1)}2$, then $X$ immerses in $\R^m$ if and only if
the image of $\Theta^m(X)$ in $\omega^{mT}_{\Z/2}(S_+)$ is trivial.

If $m>\frac{3(n+1)}2$ and $X$ immerses in $\R^m$, then the set of regular
isotopy classes of immersions of $X$ into $\R^m$ is in bijection with
$\omega^{mT-1}_{\Z/2}(S_+)$.
\endproclaim

The proof repeats those of Theorems 6.3 and 6.6, if instead of the
Haefliger--Weber Criterion 3.1 one uses its analogue for immersions
\cite{Har; Theorem 2 and footnote on p.\ 3}.
If $X$ is a closed smooth manifold, the dimensional restrictions in Theorem 6.8
can be weakened to $m\ge\frac{3n+1}2$ and $m>\frac{3n+1}2$, respectively,
by using \cite{HH} instead of \cite{Har}.

\remark{Remark}
If $X$ is a closed smooth $n$-manifold and $S$ is as in Theorem 6.8, then
$$\secat_{\Z/2}(S)=\secat_{\Z/2}(\tl X)-1=n+\nu(X),$$
where $\nu(X)$ is the largest number $k$ such that the normal Stiefel--Whitney
class $\bar w_k(X)$ does not vanish \cite{CF; 6.6} (compare \cite{Mc1},
\cite{De}), \cite{Wu}, \cite{Wu'}, \cite{Sch; Ch.\ VII, \S\S2-3}.
On the other hand, by a well-known result of Massey, $\nu(X)\le n-\alpha(n)$,
where $\alpha(n)$ is the number of ones in the dyadic expansion of $n$.
In view of Theorems 6.3 and 6.8, it would be of great interest to find
extraordinary (equivariant stable cohomotopy) versions of these results if
they exist.
(Note that the case of unoriented cobordism should be easy to check in view of
the literature on tom Dieck operations and Conner--Floyd characteristic
classes, see in particular \cite{Mc1}, \cite{BRS}, \cite{De}.)
A well-known old conjecture says that $X$ always embeds in $\R^{2n-\alpha(n)+1}$
and immerses in $\R^{2n-\alpha(n)}$; a difficult proof of the immersion
conjecture has been published by R. Cohen in mid-80s, while the embedding
conjecture remains open.
\endremark

\proclaim{Problem 6.9} Can the geometric interpretation of the extraordinary
van Kampen obstruction be used to obtain a polyhedral version of Haefliger's
generalized Whitney trick, thus proving Criterion 3.1 without induction
on simplices?
\endproclaim

\head 7. Embeddability versus disjoinability \endhead

If $f\:X\to\R^m$ is a (continuous) map, we write
$\Delta_f=\{(x,y)\in\tl X\mid f(x)=f(y)\}$ and
$\Delta_f^\eps=\{(x,y)\in\tl X\mid ||f(x)-f(y)||<\eps\}$.

\proclaim{Approximability Criterion 7.1} Let $X$ be a compact $n$-polyhedron,
suppose that $m\ge\frac{3(n+1)}2$, and let $f\:X\to\R^m$ a continuous map.

(a) $f$ is $C^0$-approximable by embeddings iff $\check\Theta(f)=0\in
\invlim\omega^{mT}_{\Z/2}(\tl X_+,\tl X_+\but\Delta_f^\eps)$.

(b) For each $\eps>0$ there exists a $\delta>0$ such that $f$ is
$C^0$-$\eps$-approximable by an embedding if
$\Theta_\delta(f)=0\in\omega^{mT}_{\Z/2}(\tl X_+,\tl X_+\but\Delta_f^\delta)$
and only if $\Theta_\eps(f)=0$.
\endproclaim

The definition of the controlled extraordinary van Kampen obstruction
$\Theta_\eps(f)$ and the proof of (b) is analogous to \cite{M2; \S3}
(see also \cite{Ah}), taking into account the controlled version of
the Haefliger--Weber Criterion \cite{RS}.
$\check\Theta(f)$ is the thread of the elements $\Theta_\eps(f)$ in
the inverse limit (compare \cite{M1}), so (a) is a special case of (b).
Taking $f$ to be a constant map, we recover Theorem 6.3.

\remark{Remark} Let $x$ be an isolated double point $f(p)=f(q)$ of
a generic PL map $f$ of an $n$-manifold into $\R^m$ in the metastable range.
At this double point, $f$ looks locally like the cone on a link
$L_x\:S^{n-1}\sqcup S^{n-1}\emb S^{m-1}$.
When $m<2n$, such a link may be non-trivial, even though each component
is null-homologous in the complement to the other one.
Therefore $f$ may fail to be $C^0$-approximable by embeddings, even though
the (first) cohomological obstruction to such an approximation, that is
the Hurewicz image of $\Theta(f)$, vanishes (see \cite{M1; Example 1}).
Indeed the component $\Theta_x(f)$ of $\Theta(f)$ in the direct summand
$\pi_{2n}(S^m)\simeq
\omega^{mT}_{\Z/2}((D^{2n}\x S^0_-)_+,(\partial D^{2n}\x S^0_-)_+)$
corresponding to $x$ is nothing but the $\alpha$-invariant of the link map $L_x$.
Indeed, $\Theta_x(f)$ is the homotopy class of the composition
$(D^{2n},\partial D^{2n})\to (S^\infty, S^{m-1})\to (S^m,*)$
which is the double suspension of the map $S^{n-1}\wedge S^{n-1}\to S^{m-2}$
representing $\alpha(L_x)$, cf\. \cite{MR}.
\endremark

\proclaim{Theorem 7.2} If $m\ge\frac{3(n+1)}2$, any map of a compact PL
$n$-manifold into $\R^m$ with countable singular set is $C^0$-approximable
by PL maps with finite singular sets.
\endproclaim

The original proof was significantly simplified by P. M. Akhmetiev, who
explained to the author how to eliminate induction by countable infinite
ordinals from that argument.

\demo{Proof} Without loss of generality $m<2n$.
Let $f\:N\to\R^m$ be the map and $\eps>0$ the desired closeness.
Let $\delta=\delta(\eps)$ be given by Criterion 7.1(b).
Since $\Delta_f$ is contained in $S_f\x S_f$, it is countable, and
since it is the preimage of $\Delta_{\R^m}$, it is closed in $\tilde X$.
Let $DN$ be an equivariant closed tubular $\delta$-neighborhood of the infinity
(i.e., of the removed diagonal) in $\tl X$.
By either proof of Lemma 7.3 below, we may assume that $\partial DN$ is
disjoint from $\Delta_f$.
By Lemma 7.3, $(\Delta_f\but DN)/t$ is contained in the interior of
a disjoint union $U/t$ of finitely many small PL balls $B_i$.
For each $i$, let $\{x_i,y_i\}$ be some point in $(\Delta_f/t)\cap B_i$.
Let $X$ be the quotient of $N$ by the equivalence relation $x_i\sim y_i$
for all $i$, so that $f$ factors as $N@>p>>X@>g>>\R^m$.
Then $\tl X$ is obtained from $\tl N$ by removing the set $P$ of all points
$(x_i,y_i)$ and $(y_i,x_i)$ and identifying each $\{x_i\}\x (N\but\{y_i\})$
with $\{y_i\}\x (N\but\{x_i\})$, as well as the symmetric sets, via a map
$q\:\tl N\but P\to\tl X$.
The cone of $q$ collapses onto an $(n+1)$-polyhedron, hence $q$ induces
an isomorphism $\omega^{mT}_{\Z/2}(\tl N_+\but P,\tl N_+\but (DN\cup U))\to
\omega^{mT}_{\Z/2}(\tl X_+,\tl X_+\but q(DN\cup U\but P))$ since $m>n+1$.
But the former group is zero.
Hence $\Theta_\delta(g)=0$ so $g$ is $\eps$-approximable by an embedding.
Thus $f$ is $\eps$-approximable by a map with finitely many double points.
\qed
\enddemo

\proclaim{Lemma 7.3} A countable closed set $S$ in a PL manifold $M$ is
contained in the interior of a disjoint union of arbitrarily small PL
balls in $M$.
\endproclaim

\demo{First proof} Let $K^{(m-1)}$ be the codimension one dual skeleton of
a sufficiently fine triangulation $T$ of $M^m$.
The interior of each simplex of $T$ of dimension $>0$ contains precisely one
vertex of $K^{(m-1)}$, and the position of $K^{(m-1)}$ is entirely determined
by the positions of its vertices in the open simplices of $T$.
Since each vertex can assume uncountably many positions, it follows by an
easy induction on the dual skeleta $K^{(i)}$ that $K^{(m-1)}$ can be made
disjoint from $S$ by an appropriate perturbation of these positions.
Then the stars of vertices of $T$ in this new derived subdivision of $T$
contain $S$ in their interiors.
They still do so after subtracting sufficiently small collars of
the boundary, which makes them disjoint. \qed
\enddemo

\demo{Second proof} Since every perfect set is uncountable,
$S\supset S'\supset\dots\supset S^{(\alpha)}=\emptyset$ for some
countable ordinal $\alpha$.
Here $S^{(\kappa+1)}=\left[S^{(\kappa)}\right]'$  is the Cantor--Bendixon
derivative, i.e\. the set of all limit points of $S^{(\kappa)}$ and,
whenever $\lambda$ is not of the form $\kappa+1$, $S^{(\lambda)}$ is
the intersection of $S^{(\kappa)}$'s for all $\kappa<\lambda$.
Since $S$ is $0$-dimensional, it is contained in a disjoint union  of compact
PL manifolds, so without loss of generality $M$ is compact.
We may also assume that $S^{(\kappa)}\ne\emptyset$ for all $\kappa<\alpha$.
Then by Cantor's theorem (every decreasing chain of compact sets has
non-empty intersection), $\alpha$ is of the form $\beta+1$ for some $\beta$.
Since $S^{(\beta)}$ is finite, it is contained in a disjoint union $V\i M$ of
arbitrarily small PL balls.
On the other hand, $T=S\but S^{(\beta)}$ satisfies $T^{(\beta)}=\emptyset$,
hence by the induction hypothesis is contained in a disjoint union
$U\i M\but S^{(\beta)}$ of arbitrarily small PL balls.
Let $\phi\:M\to Q$ be the quotient of $M$ by each of the balls $U_0\i U$
that meet $\partial V$.
Then $\phi$ can be arbitrarily closely approximated by a homeomorphism $h$,
and we may assume that $h(\partial V)$ is disjoint from both the finite set
$\phi(U_0)$ and the sufficiently distant set $\phi(U\but U_0)$.
Thus $U\cup\phi^{-1}(h(V))$ is a disjoint union of PL balls containing $S$.
\qed
\enddemo

\proclaim{Theorem 7.4} Suppose $m-n\ge 3$.
Any PL $n$-manifold that admits a map to $\R^m$ with finite singular
set embeds in $\R^m$.
\endproclaim

\demo{Proof} Let $f\:N^n\to\R^m$ be the map.
Then $f$ is the composition of a PL map $\phi$ onto a polyhedron $K^n$ and
a topological embedding $\psi\:K\emb\R^m$.
By the Chernavskij--Miller--Bryant Theorem \cite{B1} $\psi$ is
$C^0$-approximable by PL embeddings.
Let $\chi$ be one, and consider the PL map $g\:N@>\phi>>K@>\chi>>\R^m$.

What follows is a simplified version the Penrose--Whitehead--Zeeman trick.
We may assume without loss of generality that $n>1$ and that $N$ is
connected.
Then the inclusion $S_g\i N$ extends to a PL embedding of the cone
$i\:c*S_g\emb N$.
Since $m-n\ge 3$, by general position the inclusion of the polyhedron
$Q:=g(i(c*S_g))$ into $\R^m$ extends to a PL embedding of the cone
$j\:c*Q\emb\R^m$, with image disjoint from $g(N)\but Q$.
Let $D^n$ and $B^m$ be the second derived neighborhoods of $i(c*S_g)$
and $j(c*Q)$ in some triangulations of $N$ and $\R^m$ where $g$ is
simplicial.
Then $g^{-1}(Q)=P$, and since $i(c*S_g)$ and $j(c*Q)$ are collapsible,
$D$ and $B$ are balls.
Viewing them as cones on their boundaries, define $h\:N\to\R^m$ to coincide
with $g$ outside $D$ and to conewise extend
$g|_{\partial D}\:\partial D\to\partial B$ on $D$.
Since the latter map is an embedding, so is $h$. \qed
\enddemo

In conclusion, we show that the method used in the proof of Theorem 7.4
(the Penrose--Whitehead--Zeeman trick) gives an alternative proof for
a weak form of Theorem 7.2.
An advantage of this approach is that it extends to codimension $3$.

\proclaim{Theorem 7.5} If $m-n\ge 3$, any map of a compact PL $n$-manifold
into a PL $m$-manifold with countable closure of the singular set is
$C^0$-approximable by PL maps with finite singular sets.
\endproclaim

\demo{Proof} Let $f\:N\to M$ be the given map and $\eps>0$ the desired
closeness.
The open manifold $N\but\Cl(S_f)$ can be represented as the union of
a sequence $K_0\i K_1\i K_2\i\dots$ of compact manifolds with boundary.
By the Chernavskij--Miller Theorem (see \cite{B2}) $f$ is
$\frac\eps4$-homotopic with support in $K_2$ to a map which still embeds
$K_2$ and is PL on $K_1$.
This map is $\frac\eps8$-homotopic with support in $K_3\but K_0$ to a map
which still embeds $K_3$ and is PL on $K_2$.
Proceeding in this fashion, in the end%
\footnote{An extra care about the epsilonics in what follows would allow us
to do with just one application of the Chernavskij--Miller Theorem here.}
we will obtain that $f$ is $\frac\eps2$-homotopic $\rel\Cl(S_f)$ to a map
still denoted $f$, which PL embeds $N\but\Cl(S_f)$.

By Lemma 7.3, $f(\Cl(S_f))$ is contained in the interior of a disjoint
union $P_k$ of PL balls of diameters $<\eps/2$ in $M$.
Then $L:=f^{-1}(P_k)$ is a compact PL manifold.
By Lemma 7.3, the intersection of $\Cl(S_f)$ with each connected component
of $L$ is contained in a PL ball in this component.
Let $J_k$ be the union of all these PL balls.
Then $Z:=f(\Cl(N\but J_k))$ is disjoint from $f(\Cl(S_f))$.
By Lemma 7.3, $f(\Cl(S_f))$ is contained in the interior of a disjoint
union $P_{k-1}$ of PL balls in $P_k\but Z$.
By construction, $f^{-1}(P_{k-1})\i J_k$.
Proceeding in this fashion, we obtain $J_i$'s and $P_i$'s as in Lemma 7.6.
Finally, if we choose them with something to spare (i.e.\ leaving sufficient
margins), their inclusion properties will preserve for a sufficiently close
PL approximation of $f$. \qed
\enddemo

\proclaim{Lemma 7.6} Let $f\:N^n\to M^m$ be a PL map between PL manifolds,
$m-n\ge 3$, and let $k$ be such that $m\ge(1+\frac1k)(n+1)$.
Suppose that $S_f\i J_0\i\dots\i J_k\i M$, where each $J_k$ is a disjoint
union of PL $n$-balls and each $f(J_i)$ is contained in a PL $m$-ball
$P_i$ such that for $i<k$, $f^{-1}(P_i)\i J_{i+1}$.
Then $f$ is homotopic with support in $J_k$ and with image of the support
in $P_k$ to a PL map $g$ with finite $S_g$.
\endproclaim

\demo{Proof} The Penrose--Whitehead--Zeeman--Irwin trick \cite{Z}. \qed
\enddemo

\remark{Remark} If $f$ is a map of a compact PL $n$-manifold into a PL
$m$-manifold, $m-n\ge 3$, such that both $\Cl(S_f)$ and $f(\Cl(S_f))$ are
$k$-dimensional and {\it tame} in the sense of Shtan'ko, then $f$ is
approximable by PL maps whose singular sets are $k$-polyhedra.
This can be proved using the Homma--Bryant argument \cite{B2} for
the Chernavskij--Miller approximation theorem.
We omit the details.
\endremark

\proclaim{Problem 7.7} Let $f$ be a codimension $\ge 3$ continuous map between
compact PL manifolds with a $k$-dimensional $S_f$.
Is $f$ approximable by PL maps $g$ with $k$-dimensional $S_g$?
\endproclaim

One might also wonder whether there exists a map $f\:S^3\to\R^5$ with
$S_f$ the Antoine necklace, which is non-approximable by maps with finite
singular sets.

It should be mentioned here that every map $f$ from a closed $2$-manifold
to a closed $3$-manifold with $0$-dimensional $\Cl(S(f))$ is approximable
by embeddings \cite{Bra}.
It has been conjectured in mid-80s that this is the case more generally for
every map $f$ from a closed $(m-1)$-manifold to a closed $m$-manifold, $m\ge 3$,
with $0$-dimensional $S(f)$ (see \cite{RR+}).

\head 8. Polyhedra whose subsets satisfy partial Alexander duality \endhead

\proclaim{Lemma 8.1} For every compact polyhedron $X$ with
$H^i(\tl X)=0$ for $i\ge m$, there exists an equivariant map
$\tl X\to S^{m-1}$.
\endproclaim

\demo{Proof (revised in v5)\footnote{A. Skopenkov pointed out that the published version of this
proof was insufficient.
That argument, which only showed that
$H^{i+1}(\bar X;\,\pi_i(S^{m-1}))=0=H^{i+1}(\bar
X;\,\pi_i(S^{m-1})\otimes\Z_p)$,
is in fact well-known \cite{HH; p.\ 237}, \cite{Ad; Lemma 7.3} and, as observed
in \cite{HH; p.\ 236}, suffices to prove Lemma 8.1 under the (harmless)
additional hypothesis $\dim X<m$.
Indeed, the antipodal action on $\pi_i(S^{m-1})$ is trivial when $m-1$ is odd
(because $[-\id]=0$ in this case), and coincides with multiplication by $-1$
when $m-1$ is even and the suspension homomorphism
$\pi_{i-1}(S^{m-2})\to\pi_i(S^{m-1})$ is onto (because $[\id]+[-\id]=0$ in this
case, and it is easy to see that $[g(\Sigma f)+h(\Sigma f)]=[(g+h)(\Sigma f)]$).
An alternative proof of Lemma 8.1 (without the additional hypothesis) was found
in \cite{GS} using the generalized Smith sequence (see Remark 2.3).}}
Let $G$ be a finitely generated abelian group.
By the universal coefficients formula, $H^i(\tl X;\,G)=0$ for all $i\ge m$.
It follows from the two Smith sequences, by downward induction on $i$, that
$H^i(\bar X;\,G)=0=H^i(\bar X;\,G\otimes\Z_p)$ for all $i\ge m$,
where $\Z_p$ is the integral local coefficient system associated with
the double covering $p\:\tl X\to\bar X$.

Let $M$ be a finitely generated $\Lambda$-module, where $\Lambda=\Z[\Z/2]$, and
let $I$ be the augmentation ideal of $\Lambda$.
By a trivial case of \cite{CR; Theorem 74.3}, which is easy to verify directly,
$M\simeq G_1\oplus (G_2\otimes_\Z I)\oplus (G_3\otimes_\Z\Lambda)$
for some finitely generated abelian groups $G_1$, $G_2$ and $G_3$.
Since $H^i_{\Z/2}(\tl X;\,G_1)\simeq H^i(\bar X;\,G_1)$,
$H^i_{\Z/2}(\tl X;\,G_2\otimes I)\simeq H^i(\bar X;\,G_2\otimes\Z_p)$ and
$H^i_{\Z/2}(\tl X;\,G_3\otimes\Lambda)\simeq H^i(\tl X;\,G_3)$
(see \S2), we conclude that $H^i_{\Z/2}(\tl X;\,M)=0$ for all $i\ge m$.

Let $M_i$ be $\pi_i(S^{m-1})$ regarded as a $\Lambda$-module under the
action induced by the antipodal involution on $S^{m-1}$, and let $\Cal F_{M_i}$
be the corresponding local coefficient system on $\bar X$.
Then $H^{i+1}(\bar X;\,\Cal F_{M_i})\simeq H^{i+1}_{\Z/2}(\tl X;\,M_i)=0$
for all $i\ge m-1$.
The assertion now follows by the standard non-homotopically-simple obstruction
theory (as e.g.\ in the Hilton--Wylie textbook). \qed
\enddemo

\proclaim{Theorem 8.2} Every $n$-polyhedron $X$ with $H^{n-d}(X\but x)=0$
for each $x\in X$ and $d\le k$, where $k<\frac{n-3}2$, embeds in $\R^{2n-k}$.
\endproclaim

\demo{Proof}
Let $\H^{n-d}(\pi)$ be the Leray sheaf of the projection
$\pi\:\tl X_s\i X\x X\to X$ \cite{Bre}.
Since the simplicial deleted product $\tl X_s$ is compact, $\pi$ is closed,
so the stalks $\H^{n-d}(\pi)_x\simeq H^{n-d}(X\but x)$.
Consider the Leray spectral sequence \cite{Bre}
$$E_2^{pq}=H^p(X;\H^q(\pi))\Rightarrow H^{p+q}(\tl X_s).$$
Since $\H^{n-d}(\pi)=0$ for all $d\le k$, we get that $H^{2n-d}(\tl X_s)=0$
for all $d\le k$. \qed
\enddemo

V. M. Buchstaber asked the author in September 2005, whether there is
a generalization to polyhedra of the classical Penrose--Whitehead--Zeeman
Theorem that $k$-connected PL $n$-manifolds embed in $\R^{2n-k}$ in
the metastable range.
(For a proof of this theorem see \cite{Z}, which includes Irwin's extension
to codimension three.%
\footnote{Added in v5: Another proof, which also works for homologically $k$-connected 
manifolds, is given by the Haefliger--Weber Criterion 3.1 along with Lemma 8.1 and
the Poincare duality $H^i(\tl X)\simeq H_{2n-i}(X\x X,\Delta_X)$; cf.\
\cite{Hae; p.\ 66}, \cite{We; Th\'eor\`eme 4}, \cite{Ad}.}
The author learned from A. Skopenkov that the extension of the
Penrose--Whitehead--Zeeman Theorem to {\it homologically} $k$-connected
manifolds has been long known to A. Haefliger and can be proved by the methods
of \cite{GS}, but apparently has not been explicitly stated in
the literature.)

One answer to Buchstaber's question is already given by Theorem 8.2.
A minimal set of purely local and purely global conditions implying
the hypothesis of Theorem 8.2 is as follows:

\roster
\item"(i)$_k$\ \ \ \ " $H^{n-d}(X)\simeq H_d(pt)$ for $d\le k$, and

\smallskip
\item"(ii)$_{k-1}$" $H^{n-d}(X,X\but x)\simeq H_d(pt)$ for $d\le k-1$ and each
$x\in X$.
\endroster
\medskip

\proclaim{Theorem 8.3} Let $X$ be a compact $n$-polyhedron.
\smallskip

(a) The following are equivalent:

\roster
\item"$\bullet$" conditions (i)$_k$ and (ii)$_{k-1}$;

\item"$\bullet$" $\tl H_i(X\but Y)\simeq H^{n-i-1}(Y)$ for $i\le k-1$ and
all $Y\i X$.
\endroster
\smallskip

(b) If $X$ satisfies conditions (i)$_0$ and (ii)$_0$, the following
are equivalent:

\roster
\item"$\bullet$" conditions (i)$_k$ and (ii)$_k$;

\item"$\bullet$"
$H_i(X)=0$ for $i\le k$, and $H_i(X,X\but Y)=0$ for $i\le k+1$ and all
$Y\i X$ with $\dim Y<n-i$;

\item"$\bullet$"
$H_i(X)=0$ for $i\le k$, and $X\but Z$ is a homology manifold for some
$(n-k-1)$-polyhedron $Z\i X$ such that $H_i(X,X\but Y)=0$ for $i\le k+1$ and
all $Y\i Z$.
\endroster
\endproclaim

Notice that $H_i(X,X\but Y)=0$, where $\dim Y<n-i$, may be viewed as
a transversality-type condition: it says that $i$-pseudo-manifolds in $X$ can be
made disjoint, by pseudo-bordism relative to the boundary, from the subpolyhedron
$Y$ of codimension $>i$ (compare \cite{Mc3}, \cite{RoS; Theorem 4.2}).

At the same time, the original Penrose--Whitehead--Zeeman argument can be seen
to work to embed into $\R^{2n-k}$, $k<\frac{n-3}2$, every $k$-connected compact
$n$-polyhedron $X$ that is a manifold away from a subpolyhedron $Z$ of
dimension $<n-k-1$ such that $\pi_i(X,X\but Z)=0$ for all $i\le k+1$
(see \cite{Z}).

In the case $k=0$, Sarkaria \cite{Sa1} noticed that the PWZ method also works
to embed in $\R^{2n}$ every quotient $X$ of an $n$-manifold $M$, $n>2$, by
an identification on the boundary such that no two components of $M$ remain
disjoint in $X$.
Clearly, this includes all $n$-polyhedra satisfying the hypothesis of
Theorem 8.2.

\medskip
We now proceed to the proof of Theorem 8.3.
It will be convenient to work in a slightly greater generality.
To this end, the {\it orientation sheaf} $\H_n(X)$ of an $n$-polyhedron $X$
is defined by $U\mapsto H_n(X,X\but U)$ and has stalks
$\H_n(X)_x=H_n(X,X\but x)$.
The {\it fundamental cosheaf} $\H^n(X)$ is defined by
$U\mapsto H^n(X,X\but U)$ and has stalks $\H^n(X)_x=H^n(X,X\but x)$.
(It is an obvious feature of the top dimension that this presheaf, resp\.
precosheaf is indeed a sheaf, resp\. cosheaf.)
If $H^n(X,X\but x)\simeq\Z$ for each $x\in X$, the fundamental cosheaf is
clearly locally constant.
By the universal coefficients formula, the orientation sheaf
$\H_n(X)=\Hom(\H^n(X),\Z)$.
The notions of a cosheaf and cosheaf homology are defined by reversing all
arrows in those of a sheaf and sheaf cohomology, cf\. \cite{Bre}.

\proclaim{Proposition 8.4} Let $X$ be a compact $n$-polyhedron and $d$
a non-negative integer.
The following assertions are equivalent:

\roster
\item $H^{n-i}(X,X\but x)=0$ for all $x\in X$ and $1\le i\le d$;
\smallskip

\item $H_i(X\but Z,X\but Y;\,\H^n(X))$ is isomorphic to $H^{n-i}(Y,Z)$ when
$i\le d$ and to its subgroup when $i=d+1$, for all $Z\i X$.
\smallskip

\item $H_i(X,X\but Y;\,\H^n(X))=0$ for $i\le d+1$ and all $Y\i X$ with
$\dim Y<n-i$.
\endroster
\endproclaim

The proof also shows that (2) is equivalent to its versions with $Z=\emptyset$
and with $Y=X$.

\demo{Proof} Clearly, (2) implies (1).
Conversely, the Zeeman spectral sequence runs
$$E^2_{p,q}=H_q(X\but Z,X\but Y;\,\H^p(X))\imp H^{p-q}(Y,Z)$$
(see \cite{Mc3} and \cite{Bre; \S VI.14 with $f=\id$}).
Here $\H^p(X)$ is the cosheaf generated \cite{Bre} by the precosheaf
$U\mapsto H^p(X,X\but U)$.
The hypothesis implies that $E^\infty_{p,q}=E^2_{p,q}=0$ for $p-q\ge n-d$,
$p\ne n$.
So the edge homomorphism
$e_i\:H_i(X\but Z,X\but Y;\H^n(X))\to H^{n-i}(Y,Z)$ is an isomorphism
for each $i\le d$.
Moreover, all differentials $E^r_{p,q}\to E^r_{p+r-1,q+r}$ must be trivial
for $p-q\ge n-d$ (cf\. Figure 1 in \cite{Mc3}).
These include all differentials with target group $E^r_{p,q}$ where
$p-q=n-d-1$, hence $E^\infty_{p,q}=E^2_{p,q}$ for such $p,q$.
Thus $e_{d+1}$ is injective.

Clearly, (2) implies (3), so to complete the proof it suffices to show that
(3) implies (1).
Let us fix some triangulation of $X$.
Then (1) is equivalent to the assertion
\roster
\item"($*$)" The link $L$ of every simplex $A^a$ of $X$ has trivial reduced
cohomology groups in dimensions $l-1,\dots,l-d$, where $l=\dim L=n-a-1$.
\endroster
Assume inductively that ($*$) holds for all simplices of dimension $>a$.
The link in $L$ of a simplex $B^b$ of $L$ is the link in $X$ of the
simplex $A*B$ of $X$.
So $L$ itself satisfies condition (1).
Hence by the above it also satisfies (2), in particular
$H_i(L;\,\H^l(L))\simeq H^{n-i}(L)$ for $i\le d$.
By Lemma 8.5 below,
$H_i(L;\,\H^l(L))\simeq H_{i+1}(X\but\partial A,X\but A;\,\H^n(X))$
for $i>0$.
If $a<n-i-1$ and $i+1\le d+1$, the latter group is zero
since $H_{i+1}(X,X\but A;\,\H^n(X))=H_i(X,X\but\partial A;\,\H^n(X))=0$
by (3). \qed
\enddemo

\proclaim{Lemma 8.5} If $L^l$ is the link of a simplex $A$ in a fixed
triangulation of an $n$-polyhedron $X$, then
$H_i(L;\,\H^l(L))\simeq H_{i+1}(X\but\partial A,X\but A;\,\H^n(X))$
for $i>0$.
\endproclaim

This is obvious when $\H^n(X)$ is locally constant, but we need to be more
careful in the general case (which is not used anywhere except Proposition
8.4).

\demo{Proof} Let $S$ be the open star $c*L\but L$ of $A$.
Since $S\but c$ is homeomorphic to $L\x\R^1$, we have
$H_i(L;\,\H^l(L))\simeq H_i(S\but c;\,\H^{l+1}(S\but c))$.
Now $H_i(S;\H^{l+1}(S))$ is isomorphic to
$H_i(c*L;\H^{l+1}(cL\cup_L L\x I)|_{c*L})$.
The chain complex $C_*(c*L;\,\H^{l+1}(cL\cup_L L\x I)|_{c*L})$ is, apart from
an additional summand in dimension $0$, the cone of the identity chain map on
$C_*(L;\H^l(L))$.
It follows that $H_i(S;H^{l+1}(S))$ vanishes for $i>0$, and therefore
$H_i(S\but c;\,\H^{l+1}(S\but c))\simeq H_{i+1}(S,S\but c;\,\H^{l+1}(S))$
for $i>0$.

Let $V=A\but\partial A$ and $U=V*L\but L$.
By considering the projection $p\:U\to S$, we get
$H_{i+1}(S,S\but c;\,\H^{l+1}(S))\simeq H_{i+1}(U,U\but V;\,\H^n(U))$ due to
$\H^n(U)=p^*\H^{l+1}(S)$.
Finally $H_{i+1}(U,U\but V;\,\H^n(V))
\simeq H_{i+1}(X\but\partial A,X\but A;\,\H^n(X))$ by excision. \qed
\enddemo

\remark{Remark} The proof of Proposition 8.4 works to show that the obvious
dual homological conditions (1$'$), (2$'$) and (3$'$) are mutually equivalent
as well.
For example, (2$'$) says:

\roster
\item"(2$'$)"
$H^i(X\but Z,X\but Y;\,\H_n(X))$ is isomorphic to $H_{n-i}(Y,Z)$ when
$i\le d$ and to its subgroup when $i=d+1$, for all closed subpolyhedra
$Z\i Y\i X$.
\endroster

When $\H^n(X)$ is locally constant, $H_n(X;\,\H^n(X))$ is isomorphic to $\Z$
by (1$'$)$\iff$(2$'$) with $d=0$, and the iso/monomorphism in (2$'$) can be
identified (see \cite{Mc3; 6.3}) as the cap product with a generator of this
group.

It would be interesting to know if conditions (i) and (ii) can also be
rethought from the viewpoint of the slant product (see \cite{Mc2}).
\endremark

\proclaim{Lemma 8.6} An $n$-polyhedron $X$ satisfying condition (ii)$_{k-1}$ is
a homology manifold away from a subpolyhedron of dimension at most
$n-(2k+1)$.
\endproclaim

\demo{Proof} By the proof of (3)$\imp$(1) in Proposition 8.4, the link $L$
of an $(n-i-1)$-simplex $\sigma$ of $X$ has its top $k$ cohomology (and hence
also homology) groups, as well as the first $k$ homology groups isomorphic
to those of $S^i$.
Hence if $i<2k$, $L$ has the same homology as the $i$-sphere, and so $X$ is
a homology manifold at the interior points of $\sigma^{n-i-1}$. \qed
\enddemo

Theorem 8.3 follows immediately from Proposition 8.4 and Lemma 8.6.


\Refs\widestnumber\key{MTW}

\ref \key Ad \by M. Adachi
\book Umekomi to Hamekomi \publ Iwanami Shoten \publaddr Tokyo \yr 1984
\transl English transl.
\book Embeddings and Immersions
\bookinfo Trasl. Math. Mono. \vol 124 \publ Amer. Math. Soc. \publaddr
Providence, RI \yr 1993
\endref

\ref \key Ada \by J. F. Adams
\paper Prerequisites (on equivariant stable homotopy) for Carlsson's lecture
\inbook Algebraic Topology (Aarhus, 1982) 
\bookinfo Lecture Notes in Math. \vol 1051
\publ Springer \publaddr Berlin \yr 1984 \pages 483--532
\endref

\ref \key Ah \by P. M. Akhmetiev
\paper Pontrjagin--Thom construction for approximation of mappings by
embeddings \jour Topol. Appl. \vol 140 \yr 2004 \pages 133--149
\endref

\ref \key AE \by P. M. Akhmet'ev, P. J. Eccles
\paper The relationship between framed bordism and skew-framed bordism
\jour Bull. London Math. Soc. \vol 39 \yr 2007 \pages 473--481
\endref

\ref \key BF \by K. Barnett, M. Farber
\paper Topology of configuration space of two particles on a graph, I
\jour arXiv:0903.2180
\endref

\ref \key Bau \by D. R. Bausum
\paper Embeddings and immersions of manifolds in Euclidean space
\jour Trans. Amer. Math. Soc. \vol 213 \yr 1975 \pages 263--303
\endref

\ref \key Be \by J. C. Becker
\paper A relation between equivariant and non-equivariant stable cohomotopy
\jour Math. Z. \vol 199 \yr 1988 \pages 331--356
\endref

\ref \key BKK \by M. Bestvina, M. Kapovich, B. Kleiner
\paper Van Kampen's embedding obstruction for discrete groups
\jour Invent. Math. \vol 150 \yr 2002 \pages 219--235
\endref

\ref \key B\"o \by T. B\"ohme
\paper On spatial representations of graphs
\inbook Contemporary Methods in Graph Theory \ed R. Bodendieck \publ Mannheim
\publaddr Wien, Zurich \yr 1990 \pages 151--167
\endref

\ref \key Bra \by M. V. Brahm
\paper Approximating maps of $2$-manifolds with zero-dimensional
non-degeneracy sets
\jour Topol. Appl. \vol 45 \yr 1992 \pages 25--38
\endref

\ref \key Bre \by G. E. Bredon
\book Sheaf Theory (2nd edition)
\publ Springer \yr 1997
\endref

\ref \key Bro \by K. S. Brown
\book Cohomology of groups
\publ Springer-Verlag \publaddr New York \yr 1982
\endref

\ref \key B1 \by J. L. Bryant
\paper Approximating embeddings of polyhedra in codimension three
\jour Trans. Amer. Math. Soc. \vol 170 \yr 1972 \pages 85--95
\endref

\ref \key B2 \bysame
\paper Triangulation and general position of PL diagrams
\jour Topol. Appl. \vol 34 \pages 211--233 \yr 1990
\endref

\ref \key BRS \by S. Buoncristiano, C. P. Rourke, B. J. Sanderson
\book A Geometric Approach to Homology Theory
\bookinfo London Math. Soc. Lecture Note Ser. \vol 18
\publ Cambridge Univ. Press \yr 1976
\endref

\ref \key CG \by A. Cavicchioli, L. Grasselli
\paper Cohomological products and transversality
\jour Rend. Sem. Mat., Torino \vol 40 \issue 3 \pages 115--125 \yr 1982
\endref

\ref \key CL+ \by O. Cornea, G. Lupton, J. Oprea, D. Tenr\'e
\book Lusternik--Schnirelmann category
\bookinfo Math. Surv. and Mono. \vol 103 \publ Amer. Math. Soc. \yr 2003
\endref

\ref \key Cr \by M. C. Crabb
\book $\Z/2$-Homotopy Theory
\bookinfo London Math. Soc. LNS \vol 44 \publ Cambridge Univ. Press \yr 1980
\endref

\ref \key CF \by P. E. Conner, E. E. Floyd
\paper Fixed point free involutions and equivariant maps
\jour Bull. Amer. Math. Soc. \vol 66:6 \yr 1960 \pages 416--441
\endref

\ref \key CR \by C. W. Curtis, I. Reiner
\book Representation theory of finite groups and associative algebras
\publ Interscience \publaddr New York, London \yr 1962
\endref

\ref \key dL \by M. de Longueville
\paper Bier spheres and barycentric subdivision
\jour J. Combinat. Theory Ser. A \vol 105 \yr 2004 \pages 355--357
\endref

\ref \key Dax \by J.-P. Dax
\paper \'Etude homotopique des espaces de plongements
\jour Ann. scient. \'Ec. Norm. Sup., 4e s\'er. \vol 5 \yr 1972 \pages 303--377
\endref

\ref \key De \by M. Ded\`o
\paper Cobordism characteristic classes
\jour Ann. Mat. pura appl. \vol 117 \yr 1978 \pages 207--226
\endref

\ref \key EG \by P. J. Eccles, M. Grant
\paper Bordism groups of immersions and classes represented by
self-intersections
\jour Alg. Geom. Topology \vol 7 \yr 2007 \pages 1081--1097
\moreref arXiv:math.AT/0504152
\endref

\ref \key Fe \by R. Fenn
\book Techniques of Geometric Topology
\bookinfo London Math Soc. Lecture Note Ser. \vol 57
\publ Cambridge Univ. Press \yr 1983
\endref

\ref \key FS \by R. Fenn, D. Sjerve
\paper Geometric cohomology theory
\inbook Low-dimensional topology
\bookinfo Contemp. Math. \vol 20 \yr 1983 \pages 79--102
\endref

\ref \key FN \by T. Fleming, R. Nikkuni
\paper Homotopy on spatial graphs and the Sato--Levine invariant
\jour Trans. Amer. Math. Soc. \vol 361 \yr 2009 \pages  1885--1902
\moreref arXiv:math/0509003
\endref

\ref \key F1 \by A. Flores
\paper \"Uber die Existenz $n$-dimensionaler Komplexe, die nicht in den
$R_{2n}$ topologisch einbettbar sind
\jour Ergeb. Math. Kolloq. \yr 1933 \vol 5 \pages 17--24
\endref

\ref \key F2 \bysame
\paper \"Uber $n$-dimensionale Komplexe die im $R_{2n+1}$ absolut
selbstverschlungen sind
\jour Ergeb. Math. Kolloq. \yr 1935 \vol 6 \pages 4--7
\endref

\ref \key FKT \by M. H. Freedman, V. S. Krushkal, P. Teichner
\paper Van Kampen's embedding obstruction is incomplete for $2$-complexes
in $\R^4$ \jour Math. Res. Lett. \vol 1 \yr 1994 \pages 167--176
\endref

\ref \key GS \by D. Gon\c calves, A. Skopenkov
\paper Embeddings of homology equivalent manifolds with boundary
\jour Topol. Appl. \vol 153 \yr 2006 \pages 2026--2034
\endref

\ref \key Gra \by L. Grasselli
\paper Subdivision and Poincar\'e duality
\jour Riv. Mat. Univ. Parma, IV. \vol 9 \yr 1983 \pages 95--103
\endref

\ref \key Gr \by B. Gr\"unbaum
\paper Imbeddings of simplicial complexes
\jour Comm. Math. Helv. \vol 44 \yr 1969 \pages 502--513
\endref

\ref \key Hae \by A. Haefliger
\paper Plongements de vari\'et\'es dans le domaine stable
\jour Seminaire Bourbaki \vol 245 \yr 1962 \pages 63--77
\endref

\ref \key Har \by L. S. Harris
\paper Intersections and embeddings of polyhedra
\jour Topology \vol 8 \yr 1969 \pages 1--26
\endref

\ref \key HQ \by A. Hatcher, F. Quinn
\paper Bordism invariants of intersections of submanifolds
\jour Trans. Amer. Math. Soc. \vol 200 \yr 1974 \pages 327--344
\endref

\ref \key Hau \by H. Hauschild
\paper \"Aquivariante Homotopie I
\jour Arch. Math. \vol 29 \yr 1977 \pages 158--165
\endref

\ref \key HH \by A. Haefliger, M. W. Hirsch
\paper Immersions in the stable range
\jour Ann. Math. \vol 75 \yr 1962 \pages 231--241
\endref

\ref \key Hu \by Hu, Sze-Tsen
\paper Isotopy invariants of topological spaces
\jour Proc. Royal Soc. A \vol 255 \yr 1960 \pages 331--366
\endref

\ref \key Ja \by I. M. James
\paper On category, in the sense of Lusternik--Schnirelmann
\jour Topology \vol 17 \yr 1978 \pages 331--348
\endref

\ref \key Ki \by S. Kinoshita
\paper On elementary ideals of polyhedra in the 3-sphere
\jour Pacific J. Math. \vol 42 \yr 1972 \pages 89--98
\endref 

\ref \key Kl \by J. R. Klein
\paper On embeddings in the sphere
\jour Proc. Amer. Math. Soc. \vol 133 \yr 2005 \pages 2783--2793
\moreref arXiv:math.AT/0310236
\endref

\ref \key KW \bysame, B. Williams
\paper Homotopical intersection theory, II: Equivariance
\jour Math. Z. \toappear \moreref arXiv:0803.0017
\endref

\ref \key Kr \by V. S. Krushkal
\paper Embedding obstructions and $4$-dimensional thickenings of
$2$-complexes
\jour Proc. Amer. Math. Soc. \vol 128 \yr 2000 \pages 3683--3691
\moreref arXiv:math.GT/0004058
\endref

\ref \key Ko \by Cz. Kosniowski
\paper Equivariant cohomology and stable cohomotopy
\jour Math. Ann. \vol 210 \pages 83--104 \yr 1974
\endref

\ref \key LS \by L. Lov\'asz, A. Schrijver
\paper A Borsuk theorem for antipodal links and a spectral characterization
of linklessly embeddable graphs
\jour Proc. Amer. Math. Soc. \vol 126 \yr 1998 \pages 1275--1285
\endref

\ref \key Mc1 \by C. McCrory
\paper Cobordism operations and singularities of maps
\jour Bull. Amer. Math. Soc. \vol 82 \yr 1976 \pages 281--283
\endref

\ref \key Mc2 \bysame
\paper A characterization of homology manifolds
\jour J. London Math. Soc. \vol 16 \yr 1977 \pages 149--159
\endref

\ref \key Mc3 \bysame
\paper Zeeman's filtration of homology
\jour Trans. Amer. Math. Soc. \vol 250 \yr 1979 \pages 147--166
\endref

\ref \key Mc4 \bysame
\paper Geometric homology operations
\inbook Studies in Algebraic Topology \bookinfo Adv. Math. Suppl. Studies
\vol 5 \publ Academic Press \publaddr New York \yr 1979 \pages 119--141
\endref

\ref \key Man \by V. O. Manturov
\paper A proof of Vassiliev's conjecture on the planarity of singular links
\jour Izv. RAN, Ser. Mat. \vol 69 \issue 5 \yr 2005 \pages 169--178
\transl English transl. \jour Izv. Math. \vol 69\yr 2005 \pages 1025--1033
\endref

\ref \key MR \by W. S. Massey, D. Rolfsen
\paper Homotopy classification of higher dimensional links
\jour Indiana Univ. Math. J. \vol 34 \yr 1986 \pages 375--391
\endref

\ref \key Mat \by J. Matou\v{s}ek
\book Using the Borsuk--Ulam Theorem
\publ Springer \publaddr Berlin \yr 2003
\endref

\ref \key MTW \bysame, M. Tancer, U. Wagner
\paper Hardness of embedding simplicial complexes in $\R^d$
\jour arXiv: 0807.0336
\endref

\ref \key M+ \by J. P. May {\it et al.}
\book Equivariant Homotopy and Cohomology Theory
\bookinfo Regional Conference Series in Math. \vol 91
\publ Amer. Math. Soc. \publaddr Providence, RI \yr 1996
\endref

\ref \key M1 \by S. A. Melikhov
\paper Isotopic and continuous realizability of maps in the metastable range
\jour Mat. Sbornik \vol 195 \issue 7 \yr 2004 \pages 71--104
\transl English transl. \jour Sb. Math. \vol 195 \yr 2004 \pages 983--1016
\endref

\ref \key M1' \bysame
\paper On isotopic realizability of maps factored through a hyperplane
\jour Mat. Sbornik \vol 195 \issue 8 \yr 2004 \pages 47--90
\transl English transl. \jour Sb. Math. \vol 195 \yr 2004 \pages 1117--1163
\endref

\ref \key M2 \bysame
\paper Sphere eversions and realization of mappings
\jour Tr. Mat. Inst. Steklova \vol 247 \yr 2004 \pages 159--181
\transl English transl.
\jour Proc. Steklov Inst. Math. \vol 247 \yr 2004 \pages 1--20
\moreref arXiv: math.GT/0305158
\endref

\ref \key M3 \bysame
\paper Review of V. A. Vassiliev's paper ``First-order invariants and
first-order cohomology for spaces of embeddings of self-intersecting curves
in $\R^n$''
\jour Math. Reviews \issue 2179414 (2007i:57014) \pages 3pp
\endref

\ref \key MS \by S. A. Melikhov, E. V. Shchepin
\paper The telescope approach to embeddability of compacta
\jour arXiv:math.GT/0612085
\endref

\ref \key Mey \by D. M. Meyer
\paper $Z/p$-equivariant maps between lens spaces and spheres
\jour Math. Ann. \vol 312 \pages 197--214 \yr 1998
\endref

\ref \key Ni \by R. Nikkuni
\paper The second skew-symmetric cohomology group and spatial embeddings
of graphs
\jour J. Knot Theory Ram. \vol 9 \yr 2000 \pages 387--411
\endref

\ref \key RR+ \by D. Repov\v s, W. Rosicki, A. Zastrow, M \v Zeljko
\paper Constructing near-embeddings of codimension one manifolds with
countable dense singular sets
\jour arXiv:0803.4251
\endref

\ref \key RS \by D. Repov\v{s}, A. B. Skopenkov
\paper A deleted product criterion for approximability of maps by embeddings
\jour Topol. Appl. \vol 87 \yr 1998 \pages 1--19
\endref

\ref \key RST \by N. Robertson, P. Seymor, R. Thomas
\paper Sachs' linkless embedding conjecture
\jour J. Combin. Theory Ser. B \vol 64 \yr 1995 \pages 185--227
\endref

\ref \key Ro \by R. H. Rosen
\paper Decomposing $3$-space into circles and points
\jour Proc. Amer. Math. Soc. \vol 11 \yr 1960 \pages 918--928
\endref

\ref \key RoS \by C. P. Rourke, B. J. Sanderson
\paper Homology stratifications and intersection homology
\inbook Proc. of the Kirbyfest \bookinfo Geometry and Topology Monographs
\vol 2 \yr 1999 \pages 99--116
\endref

\ref \key Sa1 \by K. S. Sarkaria
\paper Embedding and unknotting of some polyhedra
\jour Proc. Amer. Math. Soc. \vol 100 \yr 1987 \pages 201--203
\endref

\ref \key Sa2 \bysame
\paper A one-dimensional Whitney trick and Kuratowski's graph planarity
criterion
\jour Israel J. Math. \vol 73 \yr 1991 \pages 79--89
\endref

\ref \key Sch \by A. S. \v Svarc
\paper The genus of a fibre space
\jour Trudy Moskov. Mat. Ob\v s\v c \vol 10 \yr 1961 \pages 217--272
\moreref \vol 11 \yr 1962 \pages 99--126
\transl English transl.
\paper Amer. Math. Soc. Transl. \vol 55 \yr 1966 \pages 49--140
\endref

\ref \key Sh \by A. Shapiro
\paper Obstructions to the embedding of a complex in a Euclidean space,
I, The first obstruction
\jour Ann. Math. \vol 66 \yr 1957 \pages 256--269
\endref

\ref \key SV \by E. V. Shchepin, A. Yu. Volovikov
\paper Antipodes and embeddings
\jour Mat. Sbornik \vol 196 \issue 1 \yr 2005 \pages 3--32
\transl English transl. \jour Sb. Math. \vol 196 \yr 2005 \pages 1--28
\endref

\ref \key ST \by R. Shinjo, K. Taniyama
\paper Homology classification of spatial graphs by linking numbers
and Simon invariants
\jour Topol. Appl. \vol 134 \yr 2003 \pages 53--67
\endref

\ref \key S1 \by A. B. Skopenkov
\paper Embedding and knotting of manifolds in Euclidean spaces
\inbook Surveys in Contemporary Mathematics \eds N. Young, Y. Choi
\bookinfo London Math. Soc. Lect. Note Ser. \vol 347 \yr 2008 \pages 248--342
\moreref arXiv:math/0604045
\endref

\ref \key S2 \bysame
\paper A new invariant and parametric connected sum of embeddings
\jour Fund. Math. \vol 197 \yr 2007 \pages 253--269
\moreref arXiv:math.GT/0509621
\endref

\ref \key Sk1 \by M. B. Skopenkov
\paper Embedding products of graphs into Euclidean space
\jour Fund. Math. \vol 179 \yr 2003 \pages 191--198
\moreref arXiv:0808.1199
\endref

\ref \key Sk2 \bysame
\paper On approximability by embeddings of cycles in the plane
\jour Topol. Appl. \vol 134 \yr 2003 \pages 1--22
\moreref arXiv:0808.1187
\endref

\ref \key St \by S. Stolz
\paper The level of real projective spaces
\jour Comm. Math. Helv. \vol 64 \yr 1989 \pages 661--674
\endref

\ref \key Ta1 \by K. Taniyama
\paper Homology classification of spatial embeddings of a graph
\jour Topol. Appl. \vol 65 \yr 1995 \pages 205--228
\endref

\ref \key Ta2 \bysame
\paper Higher dimensional links in a simplicial complex embedded in a sphere
\jour Pacific J. Math. \vol 194 \yr 2000 \pages 465--467
\endref

\ref \key Va1 \by V. A. Vassiliev
\book Topology of Complements to Discriminants
\publ Phasis \publaddr Moscow \yr 1997
\transl Partial English transl.
\book Complements of Discriminants of Smooth Maps: Topology and Applications:
Revised Edition \publ Amer. Math. Soc. \yr 1994
\endref

\ref \key Va2 \bysame
\paper First-order invariants and cohomology of spaces of embeddings of
self-intersect\-ing curves in $\R^n$
\jour Izv. RAN, Ser. Mat. \vol 69 \issue 5 \yr 2005 \pages 3--52
\transl English transl.
\jour Izv. Math. \vol 69 \yr 2005 \pages 865--912
\endref

\ref \key Um \by B. R. Ummel
\paper Some examples relating the deleted product to embeddability
\jour Proc. Amer. Math Soc. \vol 31 \yr 1972 \pages 307--311
\endref

\ref \key vK \by E. R. van Kampen
\paper Komplexe in euklidischen R\"aumen
\jour Abh. Math. Sem. Univ. Hamburg \vol 9 \yr 1932 \pages 72--78 and
152--153
\endref

\ref \key We \by C. Weber
\paper Plongements de poly\`edres dans le domain metastable
\jour Comment. Math. Helv. \vol 42 \yr 1967 \pages 1--27
\endref

\ref \key Wu' \by Wu, Tsen-Teh
\paper On the $\bmod 2$ imbedding classes of triangulable compact manifold
\jour Science Record (New Ser.) \vol 2 \issue 3 \yr 1958 \pages 435--438
\endref

\ref \key Wu \by Wu, Wen-Ts\"un (Wu Wen-Jun, Wu Wen-Ch\"un)
\paper On the realization of complexes in a euclidean space. I, II, III
\jour Acta Math. Sinica \vol 5 \yr 1955 \pages 505--552
\moreref \vol 7 \yr 1957 \pages 79--101
\moreref \vol 8 \yr 1958 \pages 79--94
\transl English transl.
\jour Sci. Sinica \vol 7 \yr 1958 \pages 251--297 and 365--387
\moreref \vol 8 \yr 1959 \pages 133--150
\moreref parts I \& III \inbook Selected works of Wen-Tsun Wu
\publ World Sci. \yr 2008 \pages 23--69, 71--83
\endref

\ref\key Ya \by I. V. Yaschenko
\paper Embedding a smooth compact manifold into $\R^n$
\jour Problems from Topology Atlas, Topology Atlas Document \# qaaa-04,
http://at.yorku.ca/q/a/a/a/04.htm \yr 1996
\endref

\ref \key Z \by E. C. Zeeman
\paper Polyhedral $n$-manifolds: II. Embeddings
\inbook Topology of $3$-manifolds and Re\-lated Topics \ed M. K. Fort, Jr.
\publ Prentice-Hall \publaddr Englewood Cliffs, NJ \yr 1962 \pages 64--70
\endref

\endRefs
\enddocument
\end